\documentclass[oprs, nonblindrev]{informs3-albert}
\OneAndAHalfSpacedXI 

\usepackage{natbib}
\bibpunct[, ]{(}{)}{,}{a}{}{,}%
\def\bibfont{\small}%
\usepackage{booktabs}

\usepackage{empheq}

\usepackage{listings}
\usepackage{color}
\usepackage{bold-extra}

\TheoremsNumberedThrough     
\EquationsNumberedThrough    

\definecolor{light-gray}{gray}{0.95}
\usepackage{amsmath}
\usepackage{blkarray}
\usepackage{amsfonts}
\usepackage{mathtools}
\usepackage{mathrsfs}
\usepackage{soul}
\setstcolor{blue}

\usepackage[bb=dsserif]{mathalpha}
\usepackage{bm}

\usepackage[normalem]{ulem}
\usepackage[dvipsnames]{xcolor}
\usepackage{tikz}
\usetikzlibrary{arrows.meta, backgrounds}
\usepackage{pgfplots}
\usepackage{graphicx}
\usepackage[singlelinecheck=false, font=scriptsize, labelfont=scriptsize]{caption}
\usepackage{subcaption}
\usepackage{comment}
\usepackage{algorithm}
\usepackage{algpseudocode}
\algnewcommand\algorithmicforeach{\textbf{for each}}
\algdef{S}[FOR]{ForEach}[1]{\algorithmicforeach\ #1\ \algorithmicdo}
\algdef{S}[FOR]{ForAll}[1]{\algorithmicforeall\ #1\ \algorithmicdo}

\begin{document}
\RUNAUTHOR{Hasturk}
\RUNTITLE{CRL for DIRP}
\TITLE{Constrained Reinforcement Learning for the Dynamic Inventory Routing Problem under Stochastic Supply and Demand}
\ARTICLEAUTHORS{%
	\AUTHOR{Umur Hasturk, Kees Jan Roodbergen, Evrim Ursavas}
	\AFF{Department of Operations, Faculty of Economics and Business, University of Groningen, Groningen, the Netherlands\\ \EMAIL{u.hasturk@rug.nl, k.j.roodbergen@rug.nl, e.ursavas@rug.nl}}
 	\AUTHOR{Albert H. Schrotenboer}
	\AFF{Operations, Planning, Accounting and Control Group, School of Industrial Engineering, Eindhoven University of Technology, Eindhoven, the Netherlands\\ \EMAIL{a.h.schrotenboer@tue.nl}}
}
\ABSTRACT{Green hydrogen has multiple use cases and is produced from renewable energy, such as solar or wind energy. It can be stored in large quantities, decoupling renewable energy generation from its use, and is therefore considered essential for achieving a climate-neutral economy. The intermittency of renewable energy generation and the stochastic nature of demand are, however, challenging factors for the dynamic planning of hydrogen storage and transportation. This holds particularly in the early-adoption phase when hydrogen distribution occurs through vehicle-based networks. 
We therefore address the Dynamic Inventory Routing Problem (DIRP) under stochastic supply and demand with direct deliveries for the vehicle-based distribution of hydrogen. To solve this problem, we propose a Constrained Reinforcement Learning (CRL) framework that integrates constraints into the learning process and incorporates parameterized post-decision state value predictions. Additionally, we introduce Lookahead-based CRL (LCRL), which improves decision-making over a multi-period horizon to enhance short-term planning while maintaining the value predictions. Our computational experiments demonstrate the efficacy of CRL and LCRL across diverse instances. Our learning methods provide near-optimal solutions on small scale instances that are solved via value iteration. Furthermore, both methods outperform typical deep learning approaches such as Proximal Policy Optimization, as well as classical inventory heuristics, such as $(s,S)$-policy-based and Power-of-Two-based heuristics. Furthermore, LCRL achieves a $10\%$ improvement over CRL on average, albeit with higher computational requirements. Analyses of optimal replenishment policies reveal that accounting for stochastic supply and demand influences these policies, showing the importance of our addition to the DIRP. 
}%
\maketitle

%
%
%
%
%

\section{Introduction}
Achieving carbon neutrality by 2050, as targeted by the European Union and other global initiatives \citep{ipcc2023}, requires a transition to renewable energy sources such as wind and solar \citep{holechek2022global}. While these sources play an essential role in decarbonizing energy systems, their output is inherently intermittent and uncertain \citep{drucke2021climatological}, creating significant challenges in maintaining the balance between supply and demand. Conventional gas-fired power plants are often used to address this imbalance, as they can be quickly activated to produce electricity when supply is low \citep{safari2019natural}. However, reliance on such fossil fuel-based solutions contradicts the carbon-neutral goals and pinpoints the need for alternative mechanisms. Hydrogen, particularly \textit{green hydrogen} produced through renewables, provides a promising solution to these challenges \citep{abe2019hydrogen}. It enables the storage of excess energy as gas during periods of high renewable generation via electrolysis, and its later conversion back to electricity when needed \citep{oliveira2021green}. Beyond energy storage, hydrogen can be transported from one location to another as an energy carrier and used directly in applications such as transportation, residential heating, and industrial processes \citep{ball2009future}. Its versatility and ability to integrate with multiple sectors make hydrogen a critical component in the development of sustainable energy systems. 

A hydrogen distribution network typically consists of a supplier and multiple demand locations. The supplier generates hydrogen from renewable energy sources, resulting in a stochastic green hydrogen supply. Demand locations, such as hydrogen refueling stations, residential areas, and industrial plants, have stochastic demand, requiring periodic replenishments to maintain operations. To meet these demands, hydrogen is transported by vehicles from the supplier to these geographically dispersed locations. In the near term, distribution is vehicle-based due to small volumes; pipelines may be added later, while last-mile deliveries remain vehicle-based as in conventional fuel logistics. For instance, the \textit{Hydrogen Energy Applications in Valley Environments for Northern Netherlands} (HEAVENN) project exemplifies this by incorporating both vehicle-based distribution and planning the future pipeline development \citep{newenergycoalition2020, heavenncite}. 
The supplier and customers can create inventories to buffer between supply and demand. This results in a complex stochastic dynamic planning problem that requires balancing between transportation, storage, and operational costs.

In this study, we conceptualize vehicle-based hydrogen distribution by introducing a dynamic inventory routing problem (DIRP) that considers stochastic supply and demand.
The supplier uses a homogeneous vehicle fleet to replenish the customers. Both the supplier and the customers have limited inventory capacity for hydrogen storage. The replenishments are, thus, subject to constraints, such as vehicle fleet size and customer inventory capacities. 
As is common in the inventory routing literature, we assume a direct delivery setting in our problem \citep{kleywegt2002stochastic, bertazzi2008analysis, archibald2009indexability, coelho2012inventory, bertazzi2013stochastic, bertazzi2015managing}. Moreover, this is highly practical considering the relatively small quantities that can be transported in a full truckload. 
If needed, multiple vehicles may be dispatched to a customer at the same time, to accommodate larger replenishments. 
When demand cannot be met, it is considered as a lost sale. In contrast, the supplier can sell excess hydrogen to other local suppliers. 
This formulation is intended as a focused operational abstraction of early-stage vehicle-based hydrogen distribution, rather than a complete representation of all hydrogen logistics. It builds on a setting in which distribution is organized around a producer and multiple demand locations \citep{staffell2019role, hasturk2024stochastic}. In such early-adoption settings, lower volumes and limited dedicated infrastructure make customer-specific vehicle deliveries a reasonable first-order planning setting; the direct-delivery and homogeneous-fleet assumptions then keep the analysis focused on the core inventory-routing trade-off under stochastic supply and demand. This problem is modeled as a Markov decision process (MDP) over an infinite time horizon, with each period representing a day. The goal of the problem is to determine an optimal, state-dependent, replenishment policy that uses the current state information, i.e. inventory levels on locations, to minimize future expected discounted costs, consisting of transportation, holding, lost sales, and profits from sales.

Finding optimal state-dependent replenishment policies is difficult due to the complexity of the action space.
First, the action space grows exponentially with the number of customers and the maximum size of the replenishment quantities, making explicit enumeration or probabilistic sampling intractable. Second, the feasible actions are not explicitly listed but are instead defined implicitly by constraints within a Mixed-Integer Programming (MIP) formulation.
Third, since the feasibility of the MIP depends on the current inventory levels, these MIP formulations change at each state, requiring redefinition at every decision.
While achieving optimality in such complex settings is unrealistic, various learning techniques show good performance in the literature \citep{vanvuchelen2020use, ulmer2022dynamic}.
Traditional reinforcement learning techniques, such as TD($\lambda$) and Proximal Policy Optimization (PPO), often struggle with such combinatorial, state-dependent action spaces \citep{greif2024combinatorial, harsha2025deep}. They either rely on assigning a probability to all feasible actions, which is intractable for DIRP due to the enormous action space, or rely on approximations like continuous action representations that are often not stable \citep{dehaybe2024deep}. Neural networks often violate the constraints, requiring post-hoc adjustments like allocation rules to project infeasible solutions onto a locally feasible solution \citep{stranieri2024performance}. 
Consistent with this, we experimented with several DRL baselines (PPO, REINFORCE, TD($\lambda$), A3C) and observed unstable or infeasible behavior on our DIRP instances; methods that require enumerated action sets (e.g., DQN/DDQN) are inapplicable, and projection/repair steps further degrade learning. These issues motivate methods that enforce feasibility by construction rather than repairing infeasible neural policies.

For this purpose, we develop a Constrained Reinforcement Learning (CRL) method tailored to the DIRP. 
Our approach belongs to a class of methods that incorporate state-dependent action feasibility by solving a Mixed-Integer Program (MIP) at each decision stage \citep{rivera2017anticipatory, jia2025scenario}. In our case, the feasible MIP region is solved with a quadratic objective that minimizes the sum of the immediate cost of the selected action, and the predicted future cost associated with the resulting post-decision state. The separation of immediate costs from the value function is common in  Approximate Dynamic Programming (ADP), as the immediate costs (e.g., transportation) are often deterministic and can be precisely calculated separately from the stochastic future costs \citep{powell2019unified}. Our CRL adopts this structure to focus the learning process on the future costs, reducing the unnecessary overhead associated with learning the immediate action cost.
These future predictions are learned via a nonlinear value function approximation of post-decision state values. Additionally, we introduce an extended version of this algorithm, Lookahead-based Constrained Reinforcement Learning (LCRL), which incorporates lookahead-based heuristics to improve decision quality with the learned parameters, a strategy that has proven effective in similar dynamic routing contexts \citep{al2020approximate, ulmer2020horizontal}. Although LCRL increases computational requirements compared to CRL, it achieves higher solution quality, making it effective in applications demanding high-accuracy decisions. Both CRL and LCRL provide robust alternatives to the DIRP, offering a trade-off between solution quality and computational time. 

This study makes several key contributions to the related literature. To the best of our knowledge, we are the first to address the DIRP under stochastic supply. 
To address this, we propose CRL, which handles the problem's combinatorial complexity by sequentially searching within an MIP with an objective function incorporating a nonlinear value function approximation (VFA) of the post-decision state. Our LCRL extension further enhances solution quality by integrating a lookahead model.
We show that learning-based methods produce near-optimal solutions. 
Also, compared to similar inventory studies focusing on RL techniques, which typically focus on systems with three to five customers \citep{kleywegt2002stochastic, kaynov2023deep}, our approach scales to problems involving up to 15 customers and six vehicles, demonstrating the applicability of our methods to more complex, realistic settings. 
Other studies with realistic-sized instances rely on significant simplifying assumptions, such as assuming a fully deterministic system \citep{bertazzi2008analysis}, fixed order-up-to levels \citep{bertazzi2013stochastic, ortega2023stochastic}, or the use of only a single vehicle \citep{berman2001deliveries, coelho2014heuristics}. 
In these realistic settings, our methods surpass the Deep Reinforcement Learning (DRL) techniques, such as PPO, and outperform inventory management heuristics inspired by classical constructs such as $(s,S)$ policies. With these solutions, we further provide insights into the optimal policy structure and show how optimal actions are highly interdependent on supply and the other customers' inventory quantities, an aspect that is vastly ignored by the relevant literature.

The remainder of this paper is organized as follows. In Section \ref{lit}, we review the relevant literature on DIRP, reinforcement learning, and other related literature. In Section \ref{probdef}, we present the problem narrative and define the MDP underlying to the DIRP. In Section \ref{CRL}, we introduce our solution method CRL, as well as its extension LCRL. In Section \ref{bench}, we define two benchmarks based on inventory literature; $(s,S)$-policy Power-of-two based heuristics. In Section \ref{sect:comp}, the computational experiments are presented, discussing optimal policy structure and comparisons of our methods with benchmarks.
Finally, in Section \ref{conclusion}, we conclude the paper and discuss directions for future research.

\section{Literature Review} \label{lit}

The literature review is structured as follows. We first discuss the Dynamic Inventory Routing Problem in Section \ref{sect:DIRP}. In Section \ref{sect:CMDP}, we focus on Reinforcement Learning (RL) methods on similar constrained systems, followed by a discussion on the closely related field of Approximate Dynamic Programming (ADP) in Section \ref{sect:ADP}.  Finally, Section \ref{sect:OWMR} briefly examines One Warehouse Multi-Retailer systems and their relevance to our work.

\subsection{Dynamic Inventory Routing Problem} \label{sect:DIRP}

The Dynamic Inventory Routing Problem (DIRP) is the problem that integrates inventory management with vehicle routing decisions in a dynamic setting. A single supplier must deliver products to a set of geographically dispersed customers over multiple periods. The goal is to decide when to deliver and in what quantities, ensuring that customer demands are met while minimizing costs. These costs include inventory holding costs at the customers, penalties for unmet demand (backorders and/or lost sales), and transportation costs incurred by sending vehicles from the supplier to the customers. See \cite{coelho2014thirty} for a review.

Most of the literature focuses on this problem with the assumption that demand is uncertain. 
\cite{kleywegt2002stochastic} were pioneers in addressing this DIRP by employing decomposition methods to approximate the value functions of the Markov decision process (MDP), focusing on individual customer subproblems. Their work used simulation to refine policy decisions by assessing probabilities of customer visits, resembling modern reinforcement learning techniques. They assume a direct delivery setting, where later extending their work with routes up to three customers in \cite{kleywegt2004dynamic}. Similarly, \cite{adelman2004price} approached the problem from a multi-product setting for vendor-managed inventory at the supermarket chain.
They employ column generation to develop daily routes by solving a set-packing problem where subproblems resemble nonlinear discrete knapsack problems. \cite{bard1998decomposition} explore decomposition based heuristic techniques, using satellite facilities for emergency shipments. The method categorizes customers into two sets: those who must be replenished and those who might only be replenished if cost savings are evident.

In most cases, direct deliveries are applied in DIRP as the volumes to be transported are large compared to the vehicle's inventory capacity, making it practical to dedicate a full truckload to a single customer.
For instance, \cite{bertazzi2008analysis} examine a scenario in which both supply and demand are deterministic, while differing from the majority of the literature by assuming supply is limited. Authors employ parametric optimization for each customer to determine the optimal single-frequency direct delivery policy. Similarly, \cite{archibald2009indexability} utilizes a parametric approach to devise optimal direct delivery schedules, accommodating scenarios where vehicles can complete multiple routes within a shift, depending on the distances covered. 

Another line of research focuses on simplifying the DIRP through Mixed-Integer Programming (MIP)-based approximations, solving a static model to approximate underlying dynamic decision. 
For example, \cite{hvattum2009scenario} adopt a finite-horizon approach, replacing infinite-horizon MDP with a finite MIP model. Concurrently, \cite{bertazzi2013stochastic} employ a deterministic mixed-integer linear programming (MILP) approach, and a branch-and-cut algorithm enhanced with valid inequalities for the direct delivery setting. Authors further simplify the setting by assuming each customer maintains a base-stock level up to their inventory capacity. 

Several DIRP studies assume cases with fleets having only one vehicle. \cite{berman2001deliveries} dynamically adjust delivery quantities based on real-time assessments of customer tank levels within the context of industrial gas distribution. This approach aims to minimize the total expected costs including penalties for early or late deliveries. Similarly, \cite{schwarz2006interactions} adopt a cyclic approach where deliveries are scheduled at a constant frequency, selected the same for all customers. \cite{cui2023inventory} explore a case with a single uncapacitated vehicle and backlogged demand over a finite time horizon. Additionally, \cite{coelho2014heuristics} propose heuristics for large-scale instances with one vehicle, focusing on its capacity utilization and demand consistency by dispatching only when the vehicle is sufficiently filled and delivering only to customers within a selected inventory range, further restricting the action space.

Another simplification technique is done via exploring order-up-to (OU) policies within the DIRP, where customer inventories are always replenished up-to some selected quantities. \cite{bertazzi2015managing} address a variant by ensuring customer inventory levels meet a predefined maximum at each replenishment. Their solution employs a matheuristic combining rollout algorithms with MIP. Similarly, \cite{ortega2023stochastic} solve a case where customers are replenished up to a given target service level using a time-windows setting via chance-constrained programming techniques. They adopt an adaptive large neighborhood search (ALNS) based matheuristic for routing, and implement policy learning for the replenishment decisions. 

Beyond these primary DIRP formulations, various studies have explored alternative problem extensions and methodological adaptations. For instance, some works address specific application settings, such as perishable products \citep{crama2018stochastic} and multi-product contexts \citep{huang2010modified}. Others extend DIRP to more complex problem settings, such as a bike-sharing system where products are temporarily rented between stations \citep{brinkmann2019dynamic}. This study considers dual demand types ---renting and returning--- without a supplier, employing a dynamic lookahead policy to maintain predetermined service levels. Additionally, \cite{ozener2013allocating} addresses an inventory routing game, integrating IRP with cooperative game theory. The stochastic cyclic inventory routing problem (SCIRP) has also been studied under an infinite horizon, converting the underlying dynamic settings into cyclic delivery schedules \citep{malicki2021cyclic,sonntag2023stochastic,raa2023shortfall}.

Our research aims to address a significant gap in the DIRP literature by incorporating supply uncertainty. 
It is not yet explored in DIRP to the best of our knowledge, despite its recent attention in broader IRP studies \citep{jafarkhan2018efficient, alvarez2021inventory, hasturk2024stochastic}. Integrating stochastic supply transforms the underlying DIRP from a weakly coupled Markov decision process to a strongly coupled system. This increases the problem complexity significantly, as the scarcity of supply complicates the application of traditional decomposition approaches, which is adopted in most of the direct delivery literature. 
Furthermore, we show that the on-hand supply quantity changes the optimal policy structure significantly, underscoring its importance in the problem scope.
Overall, existing DIRP studies ignore stochastic supply, and simplify the decision problem through decomposition, static approximations, single-vehicle settings, or restricted replenishment policies. This motivates our study of an infinite-horizon DIRP with stochastic supply and demand, without relying on such simplifications, which improves its practical application in addressing modern logistical challenges, such as those found in green-hydrogen distribution.

\subsection{Reinforcement Learning} \label{sect:CMDP}

Reinforcement learning (RL) has increasingly been applied to solve the complexities of Dynamic Inventory Routing Problems (DIRP). 
\cite{achamrah2022solving} integrate genetic algorithms with Deep Q-learning (DQ) to facilitate routing decisions in their model, which introduces substitutions and transshipments as strategies to mitigate shortages at the customer level. 
\cite{greif2024combinatorial} explore DIRP with a single vehicle with only full vehicle capacity replenishments, limiting to no more than one route per period. This study integrates historical demand data into the state representation and adopts a prize-collecting mechanism that selects customers to maximize revenue.
Moreover, \cite{guo2024efficient} apply RL within a combinatorial framework to optimize a bike-sharing system, including additional performance indicators such as user satisfaction. Their approach includes heuristic-based dispatching and demand forecasting through Gradient Boosting Machines. Recent DRL studies also address VRP variants through route-construction methods based on RL techniques such as REINFORCE and Proximal Policy Optimization (PPO)  \citep{guan2025token, guan2026heterogeneous}, although their focus remains on routing rather than the inventory-control coupling central to DIRP.

RL has also seen applications in inventory control problems, which share cost structure and settings similar to the direct delivery DIRP model. For instance, \cite{temizoz2025deep} employ approximate policy evaluation with deep learning to generate multiple scenario trees for state evaluation. 
Unlike conventional approaches that typically compare the prediction against one observation, this method enhances the learning accuracy of the underlying stochasticity. 
Other traditional RL algorithms are also tested on inventory problems, such as PPO \citep{vanvuchelen2020use}, Asynchronous Advantage Actor-Critic (A3C) \citep{gijsbrechts2022can}, and DQ \citep{oroojlooyjadid2022deep}. \cite{vanvuchelen2020use} apply PPO onto a joint replenishment problem, but with a highly constrained action space: vehicles are dispatched only when fully loaded, and the problem setting is limited to just two customers. 
\cite{gijsbrechts2022can} explore A3C in inventory management across several settings including lost sales, dual-sourcing, and multi-echelon systems. \cite{oroojlooyjadid2022deep} apply DQ to the multi-agent, cooperative beer game, emphasizing the use of transfer learning to enhance learning efficacy. See also \cite{boute2022deep} for a general roadmap of Deep Reinforcement Learning on inventory control.

Our CRL algorithms select actions from given MIP regions, ensuring feasibility throughout the learning process. This is relevant in constrained Markov decision processes (CMDP) \citep{altman2021constrained}. However, the CMDPs are subject to constraints on the long-run average. Our problem instead involves state-dependent constraints, defining a different action set at each period. The CMDP
literature field encompasses various methodologies under different terminologies. 
For example, \cite{yu2022reachability} introduce ``Reachability Constrained Reinforcement Learning," where safety constraints are embedded via a safety value function that predicts the worst-case constraint violations. \cite{chow2018risk} discuss ``Chance-constrained Reinforcement Learning," which allows for a probabilistic level of risk in constraint violations, differing from hard constrained models. \cite{miryoosefi2019reinforcement} and \cite{li2021augmented} transform constrained problems into unconstrained min-max formulations through Lagrangian multipliers. 

In the transportation and logistics context, CMDP applications are scarce. Existing approaches address constraints through chance-constrained formulations, for example, to manage service levels in dynamic delivery routing \citep{ulmer2022dynamic}, or by integrating them as penalties within the reward function, as seen in vehicle scheduling \citep{zhang2023hierarchical}. The broader challenge of designing policies for problems with hard, state-dependent constraints is discussed by \cite{powell2022designing}, who highlights the need for more robust methods. This motivates our work, which aims to address such hard constraints directly.

Our study contributes to constrained RL in three ways. First, CRL selects actions from a vast, non-enumerable set which itself is a state-dependent MIP region. Second, we learn post-decision values with a nonlinear value function approximation (VFA), enabling efficient learning without exhaustive state space search. Third, feasibility is enforced at every decision via the mathematical program, rather than through Lagrangian penalties that trade off violations against costs; such penalties can obscure whether increases in the objective stem from system costs or constraint violations and can destabilize learning. Together, these elements enable feasible policy learning for DIRP with stochastic supply. Overall, existing RL approaches for DIRP and related inventory settings typically assume simplified or enumerable action spaces, or address constraints through penalties or heuristic adjustments. These limitations motivate our CRL approach, which combines state-dependent mathematical programming and VFA to obtain better policies in a large constrained action space.

\subsection{Approximate Dynamic Programming} \label{sect:ADP}

Approximate Dynamic Programming (ADP) provides a framework of methods for solving large-scale Markov Decision Processes where exact solutions are computationally intractable \citep{powell2011approximate}. Most of the literature focuses on constraint driven MDPs as in our case. While the field is vast, this review focuses on a class of ADP algorithms particularly relevant to transportation and logistics. Moreover, we select from the literature in which an action is selected by minimizing the sum of the immediate action cost and an approximation of the post-decision state costs, the state immediately after a decision has been made but before new stochastic information has arrived.

A common algorithm used to learn these post-decision state values is Approximate Value Iteration (AVI), which iteratively learns the value function by simulating samples and updating value estimates based on observed outcomes. AVI has been effectively applied in diverse contexts, such as humanitarian relief distribution by \cite{van2023reinforcement}, which utilizes approximate policy iteration
, and in container relocation problems addressed by \cite{boschma2023approximate}, where aggregated states are used to manage the complexity of the decision space. 

A key component of ADP is the Value Function Approximation (VFA), which is used to approximate the post-decision state information as future costs. VFA architectures vary, ranging from non-parametric state aggregation via lookup tables \citep{ulmer2018budgeting}, parametric models using basis functions \citep{rivera2017anticipatory, van2019delivery}, hybrid structures \citep{ulmer2020meso}, and neural networks in Deep RL applications \citep{joe2020deep, chen2022deep}. Given a VFA, an action must be selected from the state-dependent feasible set. For problems with simpler action spaces, this is often done via simple enumerations or heuristics \citep{ulmer2019offline}. However, more complex, constrained action spaces may require solving a mathematical program at each step, such as a Linear Program \citep{al2020approximate} or a Mixed-Integer Program \citep{rivera2017anticipatory}. An alternative to relying on a VFA is to use a lookahead policy, which approximates the dynamic problem by solving a static model with a truncated, finite horizon approximation; only the action for the first period is implemented before the process is repeated at the next state \citep{soeffker2022stochastic}. The structure of minimizing the immediate action cost plus the approximation of the post-decision state value can be viewed as the simplest form, a 1-step lookahead, where the VFA serves as the terminal cost at the end of a single-period horizon. Some research extends this to n-step lookaheads, often termed rollout algorithms, where a multi-period simulation is performed over a limited horizon, and the VFA is used to approximate the terminal value of the state reached at the end of that horizon \citep{ulmer2019offline, ulmer2020horizontal}.

Our work contributes to this literature in several aspects. First, while many applications rely on state aggregation or linear approximations, our model employs a parametric, nonlinear (polynomial) VFA. This structure is more expressive than linear approximations and generalizes across the state space, avoiding the need for the exhaustive data often required by lookup-table methods. Second, we integrate this VFA directly into the action-selection problem, which necessarily formulates the decision at each stage as a Mixed-Integer Quadratic Program (MIQP). To the best of our knowledge, the quadratic program is a novel approach in the ADP context. Finally, we embed this MIQP-based decision process within an n-step lookahead framework, thereby combining the accuracy of a nonlinear value function with the foresight of multi-period planning. Overall, the ADP literature provides the main methodological basis for our setting through post-decision value approximation, mathematical programming-based action selection, and lookahead. We extend this line of work to the DIRP setting via a combination of these settings with a nonlinear VFA.

\subsection{One Warehouse Multi-Retailer Problem} \label{sect:OWMR}

The One Warehouse Multi-Retailer (OWMR) problem, also known as the two-echelon inventory problem, focuses on optimizing the distribution of inventory from a central warehouse to multiple retailers. This structure resembles DIRP with a direct delivery setting, where although lacking a transportation setting, the inventory cost structure is similarly aligned. \cite{chu2010power} investigate a power-of-two (PO2) based replenishment policy under stochastic demand, demonstrating that their algorithm performs within $\sqrt[3]{2} \approx 1.26$ times the optimal solution in the worst case. Extending this approach, \cite{li2019inventory} integrate emission costs into the objective, applying PO2 to also address environmental impacts. \cite{malmberg2023evaluation} treat emissions not as a cost component but as a constraint, aligning the replenishment strategy with specified emission targets. \cite{sakulsom2019heuristics} tackle the problem using a $(R,s,S)$ policy, with a two-stage heuristic that introduces safety stocks as opposed to the service level targets common in prior studies. Other variants of the OWMR problem focus on multi-item settings \citep{johansson2020controlling}, perishable items \citep{nguyen2014consolidation}, and enhancements of micro-retailer and consumer welfare \citep{gui2019improving}.

Applying RL has recently gained attention in the OWMR literature. \cite{kaynov2023deep} implement Proximal Policy Optimization (PPO) to manage limited supply scenarios with lead times, having an optimality gap of $3\%$ for lost sales, and ranging up to $20\%$ for partial lost sales. For the cases when the action suggested by the neural network exceeds available warehouse inventory, the authors propose two allocation methods: proportional and random. \cite{stranieri2024performance} later extend these allocation methods by introducing a ``balanced allocation rule" through PPO to further refine the strategy.

While the OWMR cost structure may also represent vehicle distribution costs in DIRP, we improve OWMR literature by integrating limited vehicle fleets and inventory capacities at each vehicle, which significantly enhances the complexity of the action space. Additionally, traditional studies further simplify the action space, resulting in less challenging scenarios. For example, \cite{kaynov2023deep} analyze various instances, with only two involving more than four customers and maximum expected demands of five units per customer. In contrast, our model includes higher inventory levels, more customers, and the option for products to be sold directly from the supplier's inventory. Furthermore, we create two benchmark solution methods based on traditional methods used in OWMR, PO2- and (s,S)-based policies, and show our algorithms improve their solution quality significantly. Overall, the OWMR literature provides relevant benchmarks and cost structures, but missing critical elements such as limited vehicles and inventories, and its action spaces are simpler than in our setting. This motivates developing CRL and LCRL for the more general problems such as DIRP.

\section{Problem Definition} \label{probdef}
In this study, we consider the dynamic inventory routing problem, where the supply chain consists of a single supplier and multiple customers. We define a location set $\mathscr{N}^0 \coloneqq \{ 0,1,\ldots,N \}$, where $0$ represents the supplier location and $i \in \mathscr{N} \coloneqq \{ 1,2,\ldots,N \}$ represents customer $i$. The supplier has storage and receives inventory through an exogenous, stochastic supply process. Each customer location has its own inventory and observes random demand. We consider an infinite discrete time horizon with periods (i.e., day). Both supply and demand realizations are modeled as random variables with known stationary distributions for each period $t \in \mathbb{T} \coloneqq \{0,1,2,\dots\}$. We assume the distributions are independent across periods and locations. The aim is to find a cost-minimizing replenishment strategy that balances inventory levels at all locations by sending inventory from the supplier to customers via vehicles in a direct delivery setting.

We model this problem as a Markov decision process. The state space is given by \( \mathscr{X} = [0, U_0] \times [0, U_1] \times \dots \times [0, U_N] \), where \( x \coloneqq \{ x_0, x_1, \dots , x_N \} \in \mathscr{X} \) represents the system's state. Here, for each location \( i \in \mathscr{N} ^0 \), \( x_i \) denotes the inventory at the beginning of a period, and  \( U_i \) is the inventory capacity. At each state $x$, the decision maker decides the quantity of inventory to replenish to each customer ($a_i$), as well as the quantity of inventory to sell directly from the supplier storage ($a_0$). These actions are subject to constraints, including the availability of inventory at the supplier, vehicles, and storage capacities. The supplier uses a homogeneous fleet of vehicles to transport inventory to customers. The fleet consists of $q$ vehicles, each having a capacity of $C$. We assume each vehicle can make at most one trip per period, and a customer may receive replenishments via multiple vehicles at a period. 
The action, thus, is defined as $a \coloneqq \{a_0, a_1,\ldots,a_N\} \in A(x) $, where $A(x)$ is the feasible action set in state $x$. For each customer $i$, the decision variable $b_i$ represents the number of vehicles dispatched to replenish $a_i$ units of inventory. The action set is defined as:

\begin{empheq}[left={A(x) \coloneqq \left\{ (a_0, a_1, \dots, a_N) \, \middle|\, \begin{aligned} } , right={\end{aligned} \right\}}]{align}
    &\sum_{i \in \mathscr{N}^0} a_i \leq x_0, \label{totAction} \\
    &x_i + a_i \leq U_i && \forall i \in \mathscr{N}, \label{capCust} \\
    &a_i \leq C b_i && \forall i \in \mathscr{N}, \label{capTruck} \\
    &\sum_{i \in \mathscr{N}} b_i \leq q, \label{limTruck} \\
    &a_i \geq 0 && \forall i \in \mathscr{N}^0, \label{nonNeg} \\
    &b_i \in \mathbb Z_{\geq 0} && \forall i \in \mathscr{N}. \label{MDP-end}
\end{empheq}

Constraint \eqref{totAction} ensures that the total inventory for replenishment and selling does not exceed the available inventory at the supplier location. Constraints \eqref{capCust} ensure that the customer inventory capacity is not exceeded after the replenishment. Constraints \eqref{capTruck} assign deliveries to the vehicles according to the vehicle capacity. Constraint \eqref{limTruck} ensures that the action does not use more vehicles than the fleet has.

The costs associated with the action, denoted as $c_x(a)$ include both transportation and sales. The transportation cost of customer $i$, $c_i^\textsc{t}$, consists of a fixed cost of $W$ per vehicle, $b_i$, and a variable cost of $w$ dependent on the distance $d_{i}$ between supplier and customer $i$; $c_i^\textsc{t}(b_i) \coloneqq b_i \left( W + 2 w d_{i} \right)$. Additionally, the supplier makes a profit of $\rho$ per each item sold externally, outside the replenishment system; $a_0$. The cost of the action is, then, defined as:
\begin{align}
    &c_x(a) \coloneqq -\rho a_0 + \sum_{i=1}^N c_i^\textsc{t}(b_i).
\end{align}

After executing the action, the system transitions to a post-decision state $s \coloneqq \{ s_0, s_1, \dots , s_N \}$, where $s_0 = x_0 - \sum_{i=0}^N a_i$ represents the remaining inventory at the supplier after selling and replenishments, and $s_i = x_i + a_i$ for each customer $i \in \mathscr{N}$ indicates the updated inventory post-delivery. At the post-decision state, the system observes the realization of stochasticity per period, denoted by $\phi \coloneqq \{ \phi_0, \phi_1, \dots , \phi_N \}$, drawn from known probability distributions that capture supply and demand uncertainties. The random variables $O \coloneqq \{ O_0, O_1, \dots , O_N \}$ represent the supply or demand at the locations. Stochasticity is modeled through discrete support sets $\Phi \coloneqq \{ \Phi_0, \Phi_1, \dots , \Phi_N \}$, where $\phi_i \in \Phi_i$ is a possible realization at location $i \in \mathscr{N}^0$. The probability of each realization $\phi_i$ is given by $\mathbb{P}(O_i = \phi_i)$.

Let $c_s(\phi)$ capture the costs associated with the stochasticity, including holding costs, lost sales due to insufficient inventory at customer locations, and excess inventory sales when the supplier's stock exceeds its capacity. At the end of a period, a holding cost of $h_s$ and $h_c$ are paid per inventory at the supplier and at the customer locations, respectively. A lost sale penalty of $\ell$ is applied for each demand that is not satisfied at customers. Finally, the supplier has to sell the product that does not fit into its inventory at the end of a period. We formally denote $c_s(\phi)$ as:
\begin{align}
    c_s(\phi) \coloneqq & -\rho \left( s_0 + \phi_0 - U_0\right)^+ + h_s \min \left\{ s_0 + \phi_0, U_0  \right\} \nonumber \\
    & +\sum_{i=1}^N \left( h_c\left( s_i - \phi_i \right)^+ + \ell \left( s_i - \phi_i \right)^- \right),
\end{align}
where the notation $(\cdot)^+$ and $(\cdot)^-$ represent the positive and negative parts, respectively. 

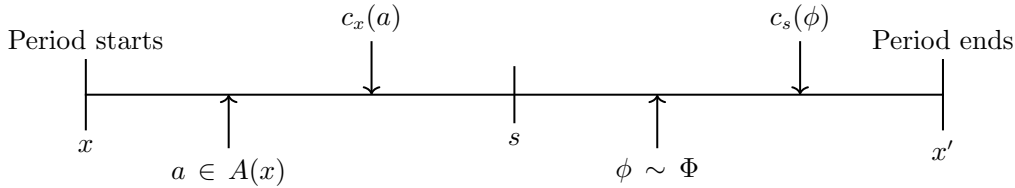
\begin{figure}[h]
\centering
\resizebox{\textwidth}{!}{%
\begin{tikzpicture}[x=2cm,nodes={text width=5cm, align=center, font=\normalsize}]
\draw[black, thick, >=latex, line cap=rect] 
  (0,0) -- (6,0); 

\foreach \Xc/\Txt in {0/Period starts} {
  \draw[black, thick] (\Xc, -0.5) -- (\Xc, 0.5); 
  \node[above] at (\Xc, 0.5) {\Txt};
  \node[below] at (\Xc, -0.5) {$x$};
}

\foreach \Xc/\Txt in {6/$x'$} {
  \draw[black, thick] (\Xc, 0.5) -- (\Xc, -0.5); 
  \node[below] at (\Xc, -0.5) {\Txt};
  \node[above] at (\Xc, 0.5) {Period ends};
}

\foreach \Xc/\Txt in {3/{$s$}} {
  \draw[black, thick] (\Xc, -0.4) -- (\Xc, 0.4); 
  \node[below] at (\Xc, -0.4) {\Txt};
  }

\foreach \Xc/\Txt in {1/{$a \in A(x)$}, 4/{$\phi \sim \Phi$ }} {
  \draw[black, thick, ->] (\Xc, -0.75) -- (\Xc, 0); 
  \node[below] at (\Xc, -0.75) {\Txt};
}

\foreach \Xc/\Txt in {2/{ $c_x(a)$}, 5/{$c_s(\phi)$}} {
  \draw[black, thick, ->] (\Xc, 0.75) -- (\Xc, 0); 
  \node[above] at (\Xc, 0.75) {\Txt};
}

\end{tikzpicture}
}
\caption{Timeline of a single period.}
\label{fig:timeline}
\end{figure}

Let $x' \coloneqq \{ x'_0, x'_1, \dots , x'_N \}$ denote the pre-decision state for the next period. Here, $x'_0 = \min \left\{ s_0 + \phi_0, U_0  \right\}$, and $x'_i = \left( s_i - \phi_i \right)^+$ for each customer $i$ represent the corresponding inventories after accounting for the stochasticity, which is the start inventory of the next period. The timeline of this transition is given in Figure \ref{fig:timeline}. 

\begin{figure}[h]
    \centering
    \begin{subfigure}[t]{0.28\textwidth}
        \centering
        \includegraphics[width=\textwidth]{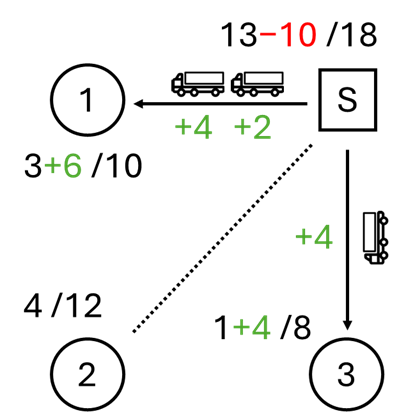}
        \caption{Action $a=\{0,6,0,4\}$: Customers 1 and 3 are replenished using vehicles with a capacity of 4 at $x=\{13,3,4,1\}$.}
        \label{fig:toy1}
    \end{subfigure}
    \hfill
    \begin{subfigure}[t]{0.29\textwidth}
        \centering
        \includegraphics[width=\textwidth]{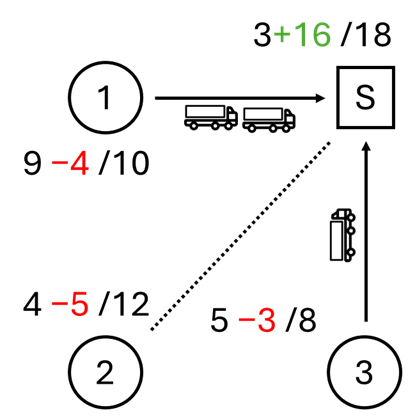}
        \caption{Stochasticity $\phi=\{16,4,5,3\}$: Supply and demands are realized for the supplier and customers at $s=\{3,9,4,5\}$.}
        \label{fig:toy2}
    \end{subfigure}
    \hfill
    \begin{subfigure}[t]{0.27\textwidth}
        \centering
        \includegraphics[width=\textwidth]{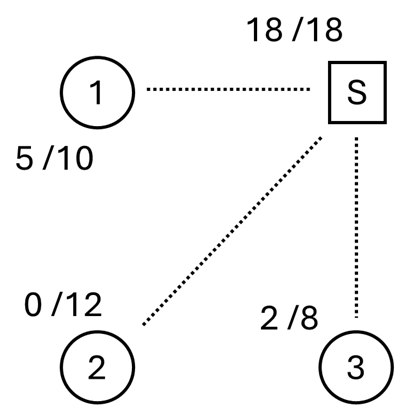}
        \caption{Holding and lost sales costs are incurred at the end of state. Excess supply is forced to be sold.}
        \label{fig:toy3}
    \end{subfigure}
    \caption{Illustrative example of DIRP at a single stage, showing transportation, inventory updates, stochastic realizations, and end-of-period inventories.}
    \label{fig:dirp-toy}
\end{figure}

Additionally, Figure~\ref{fig:dirp-toy} provides an illustrative one-period example of the DIRP. 
Figure~\ref{fig:toy1} displays the action $a = \{0,6,0,4\}$ and the associated vehicle assignments, which satisfy the capacity and fleet constraints \eqref{totAction}–\eqref{limTruck}. 
Figure~\ref{fig:toy2} shows the realization of supply and demand, corresponding to the stochasticity with $\phi=\{16,4,5,3\}$. 
Figure~\ref{fig:toy3} corresponds to the end-of-period cost calculation: the supplier realizes 16 units of supply and, due to the capacity $U_0 = 18$, is forced to sell one unit externally, while customer 2 incurs one unit of lost sales and the remaining inventories generate holding costs on $c_s(\phi)$.

Our goal is to minimize the expected long-run average cost per period, formally expressed as:
\begin{align}
    Z^* = \min_{a_t \in A_t(x)} \left( \lim_{T \to \infty} \frac{1}{T} \mathbb{E} \left[ \sum_{t=0}^{T-1} \left( c_x(a_t) + c_s(\phi_t) \right) \right] \right). \label{longrun}
\end{align}

Let $V_t(x)$ be the value function that represents the expected cumulative future cost of being at state $x$ at period $t$ under the optimal policy. The recursive Bellman equation representing this transition is defined as in Eq. \eqref{bellman}.
\begin{align}
 V_t(x)  = \min_{a_t \in A_t(x)} \left\{  c_x(a_t)  + \sum_{\phi_t \in \Phi_t} \mathbb{P}\left\{ O_t = \phi_t \right\} \left( c_s(\phi_t)+ V_{t+1} (x') \right) \right\}. \label{bellman}
\end{align}

\section{Constrained Reinforcement Learning Approaches} \label{CRL}

Traditional approaches to solving Markov decision processes (MDPs), such as value iteration, often become impractical in complex scenarios due to the curse of dimensionality. For instance, the state space size of our dynamic inventory routing problem under discrete actions is of order $\mathcal{O}(\prod_{i \in \mathscr{N}^0} U_i)$, and thus grows exponentially in the number of customers. This exponential growth means that even small-scale DIRPs with only four to five customers can result in a state space that is computationally intractable.

Reinforcement Learning (RL) is particularly capable at handling large-scale MDPs where traditional methods reach their boundaries. However, typical RL applications, such as those for Atari games or game of Go \citep{mnih2013playing,silver2016mastering}, involve enumerable action sets of a few hundred actions, enabling to estimate actions for each state effectively via techniques like discrete actor-critic or proximal policy optimization (PPO) \citep{schulman2017proximal}. In contrast, the action space encountered in our MDP is significantly larger than such traditional AI applications. Namely, we have about $100,000$ actions for a DIRP with three to four customers (assuming $a \in A(x)$ and $a_i \in \mathbb Z_{\ge 0}$). Moreover, these actions are not explicitly given but are instead defined by a constraint set of Eq. \eqref{totAction}-\eqref{MDP-end}, changing every state. These constraints complicate the identification of $A(x)$, as they define a different mixed-integer programming (MIP) region at each state \citep{jia2025scenario}. For similar RL applications, the common literature applies Lagrangian relaxation on the action space \citep{miryoosefi2019reinforcement,li2021augmented} as, for example, neural networks have no direct capability of respecting such constraints. However, this approach is suboptimal due to various reasons: 1) it does not enforce hard constraints, 2) it rewards non-binding constraints, such as underutilizing vehicle capacities in our setting, and 3) (high) Lagrangian penalties can cause divergence in the algorithm, as shown in our preliminary experiments.

To address these issues, we introduce Constrained Reinforcement Learning (CRL), a reinforcement learning algorithm specifically designed for MDPs with state-dependent constraints. This algorithm, operating online, learns the values of post-decision states, $V_t(s)$, while also enabling constrained action selection within the feasible mixed-integer space using solvers like Gurobi. Similar to Eq. \eqref{bellman}, which defines the Bellman equation for the pre-decision state $x$, we now formulate the Bellman function for the post-decision state $s$ as follows:
\begin{align}
& V_t(s) =\sum_{\phi_t \in \Phi_t} \mathbb{P}\left\{ O_t = \phi_t \right\} \left(  c_s(\phi_t) +  \min_{a_t \in A_t(x)} \left\{  c_x(a_t) + V_{t+1}(s')  \right\}  \right). \label{postDec}
\end{align}

As values of pre- and post-decision states, \( V(x) \) and \( V(s) \), there also exist traditionally adapted Q-values, \( Q(x,a) \), which represent the expected total cost when taking action \( a \) in state \( x \) and following the optimal policy thereafter. In our context, this is given by:  
\begin{align}  
Q(x,a) &=    c_x(a)  + \sum_{\phi \in \Phi} \mathbb{P}\left\{ O = \phi \right\} \left( c_s(\phi)+ V(x') \right). \label{q}  
\end{align}  

Note that \( t \) is omitted in Eq. \eqref{q} since, in the infinite-horizon setting, \( Q_t(x,a_t) \), \( V_t(x) \), and \( V_t(s) \) converge to \( Q(x,a) \), \( V(x) \), and \( V(s) \).  

Q-learning is effective in traditional RL problems because immediate action costs are often uncertain or challenging to estimate. However, as in most Approximate Dynamic Programming (ADP) literature, our costs are explicitly defined, such as \( c_x(a) \), including the selling profit and transportation costs. Instead of following the conventional approach of policy improvement through Q-learning, we leverage this structure by separating the Q-value into its easily calculable immediate costs of \( c_x(a) \), and predicting only \( V(s) \), given as:  
\begin{align}  
Q(x,a) &=  c_x(a)  + V(s).  
\end{align}  

This separated approach reduces unnecessary computational overhead in estimating current action costs, and helps the model to solely focus learning $V(s)$. Moreover, as demonstrated in previous research \citep{sun2022hardware}, learning post-decision state values is generally a more effective approach compared to learning of the action values.    

The proposed algorithm is detailed in Algorithm \ref{RLalgo}. It focuses on learning an approximation to $V(s) \approx \hat{v}_w(s) \coloneqq  w^\intercal \psi(s)$, using a differentiable value function parameterization with features, $\psi(s)$. The estimated value for the post-decision state is calculated as a weighted linear combination of these features. The algorithm iteratively refines the feature weights, $w$, to improve the action selection process. The learning step follows a TD($\lambda$)-based update, where the temporal-difference error, $\delta$, and eligibility trace, $z$, are used to update the value-function weights, $w$.

\begin{algorithm}
    \caption{Constrained Reinforcement Learning} \label{RLalgo}
    \begin{algorithmic}[1] 
        \State \textbf{Input:} A differentiable value function parameterization for post-decision state $\hat{v}_w (s)$.
        \State \textbf{Algorithm Parameters:} $\lambda \in [0,1], \alpha > 0$, $\epsilon \in [0,1]$.
        \State Initialize $\bar{c} \in \mathbb{R} \coloneqq 0$, $w \in \mathbb{R}^d \coloneqq \textbf{0}$, $z \in \mathbb{R}^d \coloneqq \textbf{0}$, and $s \coloneqq s_0$.
        \While {convergence is not satisfied}
            \State On $s$, observe $\phi$, $c_s(\phi)$, and $x(s,\phi)$ as the next pre-decision state.
            \If{rand $ < \epsilon$}
            \State Select $a \in A(x)$ randomly
            \Else
            \State $a \coloneqq \argmin_{a \in A(x)} \left( c_x(a) + \hat{v}_w(s'(x,a)) \right)$ \label{RL-MIP}
            \EndIf
            \State $\delta \coloneqq c_s(\phi) + c_x(a)  + \hat{v}_w (s') - \bar{c} - \hat{v}_w (s)$ \Comment{TD error}
            \State $\bar{c} \coloneqq \bar{c} + \alpha \delta$ \Comment{Average-cost estimate}
            \State $z \coloneqq \lambda z + \nabla \hat{v}_w(s)$ \Comment{Eligibility trace}
            \State $w \coloneqq w + \alpha \delta z$ \Comment{TD($\lambda$) weight update}
            \State $s \coloneqq s'$
        \EndWhile
    \end{algorithmic}
\end{algorithm}

In this algorithm, $\bar{c}$ serves as a baseline for the expected cost per period, providing a reference point for observed costs. The temporal difference, $\delta$, measures the difference between the actual observation, $c_s(\phi) + c_x(a)  + \hat{v}_w (s') $, and the prior prediction, $\bar{c} + \hat{v}_w (s)$. Here, $s'$ represents the next post-decision state. $\alpha$ is the learning rate, and $\lambda$ is the forgetting factor for the eligibility trace. An eligibility traces vector, $z$, 
assists in adjusting the feature weights, $w$. Furthermore, the action selection mechanism employs an $\epsilon$-greedy approach, balancing exploration and exploitation. In line \ref{RL-MIP}, the algorithm solves a problem with a MIP feasible region where the objective function includes the weighted sum of features. 
We selected features based on the inventory levels as $s_i$, $s_i^2$, $s_i^3$, and $\sqrt{s_i}$ for each $i \in \mathscr{N}^0$. This selection is based on the observation that an optimal or near-optimal policy generally involves maintaining some target inventory levels for each customer to maximize system efficiency. Deviations from these targets can lead to increased costs. Such a cost function may be approximated using a third-degree quadratic function, and $\sqrt{s_i}$ helps further for this prediction without significantly increasing the solution time.

The aim of the algorithm is to determine a fixed-set of weights, $w$. For any given state $x$, a policy can be derived by solving the optimization problem at line \ref{RL-MIP} using these fixed weights to achieve a near-optimal solution. In our computational results, the long-term performance of the derived policies are evaluated through simulation, as detailed in Section \ref{sect:comp}. 

While the CRL algorithm offers a powerful method for solving MDPs with state-dependent constraints, further refinements can improve decision-making in highly dynamic environments like DIRPs. One such refinement involves incorporating lookahead strategies to optimize short-term actions further and adjust current decisions accordingly. In the following subsection, we explore this lookahead-based CRL approach as an alternative or complement to the CRL framework.

\subsection{Lookahead-based CRL} \label{sect:LCRL}
In this section, we introduce lookahead-based constrained reinforcement learning (LCRL) as an enhancement to CRL. While CRL efficiently handles large-scale MDPs with state-dependent constraints, a lookahead strategy improves planning by optimizing decisions further over a longer finite horizon. 

Dynamic IRPs may be addressed using such lookahead algorithms \citep{brinkmann2019dynamic, cuellar2024adaptive}, which transform the infinite-horizon DIRP into a tractable finite-horizon problem. In these methods, the optimal action for the current period is derived by solving a static MIP model over, for example, $T=10$ periods, using a two-stage decision tree with a finite set of scenarios. Although the horizon is truncated, the action at $t=0$ serves as an approximation of the infinite-horizon solution. We tested this approach in our setting and found its performance to be significantly worse compared to CRL. 
On instances sized at $N=5$ and $q=3$, using a truncated horizon of $T=5$ and a two-stage decision tree with 20 scenarios, the average cost was $15\%$ higher than that of CRL, while the solution time for one state was nearly $100$ times longer. Due to this poor performance, we decided to exclude this method from further use, also as a benchmark.

A more sophisticated method is a hybrid approach, where the learned value function is combined with a multi-period lookahead model, a strategy that has proven effective in similar dynamic routing contexts \citep{ulmer2019offline, ulmer2020horizontal}. In Algorithm \ref{RLalgo}, line \ref{RL-MIP}, the action selected by CRL minimizes immediate costs while accounting for predicted future costs. This can be seen as a one-step lookahead approach, where immediate actions optimize costs for one period ahead ($c_x(a)$) while considering long-term state predictions ($\hat{v}_w (s)$). 
Extending this to a multi-period lookahead would possibly improve the immediate action by considering a broader horizon, while still implementing the learned future cost predictions guiding again the long-term planning.
Let $t \in \mathscr{T} = \{1,2,\dots, T\}$ represent the decision epoch in the lookahead horizon, where at $t=0$, CRL takes its first action as usual. 
Additionally, for each $t \in \mathscr{T}$, we select actions for each scenario within a finite scenario set. These actions optimize short-term decisions within the planning horizon $\mathscr{T}$. For $t > T$, representing the infinite horizon beyond $\mathscr{T}$, the model still uses the value function approximation $\hat{v}_w (s)$ to account for future costs, ensuring that the state reached at the end of $\mathscr{T}$ reflects its impact on long-term planning. Let $\omega \in \Omega$ be a scenario in the two-stage decision tree. As with the decision variables in CRL, let $a_{it\omega}$ and $b_{it\omega}$ represent the inventory delivered to customer $i$ at epoch $t$ under scenario $\omega$, and the corresponding number of vehicles used, respectively. Similarly, let $x_{it\omega}$ denote the pre-decision inventory at location $i$ for epoch $t$ under scenario $\omega$. For managing shortages, a decision variable $\eta_{it\omega}$ is introduced, representing lost sales for customer $i$ at epoch $t$ under scenario $\omega$. Lastly, $\phi_{it\omega}$ represents the random variables of supply and demand at location $i$ in epoch $t$ under scenario $\omega$, with $\phi_0$ representing supply and $\phi_i$ for $i \in \mathscr{N}$ representing demand. The LCRL model is defined as follows:
\begin{align}
    \min\;& -\rho a_0 + \sum_{i=1}^N c_i^\textsc{t}(b_i) + \frac{1}{|\Omega|}\sum_{t\in\mathscr{T}}\sum_{\omega\in\Omega}
    \Big( -\rho a_{0t\omega} + \sum_{i=1}^N c_i^\textsc{t}(b_{it\omega}) \Big) \nonumber\\
    & + \frac{1}{|\Omega|}\sum_{t\in\mathscr{T}}\sum_{\omega\in\Omega}
    \Big( h_s x_{0t\omega} + \sum_{i=1}^N (h_c x_{it\omega} + \ell \eta_{it\omega}) \Big) + \frac{1}{|\Omega|}\sum_{\omega\in\Omega} \hat{v}_w(s'_\omega), \label{LMIP-start} \\
    \text{s.t. }  &\sum_{i \in \mathscr{N}^0} a_i \leq x_0,  \label{old-start}\\
    &x_i + a_i \leq U_i, && \forall i \in \mathscr{N},  \\
    &a_i \leq C b_i, && \forall i \in \mathscr{N},  \\
    &\sum_{i \in \mathscr{N}} b_i \leq q,  \label{old-end}\\
    & x_{it\omega} + a_{it\omega}  \leq U_i, && \forall i \in \mathscr{N}, t \in \mathscr{T}, \omega \in \Omega , \label{new-start}\\
    &a_{it\omega} \leq C b_{it\omega}, && \forall i \in \mathscr{N}, t \in \mathscr{T}, \omega \in \Omega,  \\
    &\sum_{i \in \mathscr{N}} b_{it\omega} \leq q, &&  t \in \mathscr{T}, \omega \in \Omega, \label{new-end}\\
    & x_{0t\omega} = x_0 + \phi_{0t\omega} - \sum_{i=0}^N a_i, && \forall t = 1, \omega \in \Omega, \label{inv-start}\\
    & x_{0t\omega} = x_{0,t-1,\omega} + \phi_{0t\omega} - \sum_{i=0}^N a_{i,t-1,\omega}, && \forall t \in \mathscr{T} \setminus \{1\}, \omega \in \Omega,\\
        & x_{it\omega} = x_i - \phi_{it\omega} + a_i + \eta_{it\omega}, && \forall i \in \mathscr{N} , t = 1 , \omega \in \Omega,\\
    & x_{it\omega} = x_{i,t-1,\omega} - \phi_{it\omega} + a_{i,t-1,\omega} + \eta_{it\omega}, && \forall i \in \mathscr{N}  , t \in \mathscr{T} \setminus \{1\}, \omega \in \Omega, \label{inv-end}\\
    & s'_{0\omega} = x_{0,T,\omega} - \sum_{i=0}^N a_{i,T,\omega}, && \forall   \omega \in \Omega, \label{sinv-start}\\
    & s'_{i\omega} = x_{i,T,\omega}  + a_{i,T,\omega}, && \forall i \in \mathscr{N} ,  \omega \in \Omega,\label{sinv-end}\\
    & a_i \geq 0, && \forall i \in \mathscr{N}^0,   \\
    & a_{it\omega} \geq 0, && \forall i \in \mathscr{N}^0, t \in \mathscr{T} ,  \omega \in \Omega , \\
    & b_i \in \mathbb Z_{\ge 0}, && \forall i \in \mathscr{N}, \\ 
    & b_{it\omega} \in \mathbb Z_{\ge 0}, && \forall i \in \mathscr{N}^0, t \in \mathscr{T} ,  \omega \in \Omega, \\ 
    & x_{it\omega} \geq 0, && \forall i \in \mathscr{N}, t \in \mathscr{T} , \omega \in \Omega,\\
    & x_{it\omega} \leq U_i, && \forall i \in \mathscr{N}, t \in \mathscr{T} , \omega \in \Omega,\\
    & \eta_{it\omega} \geq 0, && \forall i \in \mathscr{N}, t \in \mathscr{T} , \omega \in \Omega, \label{LMIP-end}
\end{align}

The objective function \eqref{LMIP-start} minimizes the total cost of all planned actions and future prediction, which includes the cost of stochasticity ($c_s(\phi)$) with holding and lost sales, previously was as constant and thus eliminated in line \ref{RL-MIP} of Algorithm \ref{RLalgo}. Constraints \eqref{old-start} – \eqref{old-end} are the same of Constraints \eqref{totAction} – \eqref{limTruck}. Similar constraints are needed for each subsequent decision epoch and scenario, as defined in Constraints \eqref{new-start} – \eqref{new-end}. The inventories between periods and scenarios are correctly assigned with Constraints \eqref{inv-start} – \eqref{inv-end}, where lost sales are ensured to cover at least the negative inventories. Finally, for the future cost prediction, the end horizon post-decision inventories are assigned in Constraints \eqref{sinv-start} – \eqref{sinv-end}. According to these values, the future cost prediction is computed as the expected value at the end of the lookahead horizon, given by $\frac{1}{|\Omega|}  \sum_{\omega \in \Omega} \hat{v}_w (s'_\omega)$ in the objective function. 
The reader should note that setting $\mathscr{T} \coloneqq \emptyset$ reduces the model to the initial action set defined by Eq. \eqref{totAction} – \eqref{MDP-end}. For $|\mathscr{T}| > 0$, the model allows for a more structured action selection at $t=0$, considering short-term plans over multiple epochs. Increasing $T$ increases the number of decision variables and constraints, requiring a balance between solution quality and computational time. 

As mentioned previously, a lookahead approach without the future cost prediction term ($ \frac{1}{|\Omega|} \sum_{\omega \in \Omega} \hat{v}_w (s'_\omega)$) performed substantially worse in both solution quality and computation time. This hybrid approach, however, proved effective. It can be applied during either the training phase (line \ref{RL-MIP} of Algorithm \ref{RLalgo}) and/or the simulation phase of CRL. Our features make objective function value quadratic. Although solving the mixed-integer quadratic programming (MIQP) for one period is manageable in terms of computation time for our instances, incorporating the lookahead approach during training was intractable due to the large number of solutions required for a proper training of weight vector $w$. Therefore, we propose training RL without lookahead logic using Eq. \eqref{totAction} – \eqref{MDP-end}, and then applying the model in \eqref{LMIP-start} – \eqref{LMIP-end} only during execution of the simulation once the weights are fixed
. This method improves solutions with an acceptable increase in solution time, as discussed in Section \ref{sect:comp}.

\section{Benchmark policies} \label{bench}
In this section, we introduce benchmark policies for comparison with (L)CRL. The selected benchmark policies represent standard approaches traditionally used in related literature. In addition to the two inventory control heuristics, we tested several reinforcement learning baselines, including Proximal Policy Optimization (PPO) \citep{van2023using, dehaybe2024deep}, REINFORCE, TD($\lambda$), and A3C. Despite variations in architecture and training setup, none of these DRL methods produced feasible or stable policies. This reflects a broader limitation of standard DRL in our context: algorithms such as DQN, PPO, or A3C are not equipped to handle large, state-dependent action spaces with hard constraints. While constrained DRL variants exist, they rely on either expectation-based or penalty-driven constraint satisfaction (e.g., Lagrangian or chance-constrained formulations), and actually do not support a state-dependent constrained setting. This contradicts with the DIRP, where actions are implicitly defined by a mixed-integer program and a different feasibility must be enforced at every decision state. Rather than adapting general-purpose DRL to such settings, it is more appropriate to draw from operations research methods developed for constrained, multi-period planning problems. Because of these, we picked two classical inventory control policies; (s,S) and Power-of-Two, that aligns better in the DIRP with a direct deliver setting.

We enhance these inventory control benchmarks to account for the characteristics of our DIRP. In Section \ref{sect:SS}, we propose an $(s,S)$-policy based heuristic, which explores feasible $(s,S)$ pairs for each customer and solves a mixed-integer program to find a set of inventory levels signaling a replenishment. In Section \ref{sect:PO2}, we adapt a Power-of-Two (PO2) heuristic, where replenishments are scheduled at cyclic intervals that are powers of two, and a mixed-integer program ensures vehicle capacity limits are respected. The PPO construction and additional discussion of DRL limitations are provided in Appendix~\ref{sec:DRL}.

\subsection{$(s,S)$-policy based heuristic} \label{sect:SS}

In this section, we propose an iterative $(s,S)$-policy heuristic. Each customer has a set of candidate $(s,S)$-policies; $(s_i^k,S_i^k)$, each with an estimated cost $c_i^k$ and average vehicle usage $q_i^k$ per period. We select one $(s,S)$-policy per customer via a MIP that enforces a total usage threshold $q_{\text{MIP}}$. We then iteratively adjust $q_{\text{MIP}}$ to search the solution space in order to find an efficient combined solution, addressing limitations of the vehicle fleet. The MIP is defined as;
\begin{align}
    \text{min  } &  \sum_{i \in \mathscr{N}} \sum_{k \in \mathscr{K}_i} c_i^k z_i^k , \label{SSMIPStart} \\
    \text{s.t  } & \sum_{k \in \mathscr{K}_i} z_i^k = 1, \label{SScust} \\
    & \sum_{i \in \mathscr{N}} \sum_{k \in \mathscr{K}_i} q_i^k z_i^k \leq q_{\text{MIP}}, \label{constQMIP}\\
    &  z_i^k \in \{0,1\}, \label{SSMIPEnd}
\end{align}
where $z_i^k$ is a binary variable that is 1 if the $k^\text{th}$ $(s,S)$ policy is chosen for customer~$i$, $(s_i^k,S_i^k)$. Constraint \eqref{SScust} ensures exactly one policy is selected per customer. Constraint \eqref{constQMIP} enforces the total average vehicle usage limit $q_{\text{MIP}}$. 

The $(s,S)$-policy based heuristic algorithm that implements this MIP, given in Algorithm~\ref{SSalgo}, starts with an evaluation of all feasible $(s,S)$ pairs for each customer $i$, each denoted as $(s_i^k,S_i^k)$ where $k \in \mathscr{K}_i$ and $0 \leq s_i^k < S_i^k \leq U_i$. For each pair, we simulate the individual customer policy as if supply and vehicles were unlimited. From this simulation, we record the holding costs at customer locations, lost sales, and transportation cost $c_i^k$ per period, as well as the average vehicle usage $q_i^k$ per period. We then identify the most cost-effective $(s,S)$ pair for each customer and calculate $ q_{\text{MIP}} \coloneqq \sum_{i \in \mathscr{N}} q_i^{\argmin_k{c_i^k}},$ which serves as an initial threshold on total vehicle usage. These pairs of $\argmin_k{c_i^k}$ do not necessarily provide a near-optimal system solution, as the combined policy ignores limited supply availability, the capacity of the vehicle fleet, and the replenishment planning of customers. This is correct even for the cases where $q_{\text{MIP}} < q$, as this threshold only focuses on average levels of vehicle use and does not consider the random distribution of the number of vehicles used on different periods.

\begin{algorithm}[h]
    \caption{$(s,S)$-policy based heuristic} \label{SSalgo}
    \begin{algorithmic}[1] 
        \State \textbf{Algorithm Parameters:} An increment value of $\xi \coloneqq \xi_0$, a multiplier of $m > 1$, and $t^*$ for stopping condition.
        \For{each customer $i \in \mathscr{N}$} \label{ssforstart}
            \State $k \coloneqq 0$
            \For{each $s \in \mathbb Z$ in $[0,U_i)$}
                \For{each $S \in \mathbb Z$ in $(s,U_i]$}
                    \State Simulate $(s,S)$ for customer $i$ as enough inventory and vehicles are always available.
                    \State Record $c_i^k$ as the cost of holding, lost sales, and transport per period. 
                    \State Record $q_i^k$ as the average number of vehicles used per period.
                    \State $k++;$
                \EndFor
            \EndFor
        \EndFor \label{ssforend}
        \State $q_{\text{MIP}} \coloneqq \sum_{i \in \mathscr{N}} q_i^{\argmin_k{c_i^k}}$
        \While {no improvement for the last $t^*$ iterations}
            \State Solve MIP in Eq. \eqref{SSMIPStart} - \eqref{SSMIPEnd} with  $q_{\text{MIP}}$.
            \State Simulate the selected policies for all customers and observe $V(s)$ per turn
            \If{the best solution so far}
            \State $\xi \coloneqq m \xi  $
            \Else
            \State $\xi \coloneqq \xi_0$
            \EndIf
            \State $q_{\text{MIP}} \coloneqq q_{\text{MIP}} - \xi$.
        \EndWhile
    \end{algorithmic}
\end{algorithm}

In order to tackle this, as also shown by preliminary experiments, we explore other solutions where the average number of vehicle usage is below the initially selected threshold value. For this search, we iteratively lower $q_{\text{MIP}}$ and solve the MIP \eqref{SSMIPStart}--\eqref{SSMIPEnd} for each new value. The step size for reducing $q_{\text{MIP}}$ is given by an increment $\xi$, initially $\xi_0$. If an iteration improves the overall cost, we multiply $\xi$ by a factor $m>1$ to reduce $q_{\text{MIP}}$ more aggressively in subsequent iterations; otherwise, we reset $\xi$ to $\xi_0$. After each MIP solution, we simulate the implied $(s,S)$-policies to observe the system cost $V(s)$ per period. If the chosen policies attempt to use more vehicles than are actually available, we randomly skip some replenishments while keeping the fleet at capacity in use. The iterative search terminates when no further improvements are observed over $t^*$ consecutive iterations, and the best $(s,S)$ selection found is returned as the final solution. 

\subsection{Power-of-Two based heuristic} \label{sect:PO2}

Power-of-two (PO2) policies are common in inventory control literature, especially for a single-warehouse multi-retailer system \citep{chu2010power, li2019inventory}. These policies focus on the scheduling of replenishments based on a frequency for two consecutive orders, rather than specific inventory levels signaling a replenishment as in $(s,S)$-policies. Typically, under PO2, each customer is visited cyclically at intervals that are powers of two, such as at every 1, 2, 4, and 8 periods, differing for each customer. In this section, we adapt a PO2-based heuristic to suit the dynamics of the DIRP.

Let $\mathscr{T}$ be the set of periods between two consecutive replenishments. Ideally, $\mathscr{T} \coloneqq \mathbb Z_{> 0}$. However, adhering to PO2, $\mathscr{T}$ is limited to powers of two, thus defined as $\mathscr{T} \coloneqq \{1, 2, 4, 8, \dots, 2^{\tau}\}$, where $\tau$ is a predetermined upper bound. For each customer $i$ and each replenishment interval $t \in \mathscr{T}$, we find the optimal order-up-to level. For this, we calculate the exact expected total cost per customer per period, including replenishment, holding, and lost sales costs. These exact expected costs are obtained using the discrete support of the demand distributions, $\Phi$, explained in Section \ref{probdef}. The process is repeated for all discrete feasible order-up-to levels in $[0, U_i]$, and the among these the lowest per-period cost is selected as $c_i^t$, corresponding to customer $i$ and interval $t$. To integrate these individual replenishment schedules into a combined plan, we formulate a MIP that selects exactly one replenishment interval for each customer while ensuring the vehicle availability.  Let $z_i^t$ be a binary decision variable, having a value of $1$ if customer $i$ is scheduled for visits every $t$ periods. The MIP is given as follows:

\begin{align}
    \text{min  } &  \sum_{i \in \mathscr{N}} \sum_{t \in \mathscr{T}} c_i^t z_i^t \label{po2start} , \\
    \text{s.t  } & \sum_{t \in \mathscr{T}} z_i^t = 1, \label{custAssign}\\
    & \sum_{i \in \mathscr{N}} \sum_{t \in \mathscr{T}} \frac{1}{t} z_i^t \leq q, \label{const:PO2} \\
    &  z_i^t \in \{0,1\}. \label{po2end}
\end{align}

The objective is to minimize total cost associated with customers, ensuring each customer is assigned to one replenishment schedule by Constraint \eqref{custAssign}. Constraint \eqref{const:PO2} ensures that we have always enough vehicles available for selected PO2 policies. 

PO2 policies not only provide near-optimal solutions in inventory control \citep{federgruen1992simple}, but are also crucial for defining a constraint in the MIP structure that reflects the vehicle usage limit, i.e. Constraint \eqref{const:PO2}. The constraint holds when  $\mathscr{T}$ consists solely of PO2 intervals, which ensures there exists a cyclic replenishment schedule adhering to the vehicle fleet size limit, $q$. However, if non-PO2 intervals are allowed, this property does not necessarily hold. For example, let $\mathscr{T} \coloneqq \{1, 2, 3, 5\}$. Consider a fleet of \( q = 2 \) vehicles and three customers \( \mathscr{N} \coloneqq \{1, 2, 3\} \), where each customer \( i \) is scheduled for visits once every \( i \) periods. Constraint \eqref{const:PO2} implies \( \sum_{i \in \{1, 2, 3\}} \frac{1}{i} \leq 2 \), which is satisfied. However, no cyclic schedule adheres to using at most \( q = 2 \) vehicles per period. Specifically, customer 1 uses exactly one vehicle each day. Furthermore, customers 2 and 3 both require replenishment simultaneously every six periods for all possible cyclic schedules. As a result, once every six periods, the solution requires three vehicles, not adhering to the vehicle fleet size limit.

Our computational experiments confirm that the MIP solution produces feasible cyclic schedules when using only PO2 intervals. Once the MIP is solved, the PO2-based policy defines one such feasible cyclic schedule, which we then simulate to estimate the average cost per period, including holding costs at the supplier that are not explicitly modeled in the MIP.

\section{Computational Experiments} \label{sect:comp}

In this section, we evaluate our Constrained Reinforcement Learning (CRL) approach, and its enhanced Lookahead-based CRL (LCRL) in terms of solution quality and computational time. The experiments are conducted on a platform with an AMD 7763 CPU @ 2.45 GHz and 16GB RAM. For LCRL, we use 16 CPU cores. All the experiments are executed using C++20 with Gurobi 11.0.0 for MIQP optimization. We employ Proximal Policy Optimization (PPO) with neural networks defined using Pytorch's C++ API, LibTorch 2.2.2, see Appendix \ref{sec:DRL} for its discussion. For all the benchmark methods, the suggested policies are tested via simulations to derive an average total cost per period. There exists one simulation per instance for each solution method, except for the $(s,S)$-based heuristic, which requires multiple simulations as detailed in line 14 of Algorithm \ref{SSalgo}. We report solution times for all benchmarks by excluding simulation durations to derive these costs, and including only the time taken to derive a policy. See Appendix \ref{sec:instGen} on how the test instances are generated, as well as other problem specific parameters.

In Section \ref{sec:VI}, we compare our CRL and LCRL with value iteration (VI) algorithm to test its performance with the optimal solution. 
In Section \ref{sec:opt}, we provide a numeric analysis of the optimal policy structure of the VI results, as well as the effects of stochasticities of supply and demand on these structures.
In Section \ref{sec:SSPO2}, we compare the algorithm with our proposed benchmarks; $(s,S)$-policy-based heuristic, and Power-of-Two based heuristic (PO2) for realistic-sized instances in terms of the solution quality and computational time. 
In Section \ref{sec:sensitivity}, we investigate the impact of learning parameters and training duration on the performance and convergence of the CRL algorithm, benchmarking these results against VI to assess their effectiveness. 

\subsection{Value iteration on small instances} \label{sec:VI}

Value iteration is a known dynamic programming method for solving Markov decision processes to optimality. In our setting, we adapt value iteration with a decomposition approach for tractable computations. The values, $V_i$, are updated sequentially for each $i \in \mathscr{N}^0$ before computing the overall Bellman update, $V_{N+1}$, which represents the predicted value of the Bellman equation \eqref{bellman}. Given our undiscounted horizon, these values increase indefinitely, so we monitor marginal differences between consecutive iterations. The process continues until the maximum marginal change falls below a predefined threshold. The algorithm then returns an optimal action policy for each state, $\pi(x)$. See Algorithm \ref{VI} for details.

\begin{algorithm}[h]
    \caption{Value Iteration} \label{VI}
    \begin{algorithmic}[1] 
        \State \textbf{Input:} Instance and discrete support $\Phi$.
        \State \textbf{Output:} Policy vector $\pi$. 
        \While {$\max_x \Delta V_{N+1}(x) > \epsilon$} 
            \For{each customer $i \in \{ N, N-1,\dots, 1 \}$}
                \For{each state $s$}
                    \State $V_i(s) \coloneqq \mathbb{E}\left[ \sum_{\phi_i \in \Phi_i} \left( c_s(\phi_i) + V_{i+1}(s_0, \ldots, s_{i-1}, s_i - \phi_i, s_{i+1}, \ldots, s_N) \right) \right]$
                \EndFor
            \EndFor
            \For{each state $s$}
                \State $V_{0}(s) \coloneqq \mathbb{E}\left[ \sum_{\phi_0 \in \Phi_0} \left( c_s(\phi_0) + V_{1}(s_0 + \phi_0, s_1, s_2, \ldots, s_N) \right) \right]$
            \EndFor
            \For{each state $x$}
                \State $V_{N+1}(x) \coloneqq \min_{a \in A(x)} \left\{  c_x(a)  +  V_{0}(s(x,a))  \right\} $
                \State $\pi(x) \coloneqq \argmin_{a \in A(x)} \left\{  c_x(a)  +  V_{0}(s(x,a))  \right\} $
            \EndFor      
        \EndWhile
    \end{algorithmic}
\end{algorithm}

Even in our decomposed setting, value iteration requires substantial computational resources, which becomes impractical for moderate-sized instances of our generator, described in Appendix \ref{sec:instGen}. 
To manage this, we have formulated toy instances inspired by those used by \cite{kaynov2023deep} where we set $N=3$ and $q=2$, with several modifications for our DIRP. We set discrete demand means per customer uniformly within the range of $[2,4]$
. Furthermore, the inventory capacities, $U_i$, are selected as twice the expected demand per customer and $1.5$ times the expected supply for the supplier. 
Thirdly, the vehicle capacity is set as $C \coloneqq 1.25  \mathbb{E}\left[ \sum_{i \in \mathscr{N}} \phi_i\right] / q$, limiting the amount of feasible discrete actions. Finally, a cost vector of $(W,w,h_s,h_c,\ell,\rho) = (15,1.5,2,4,15,2.5)$ is selected to ensure a balance between inventory costs and other costs.

In an average of $10$ instances sized at $N=3$ and $q=2$, the Constrained Reinforcement Learning (CRL) and lookahead-based CRL (LCRL) method demonstrated an optimality gap of $1.8\%$ and $2.3\%$. This performance is comparable to the gaps of $1-3\%$ reported by \cite{kaynov2023deep} for their Deep Reinforcement Learning (DRL) approach in instances with lost sales, although their inventory control problem does not have the complexities of limited vehicle availability and vehicle inventory capacities. Furthermore, the heuristic policies $(s,S)$ and PO2 show optimality gaps of $3.9\%$ and $6.9\%$, respectively. Lower optimality gaps are reported in the inventory management literature \citep{temizoz2025deep}; however, these are achieved where the action space is enumerable, and with no vehicle limits. In contrast,  our model deals with an intractably large number of actions, at least requiring modifications to conventional deep learning approaches or methods such as CRL; see Appendix \ref{sec:DRL} for further discussion. 

\subsection{Optimal policy structure} \label{sec:opt}

\begin{figure}[htp]
    \centering
    \begin{subfigure}{0.3\textwidth}
        \centering
        \includegraphics[width=\textwidth]{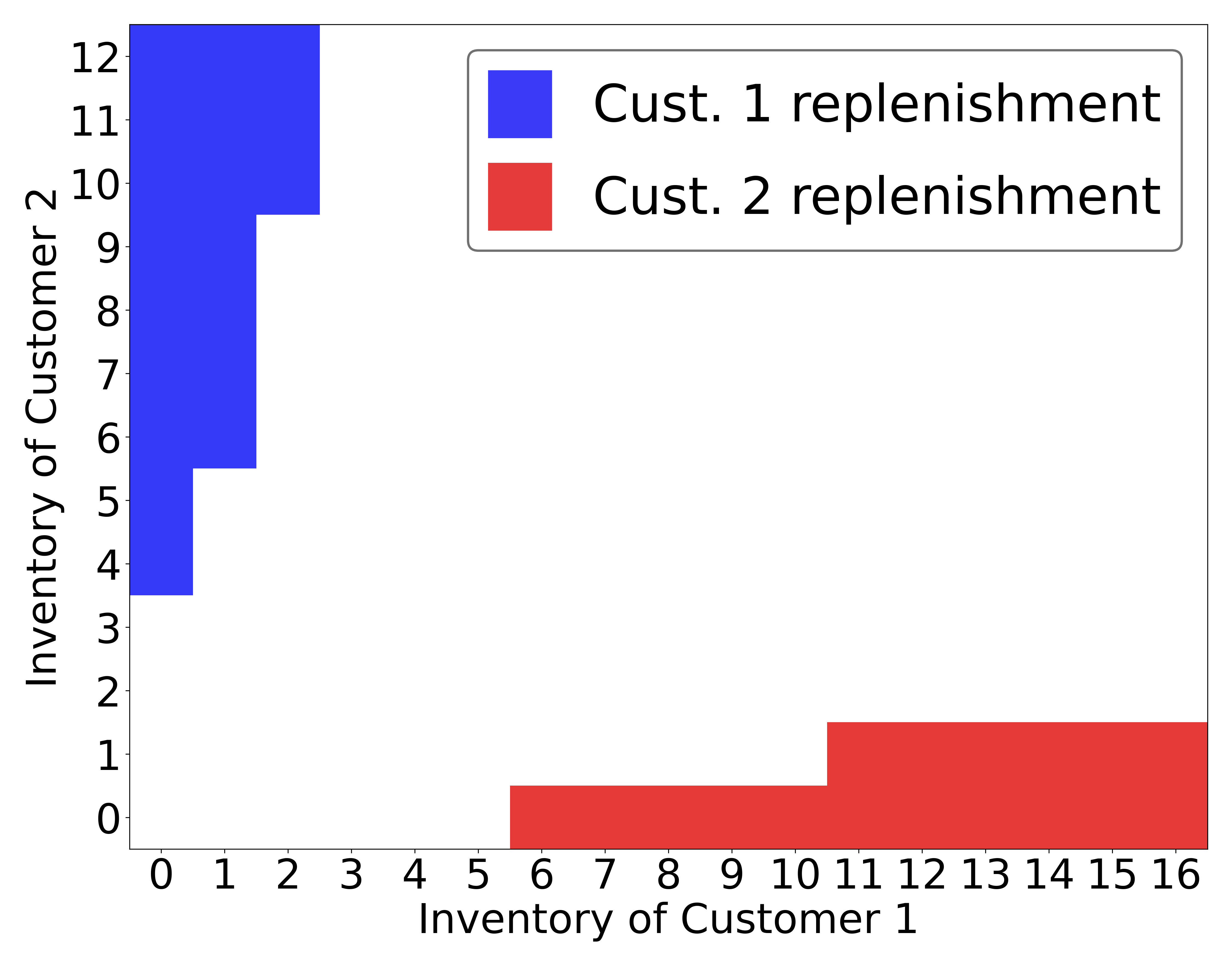}
        \caption{$x_0 = 3$} \label{x03}
    \end{subfigure}
    \begin{subfigure}{0.3\textwidth}
        \centering
        \includegraphics[width=\textwidth]{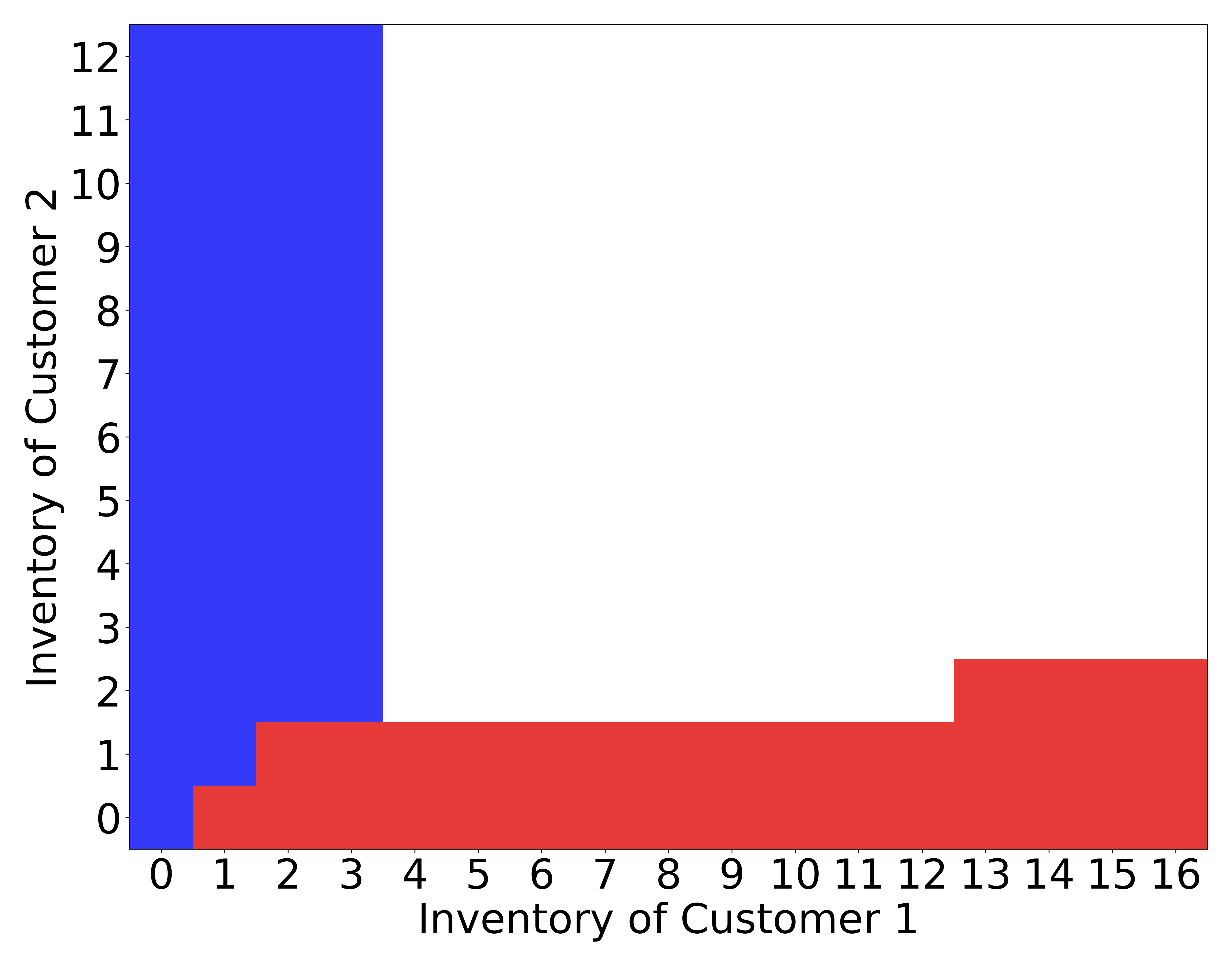}
        \caption{$x_0 = 5$} \label{x05}
    \end{subfigure}
    \begin{subfigure}{0.3\textwidth}
        \centering
        \includegraphics[width=\textwidth]{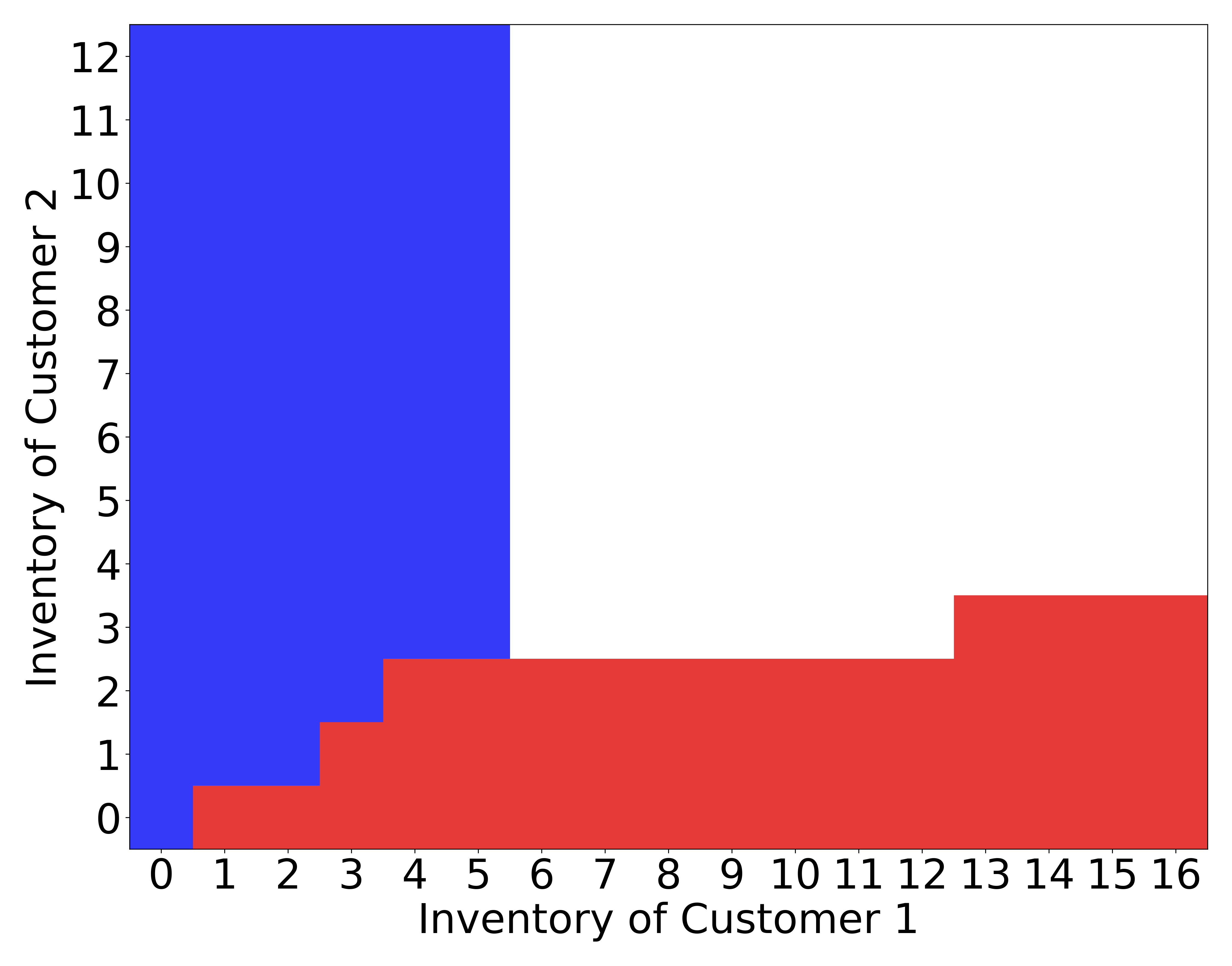}
        \caption{$x_0 = 8$} \label{x08}
    \end{subfigure}

    \begin{subfigure}{0.3\textwidth}
        \centering
        \includegraphics[width=\textwidth]{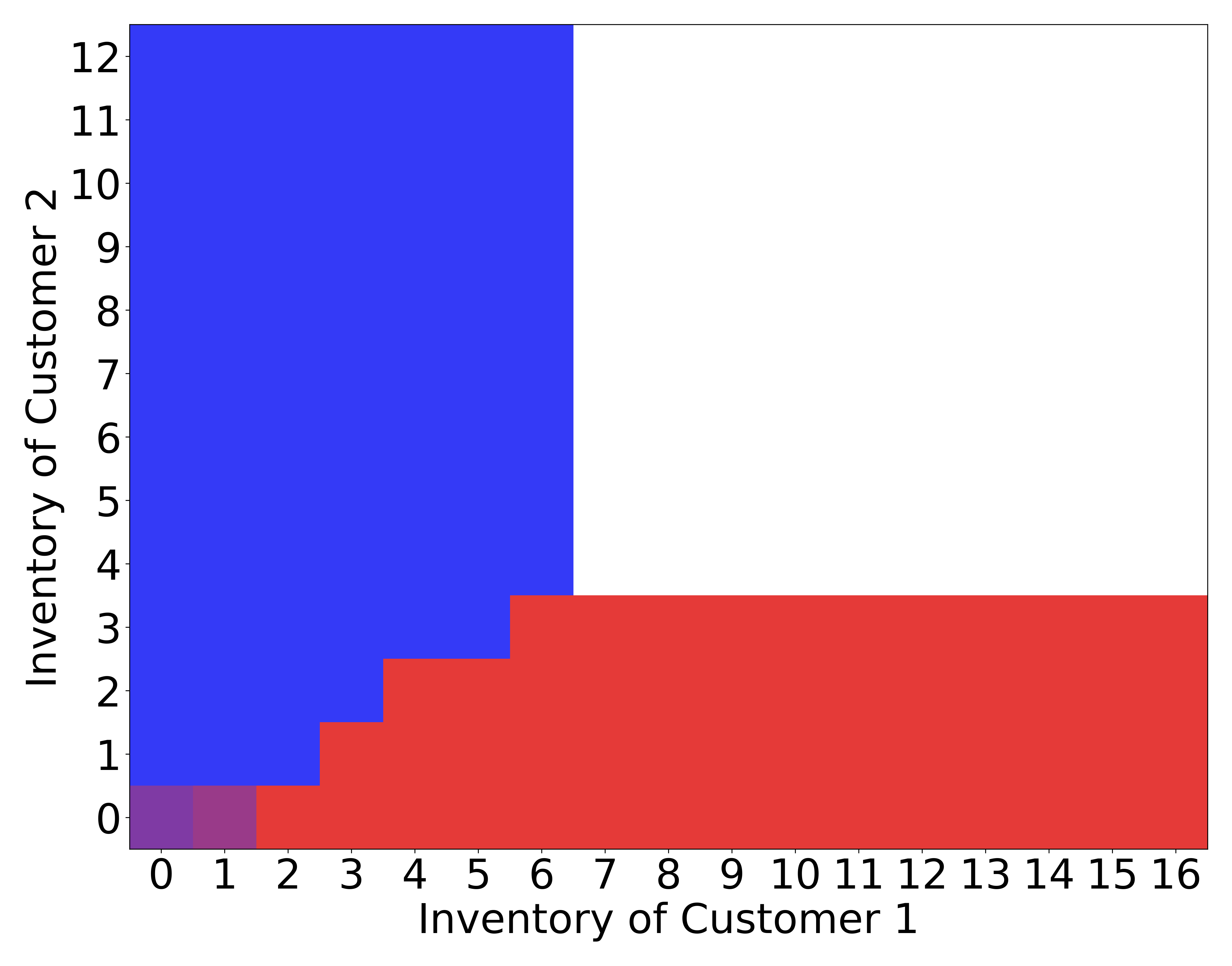}
        \caption{$x_0 = 11$}  \label{x011}
    \end{subfigure}
    \begin{subfigure}{0.3\textwidth}
        \centering
        \includegraphics[width=\textwidth]{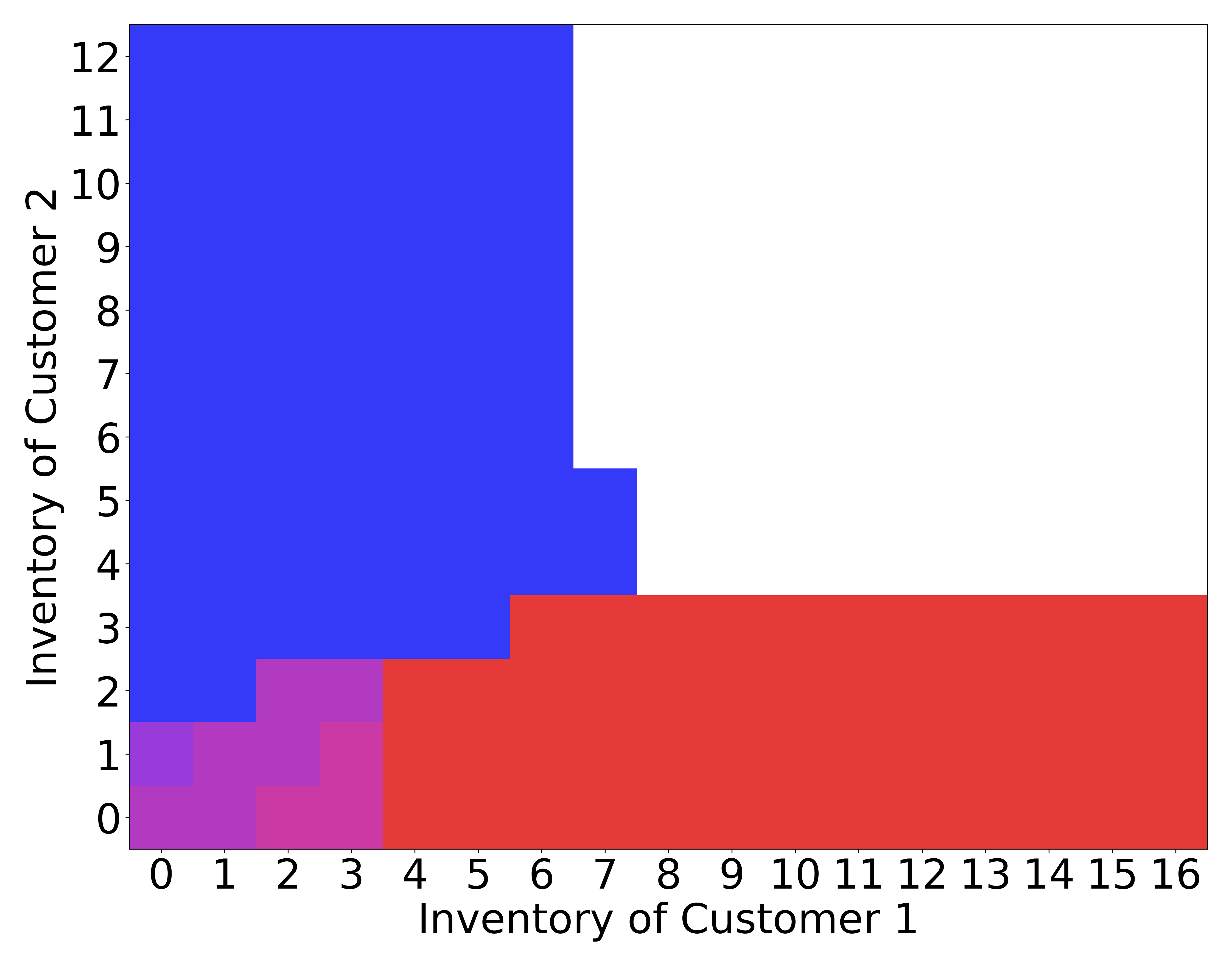}
        \caption{$x_0 = 14$} \label{x014}
    \end{subfigure}
    \begin{subfigure}{0.3\textwidth}
        \centering
        \includegraphics[width=\textwidth]{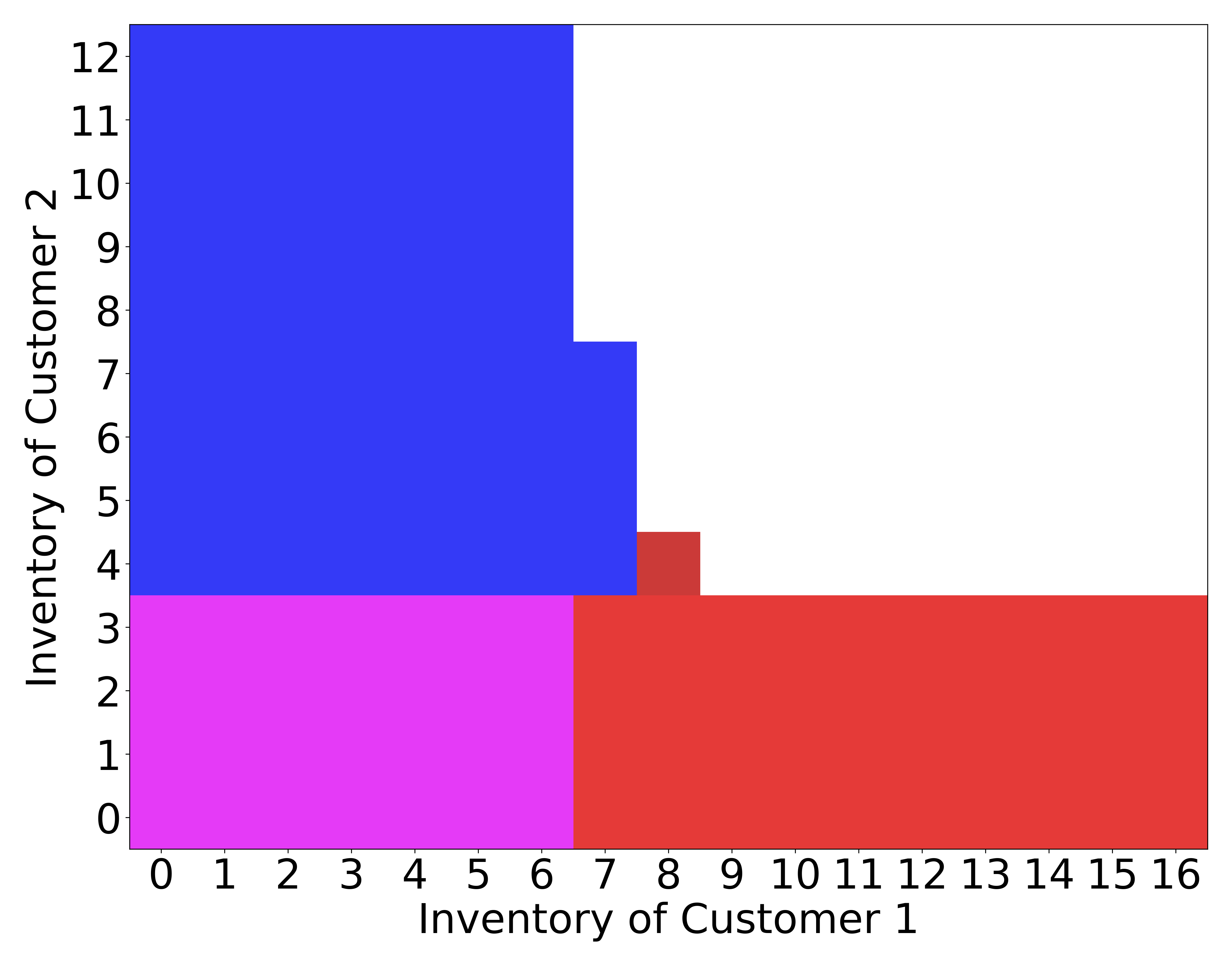}
        \caption{$x_0 = 23$} \label{x023}
    \end{subfigure}
    
    \caption{Optimal replenishment policies at varying supplier inventory levels. Blue and red regions indicate replenishments directed to customer 1 and customer 2, respectively, with color intensity reflecting the replenishment quantity. Purple areas represent shipments to both customers.}
    \label{figoptpol}
\end{figure}

In this section, we examine the structure of the optimal replenishment policies obtained via value iteration. We generate an instance identical to the ones in Section \ref{sec:VI}, but with increased inventory capacities to better capture policy effects. In a three-customer setup with two vehicles, we analyze the states where customer 3 has full inventory and thus does not require replenishment. Figure \ref{figoptpol} provides the optimal replenishment policies for various levels of supplier inventory, where the axes are inventory quantities of customer 1 and 2. Blue and red regions indicate replenishment directed to customer 1 and customer 2, respectively. Purple is a mix of both replenishments, and white region represents no vehicles dispatched to any customers. Color intensity reflects the quantity allocated to each customer.

In Figure \ref{x03}, where supplier's inventory is three, $x_0 = 3$, the optimal policy allocates replenishments only when the replenished customer inventory is too low and the other customer's inventory is high enough. For example, it is observed that in the upper left part customer 1 is replenished when $x_1 = 0$ and $x_2 \in [4,12]$. The range shrinks to $x_2 \in [6,12]$ when $x_1 = 1$, since more inventory at customer 1 decentivize the replenishment. In general, at $x_0 = 3$, replenishments are infrequent, as the fixed transportation cost makes sending three units less favorable. Additionally, with this limited supply, each unit is more valuable, creating an incentive to further skip shipments for next replenishments. Due to this, the larger the supply inventory, the more replenishments begin to occur, as seen in Figure \ref{x05} and \ref{x08}.
At $x_0 = 11$ in Figure \ref{x011}, the optimal policy suggests sending two vehicles when $x_1 \in \{0,1\}$ and $x_2 = 0$, dividing the available supplier stock between customers as both need urgent replenishment. This intersection of replenishments grows further with larger supply inventory as seen in Figure \ref{x014}. 
In Figure \ref{x023}, the replenishment quantity for one customer becomes less dependent on the inventory of the other, as the abundant supply and vehicles allow the needs of each customer to be met independently. In general, the differences among these six graphs in Figure \ref{figoptpol} highlight the need to account for stochastic and limited supply in DIRP models for accurate policy design, a factor often overlooked in previous DIRP research.

\begin{figure}[htp]
    \centering
    \begin{subfigure}{0.3\textwidth}
        \centering
        \includegraphics[width=\textwidth]{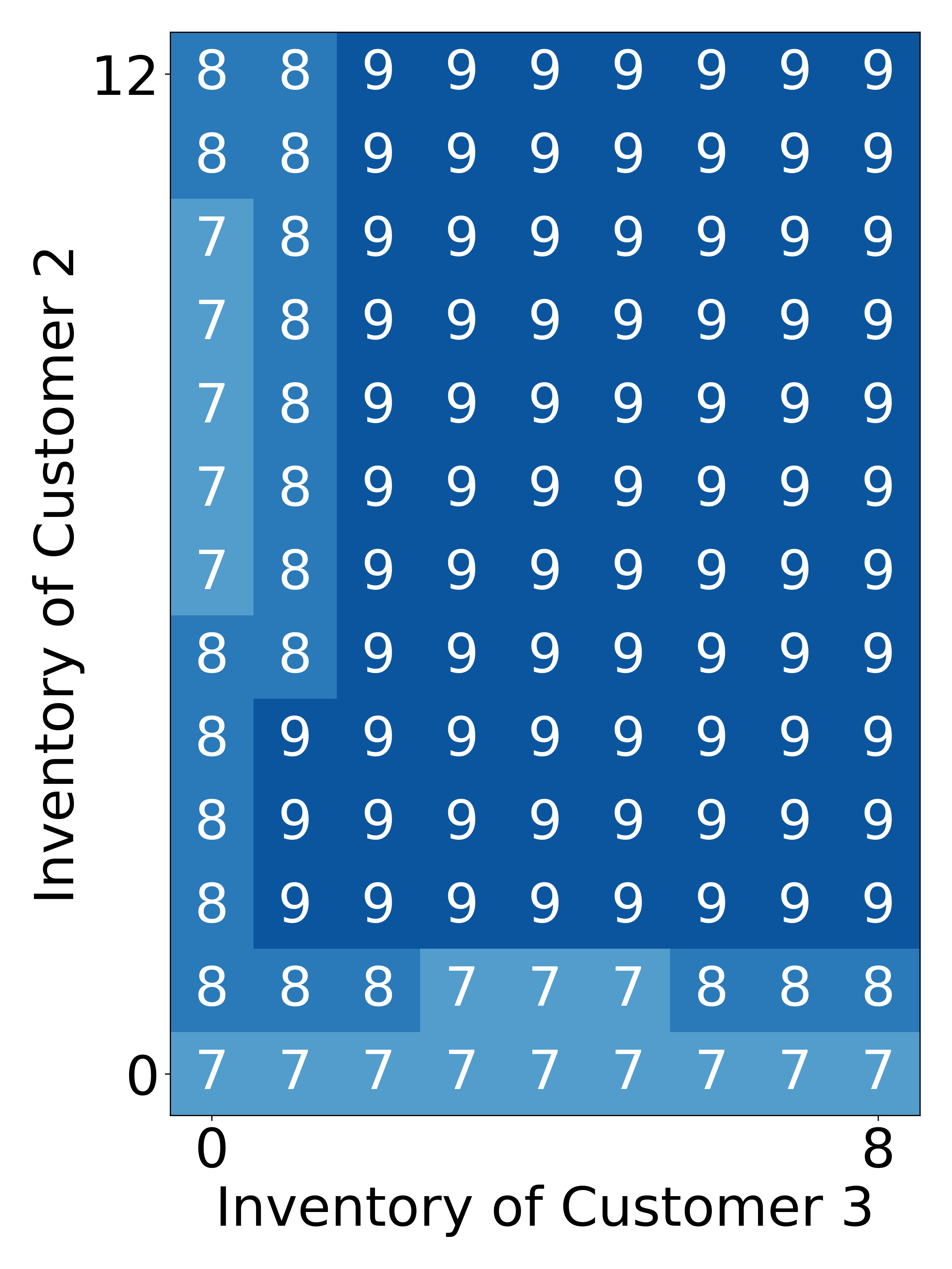}
        \caption{Replenishment to Customer 1 ($a_1$)} \label{4a}
    \end{subfigure}
    \begin{subfigure}{0.3\textwidth}
        \centering
        \includegraphics[width=\textwidth]{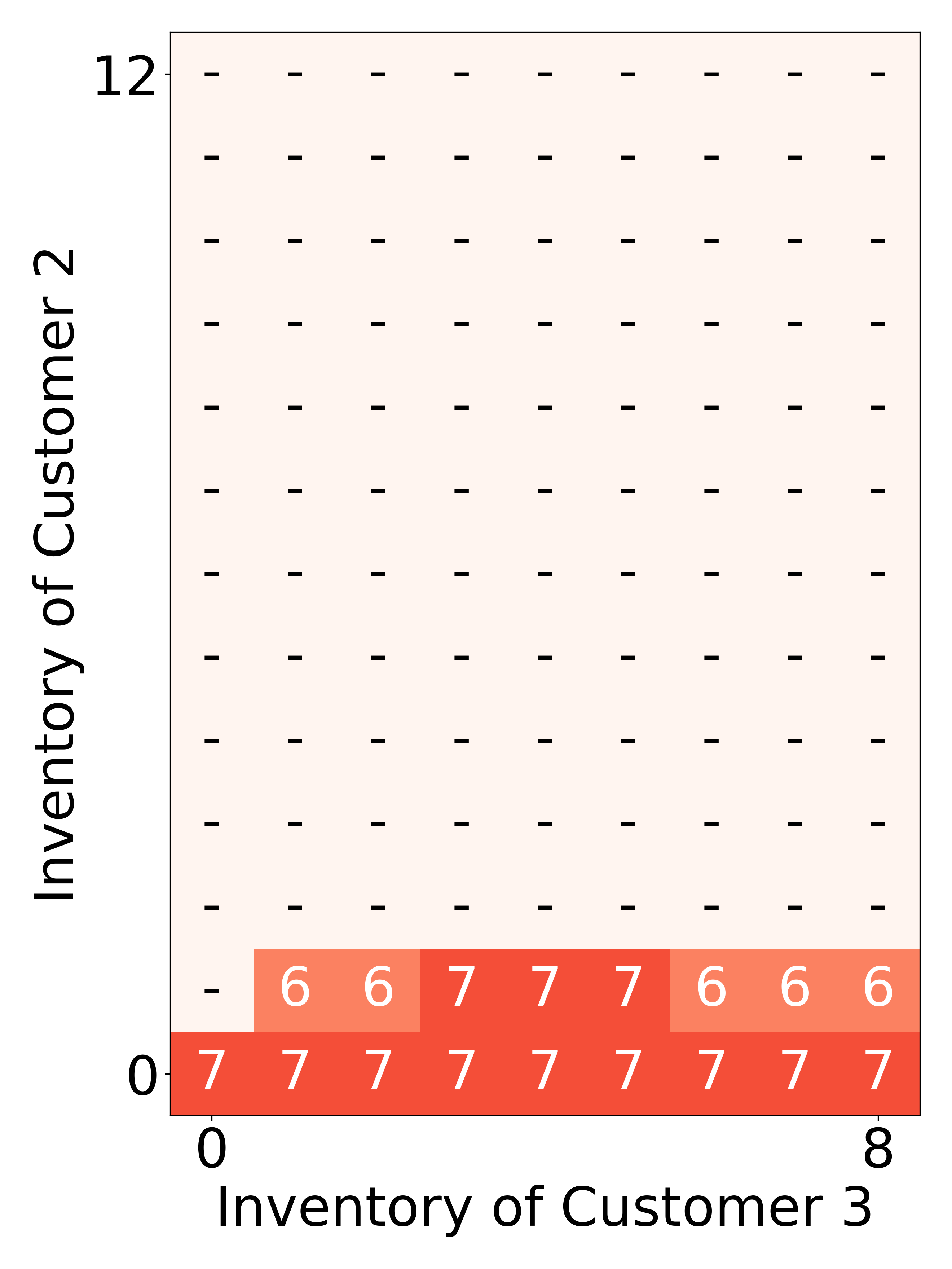}
        \caption{Replenishment to Customer 2 ($a_2$)} 
    \end{subfigure}
    \begin{subfigure}{0.3\textwidth}
        \centering
        \includegraphics[width=\textwidth]{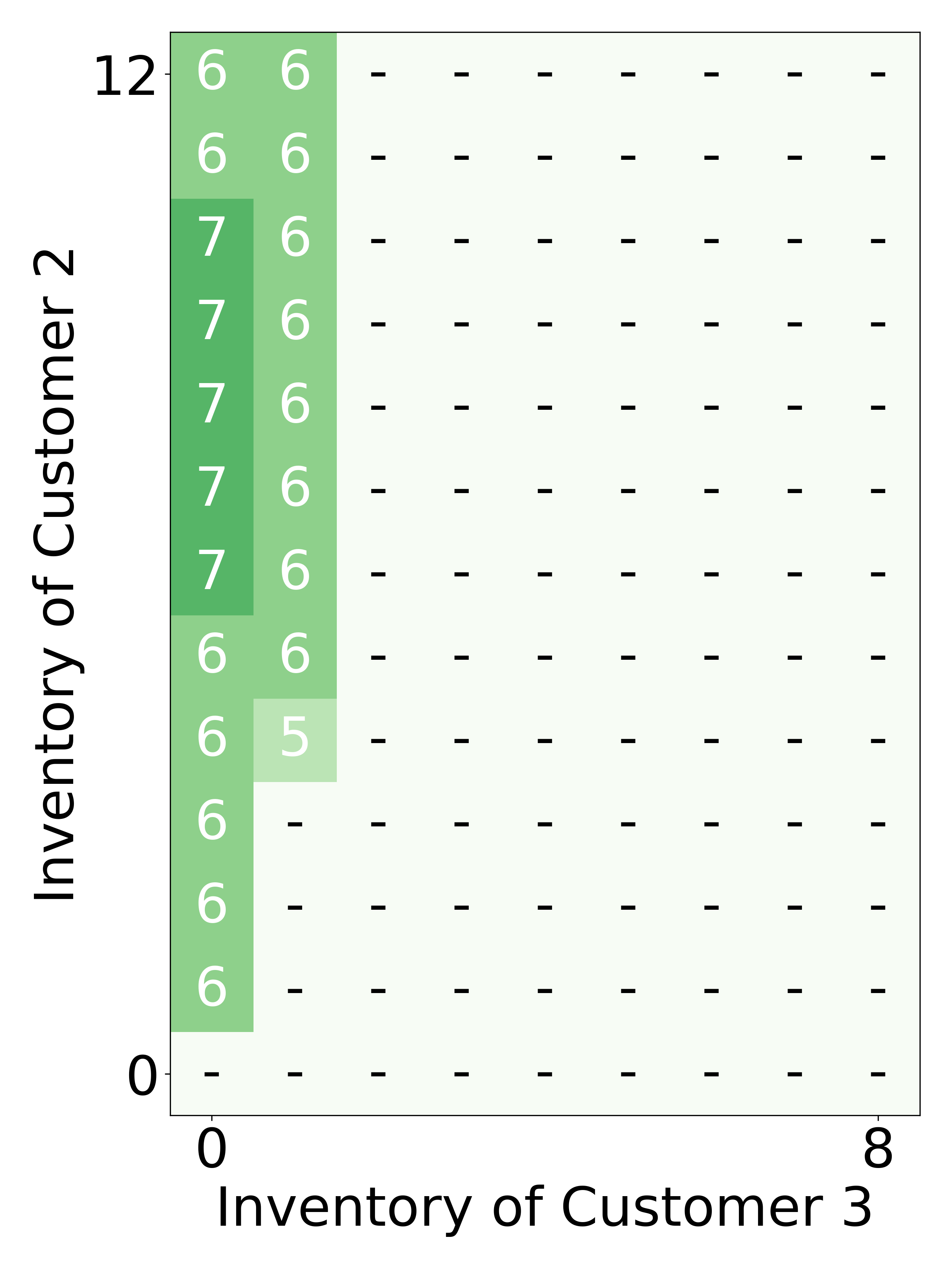}
        \caption{Replenishment to Customer 3 ($a_3$)} \label{4c}
    \end{subfigure}
    
    \caption{Optimal replenishment policies at $(x_0, x_1) = (14,0)$.  Replenishment quantities are given in cells. Darker color intensity reflects larger replenishment quantities.}
    \label{figoptpol2}
\end{figure}

Next, Figure \ref{figoptpol2} analyzes the states where the supplier and customer 1 have $(x_0, x_1) = (14,0)$, while customer inventories $x_2$ and $x_3$ vary across their feasible levels. 
Each subfigure \ref{4a} - \ref{4c} illustrates the optimal replenishment policy for one customer, respectively. 
Notably, the boundaries of replenishment regions do not extend from the extreme corners as in the previous figure; instead, they initiate along the edges. For instance, as inventory shifts from $(x_2, x_3) = (5,0)$ to $(x_2, x_3) = (11,0)$, looking all the subfigures, it is seen that the policy first recommends $(a_1, a_2, a_3) = (8,0,6)$, then $(7,0,7)$, and subsequently returns to $(8,0,6)$. This interaction between replenishment quantities for different customers reflects the inherent complexity of DIRP even in a small-scale instance, as the inventory level of one customer influences the replenishment allocation for others in a nonlinear manner.

\begin{figure}[htp]
    \centering
    \begin{subfigure}{0.3\textwidth}
        \centering
        \includegraphics[width=\textwidth]{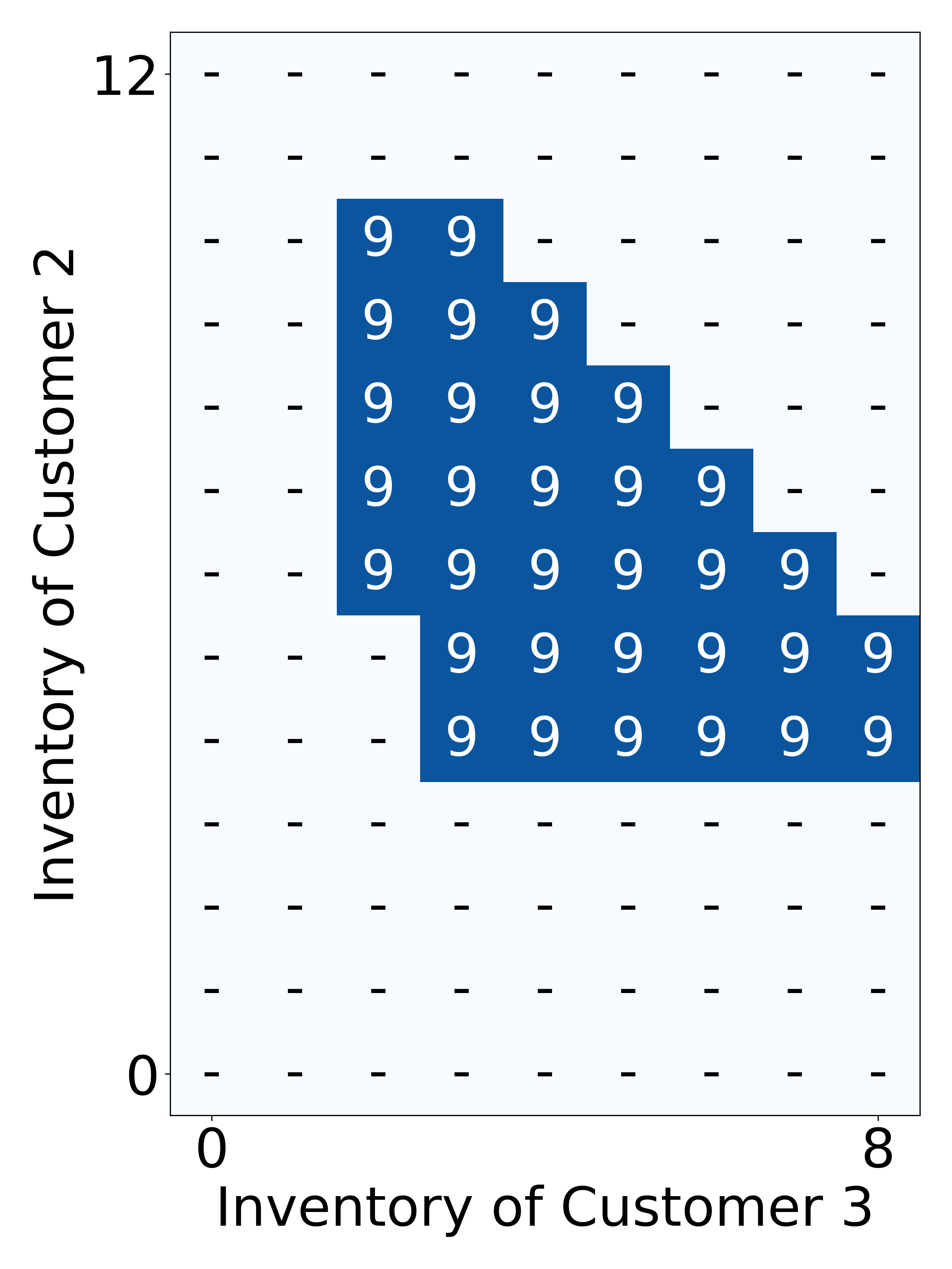}
        \caption{Replenishment to Customer 1}  \label{figoptpol31}
    \end{subfigure}
    \begin{subfigure}{0.3\textwidth}
        \centering
        \includegraphics[width=\textwidth]{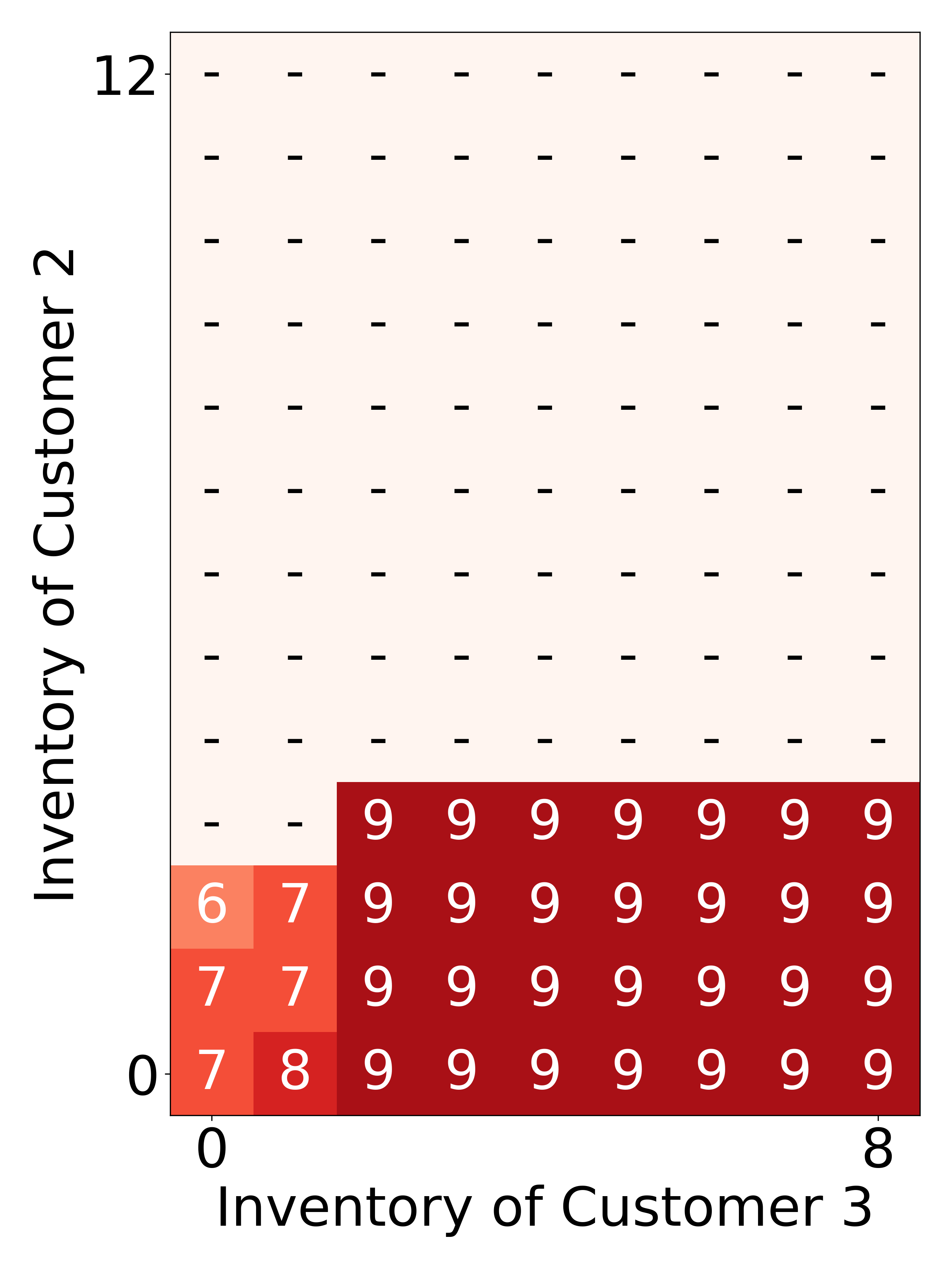}
        \caption{Replenishment to Customer 2}  \label{figoptpol32}
    \end{subfigure}
    \begin{subfigure}{0.3\textwidth}
        \centering
        \includegraphics[width=\textwidth]{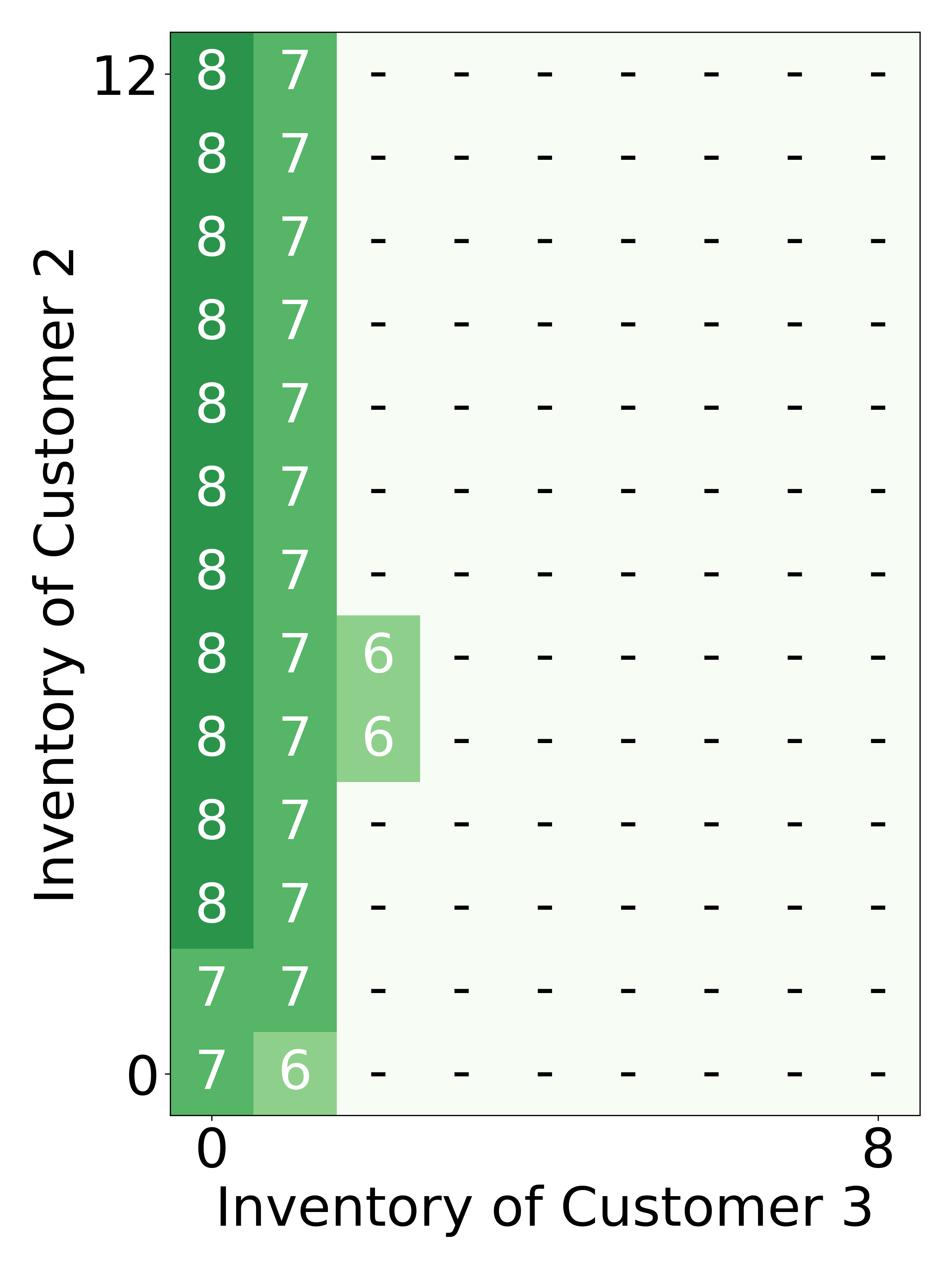}
        \caption{Replenishment to Customer 3}  \label{figoptpol33}
    \end{subfigure}
    
    \caption{Optimal replenishment policies at $(x_0, x_1) = (14,7)$.}
    \label{figoptpol3}
\end{figure}

In Figure \ref{figoptpol3}, we further increase customer 1's inventory level and analyze the states at $(x_0, x_1) = (14,7)$. Here, in Figures \ref{figoptpol32} and \ref{figoptpol33}, the replenishment policies resemble those observed in Figure \ref{figoptpol}, with replenishment boundaries shifting from Figure \ref{figoptpol2} to Figure \ref{figoptpol3} in a way that mirrors the patterns seen in Figure \ref{figoptpol}. 
We observe that optimal replenishment policies of customer 2 and 3 between Figure \ref{figoptpol2} and Figure \ref{figoptpol3} are significantly different. Since this difference occurs on the change of inventory of customer 1, it supports our argument that our MDP construction is hard to decouple between customers.
Moreover, a distinct boundary appears in Figure \ref{figoptpol31}, where replenishment to customer 1 ceases approximately if $x_2 + x_3 > 13$. 
This seems counterintuitive, since higher inventory at customers 2 and 3 might suggest a focus on replenishing customer 1. However, the policy likely aims to prevent a scenario in which all three customers require simultaneous replenishment in the next periods. With only two vehicles available, such a situation would force one customer to be skipped, potentially incurring high lost sales costs. Thus, when $x_2 + x_3 \leq 13$, customer 1 is replenished to prevent skipping a future replenishment, but this results in higher inventory holding costs. Ideally, when $x_2 + x_3 > 13$, the optimal policy suggests delaying replenishment to customer 1, as it becomes more cost-effective to meet its needs in the next periods, given that sufficient vehicle capacity will likely be available. We test this logic by assuming there exists only one vehicle in the system instead of two. In this new problem instance, the threshold shifts to about $x_2 + x_3 = 18$, extending the replenishment boundary and further reducing the likelihood of vehicle shortages in future periods, aligning with our explanation of prioritizing future vehicle availability over ideal replenishment schemes.

\begin{figure}[htp]
    \centering
    \begin{subfigure}{0.225\textwidth}
        \centering
        \includegraphics[width=\textwidth]{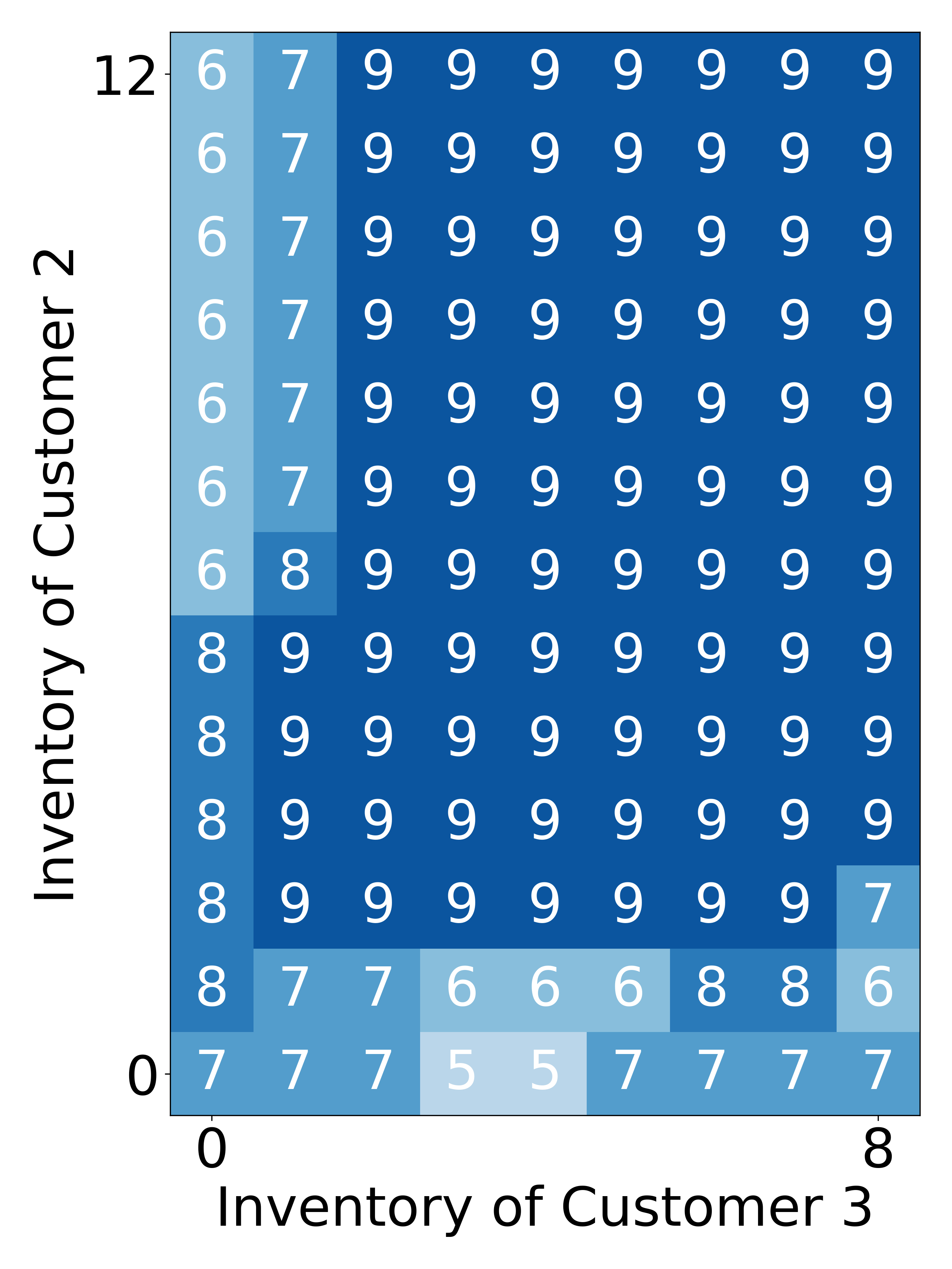}
        \caption{Customer 1 with LS-LD}
    \end{subfigure}
    \begin{subfigure}{0.225\textwidth}
        \centering
        \includegraphics[width=\textwidth]{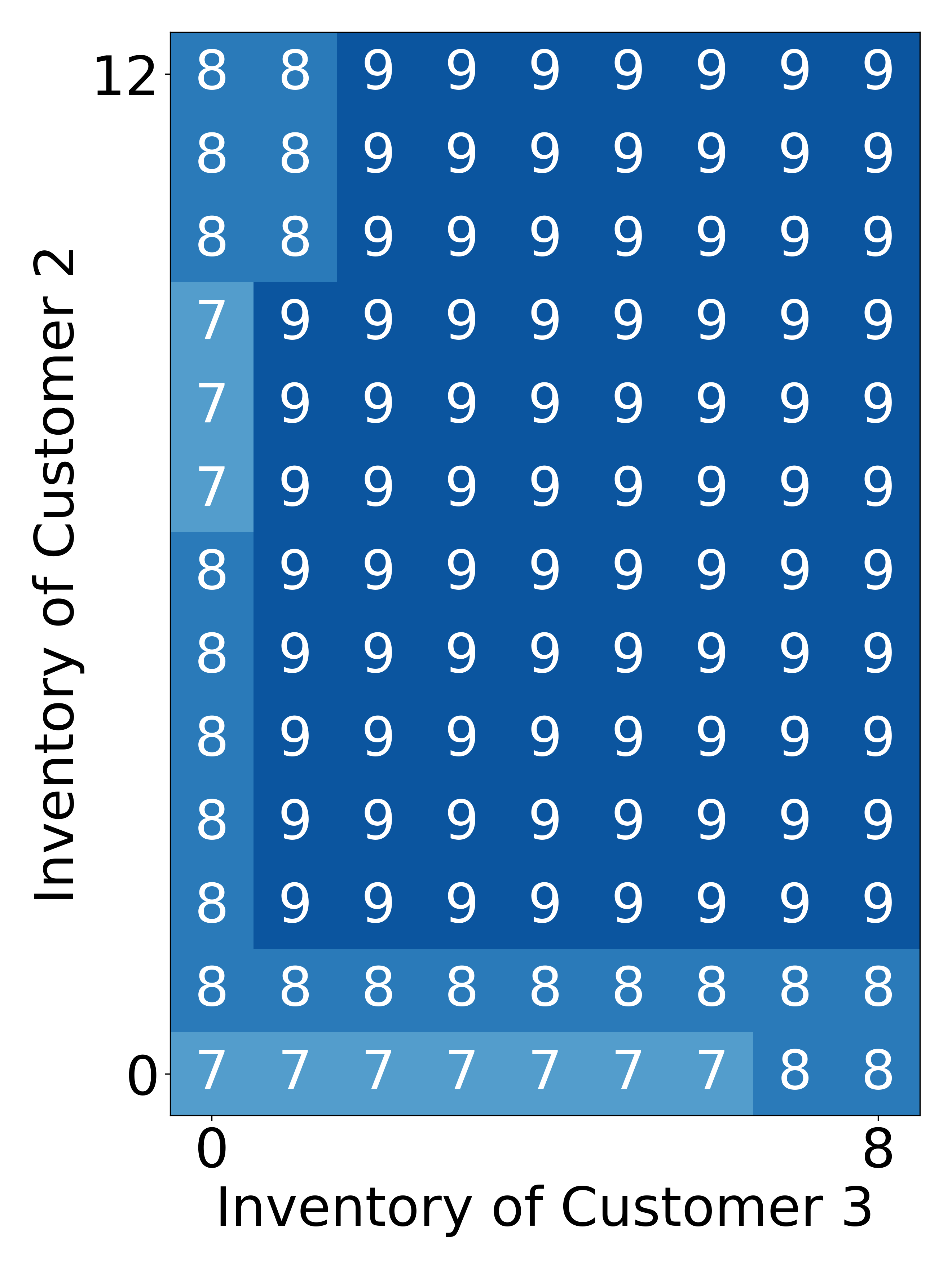}
        \caption{Customer 1 with LS-HD}
    \end{subfigure}
    \begin{subfigure}{0.225\textwidth}
        \centering
        \includegraphics[width=\textwidth]{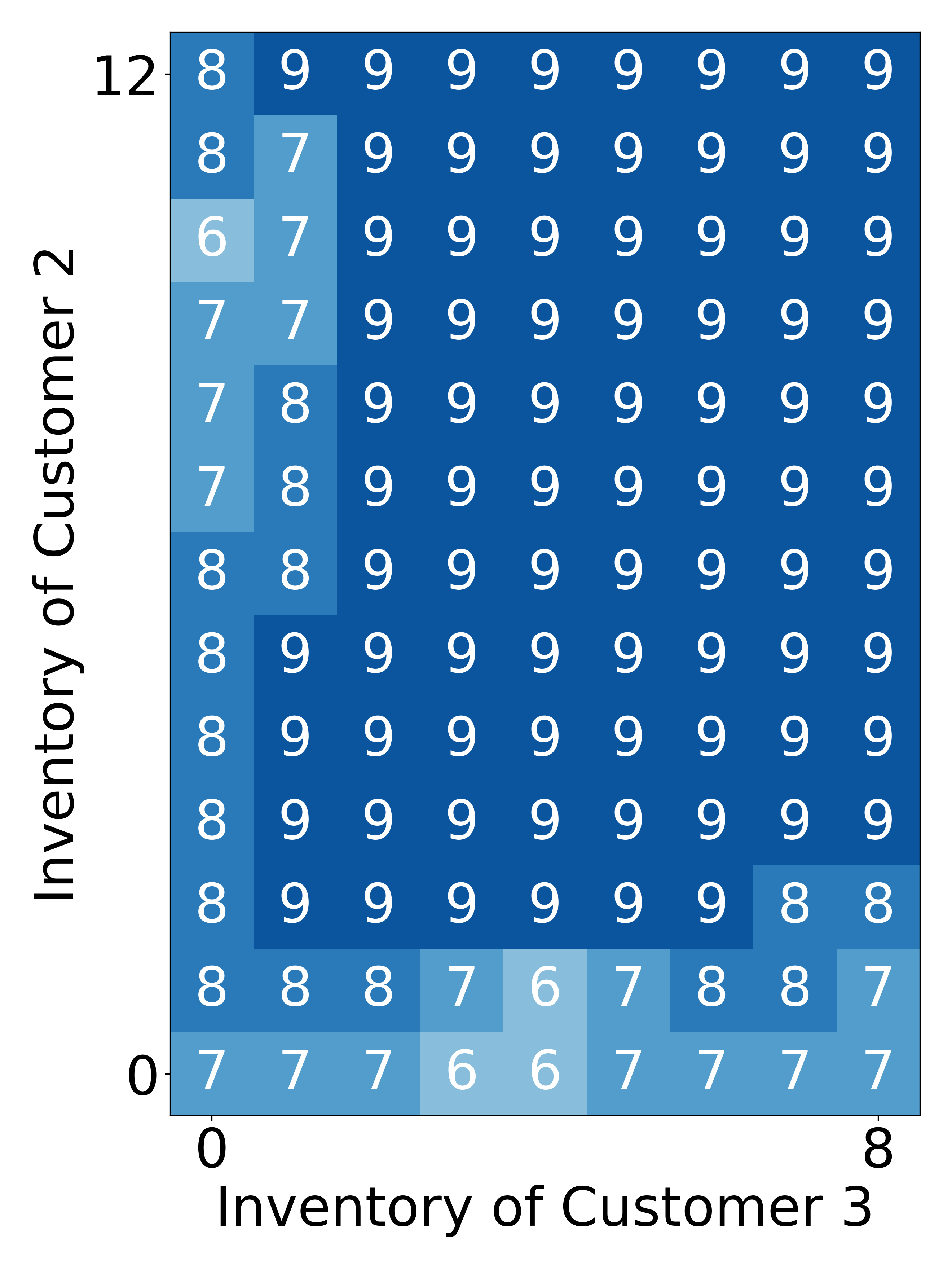}
        \caption{Customer 1 with HS-LD}
    \end{subfigure}
    \begin{subfigure}{0.225\textwidth}
        \centering
        \includegraphics[width=\textwidth]{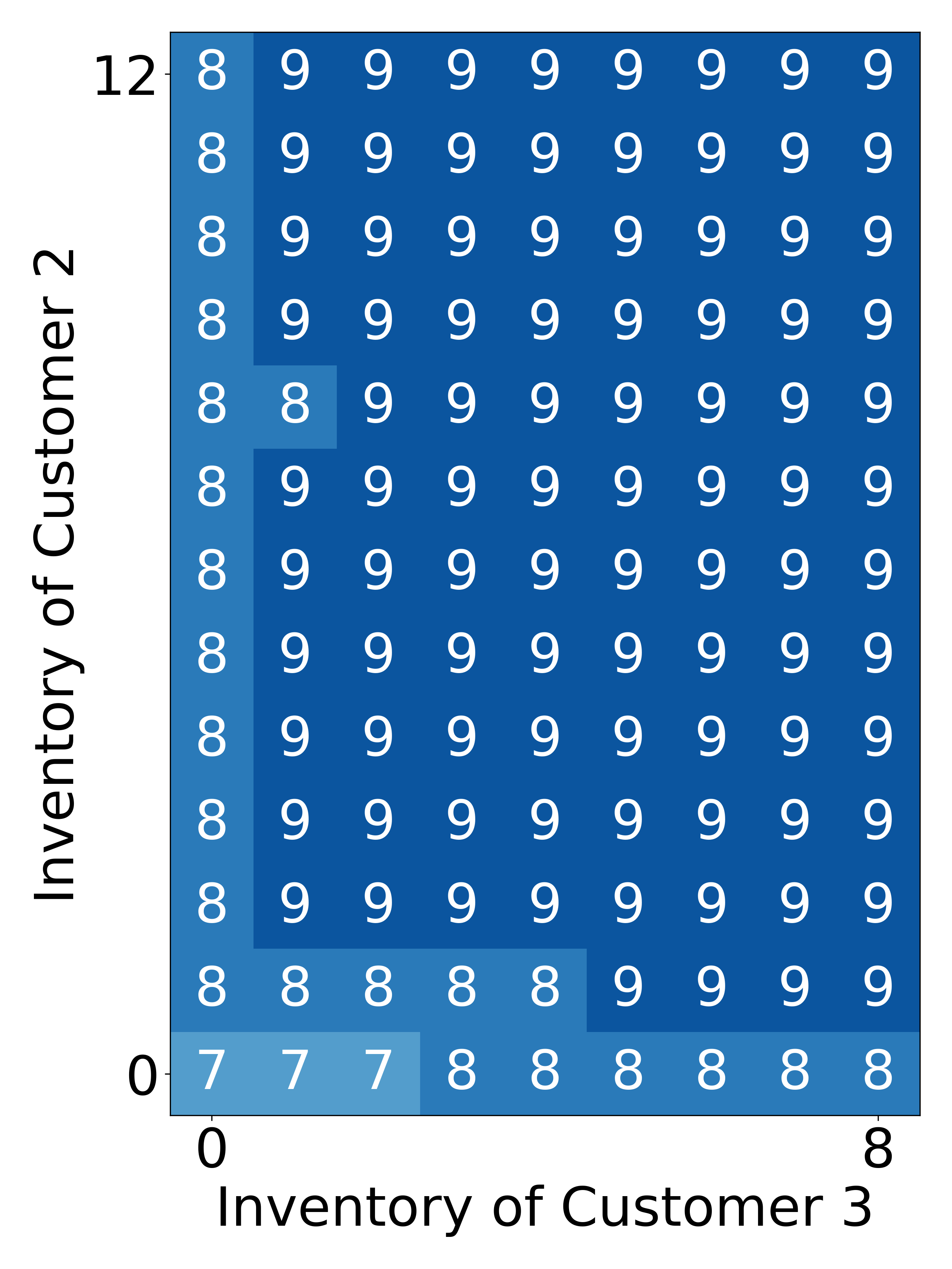}
        \caption{Customer 1 with HS-HD}
    \end{subfigure}
    
    \begin{subfigure}{0.225\textwidth}
        \centering
        \includegraphics[width=\textwidth]{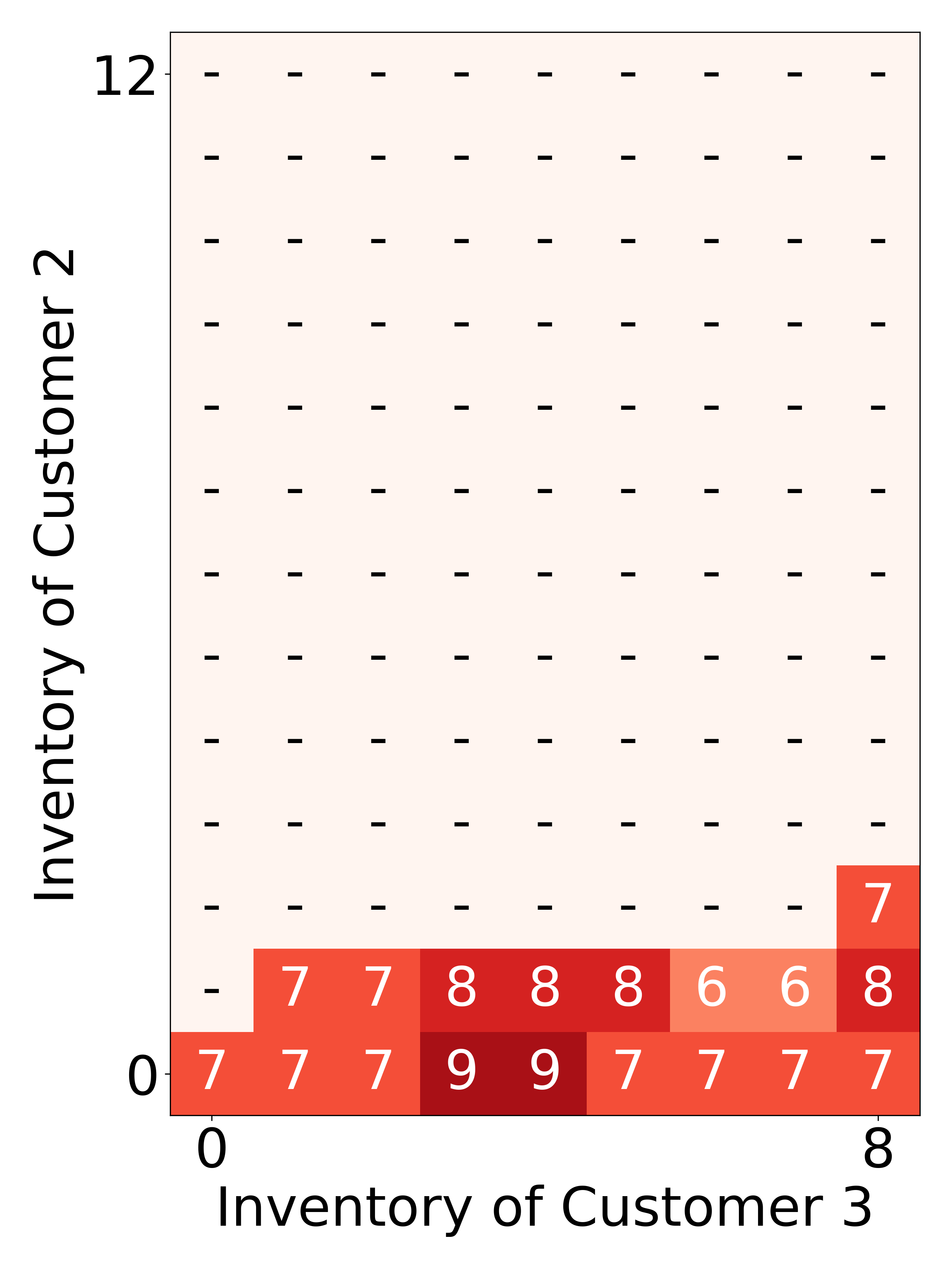}
        \caption{Customer 2 with LS-LD} \label{5e}
    \end{subfigure}
    \begin{subfigure}{0.225\textwidth}
        \centering
        \includegraphics[width=\textwidth]{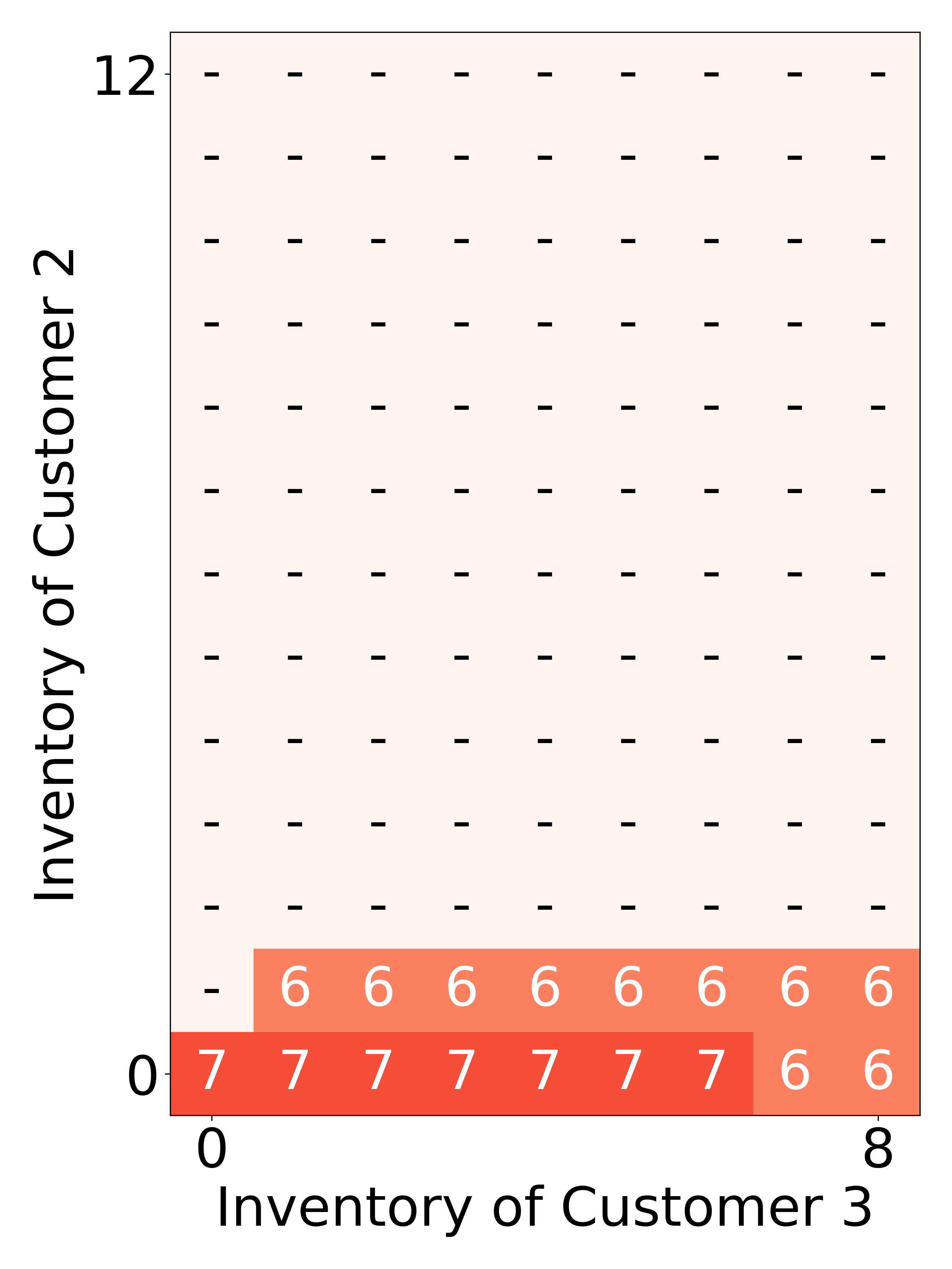}
        \caption{Customer 2 with LS-HD}
    \end{subfigure}
    \begin{subfigure}{0.225\textwidth}
        \centering
        \includegraphics[width=\textwidth]{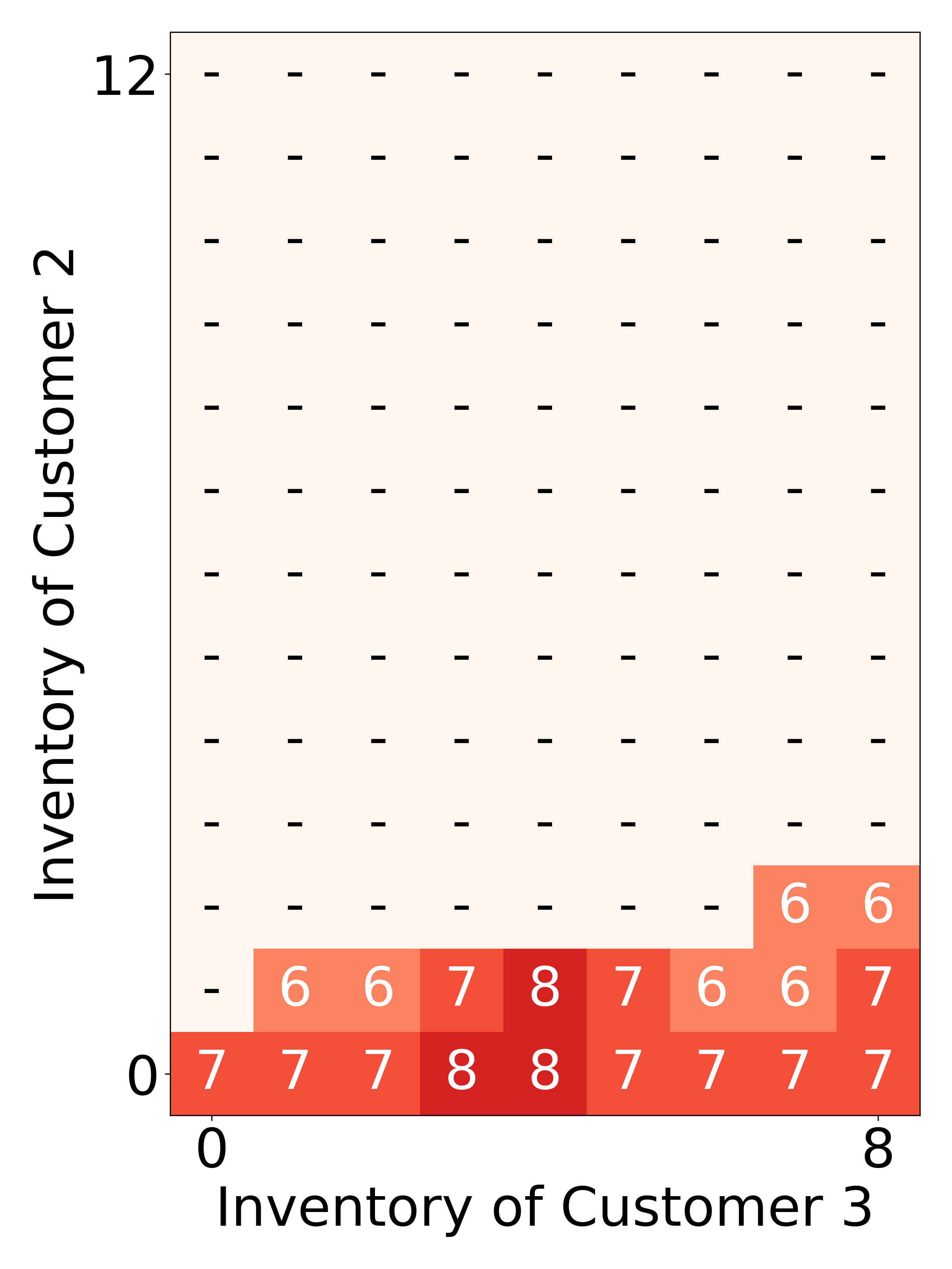}
        \caption{Customer 2 with HS-LD}
    \end{subfigure}
    \begin{subfigure}{0.225\textwidth}
        \centering
        \includegraphics[width=\textwidth]{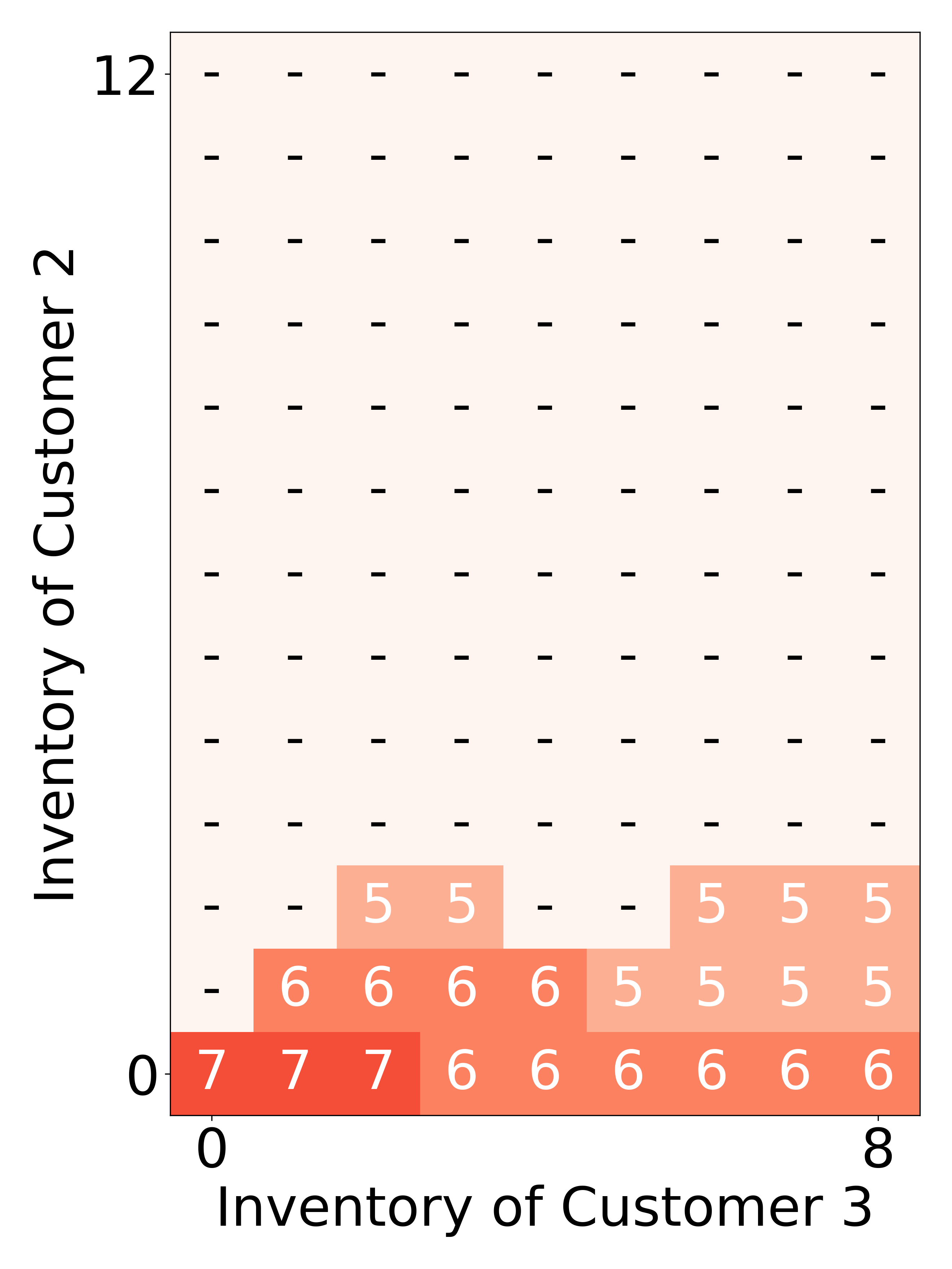}
        \caption{Customer 2 with HS-HD}  \label{5h}
    \end{subfigure}
    
    \begin{subfigure}{0.225\textwidth}
        \centering
        \includegraphics[width=\textwidth]{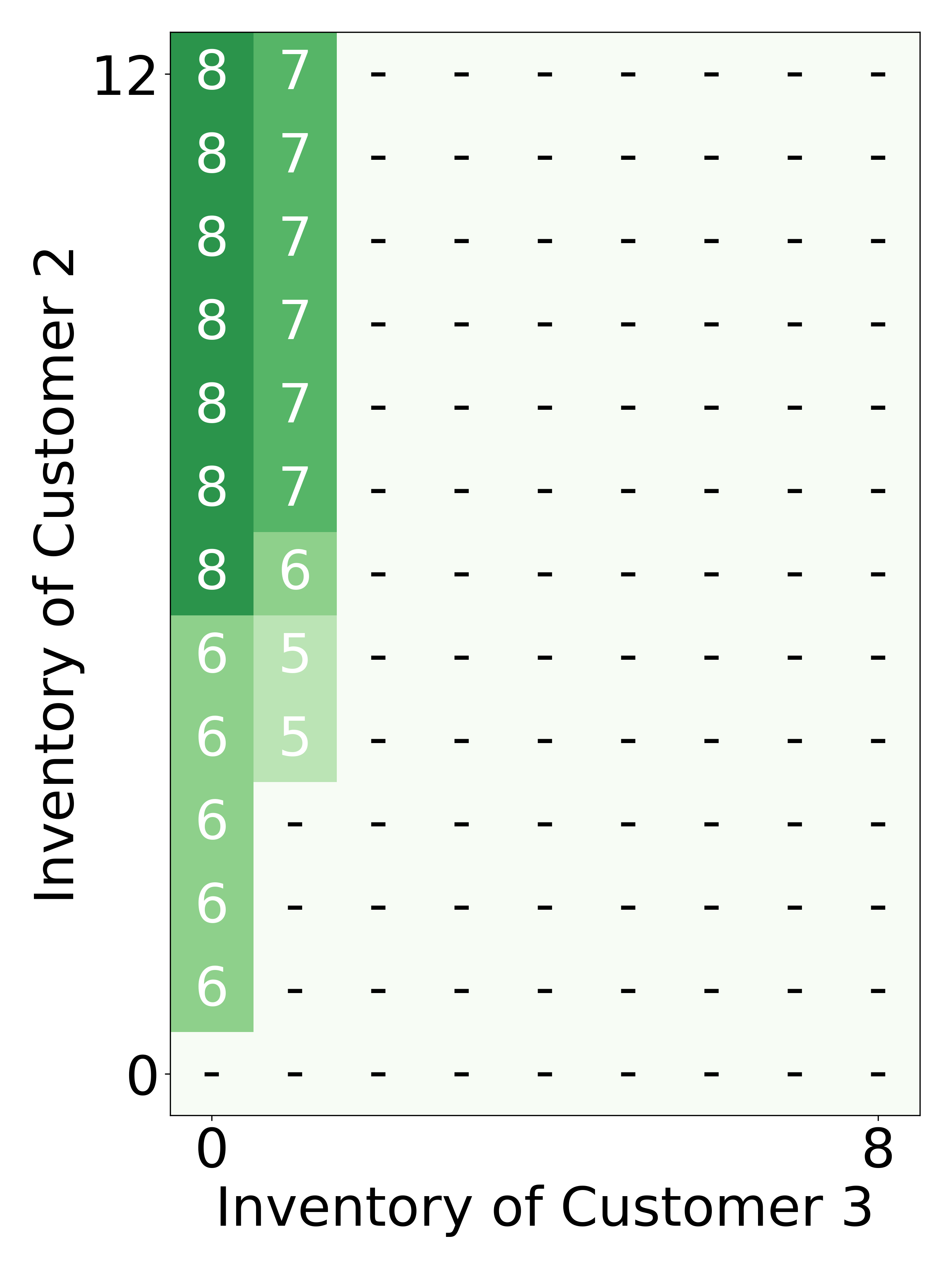}
        \caption{Customer 3 with LS-LD}
    \end{subfigure}
    \begin{subfigure}{0.225\textwidth}
        \centering
        \includegraphics[width=\textwidth]{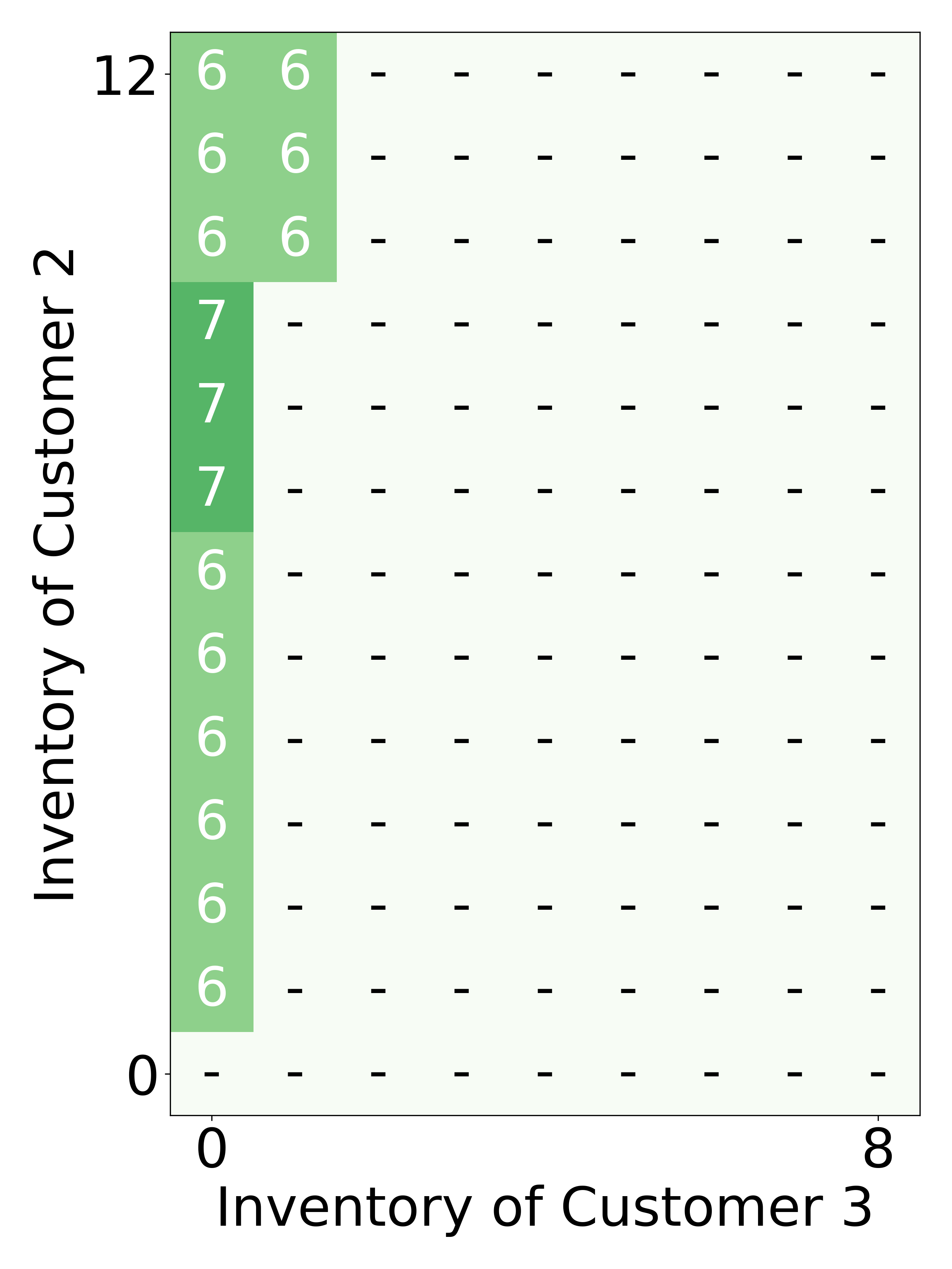}
        \caption{Customer 3 with LS-HD}
    \end{subfigure}
    \begin{subfigure}{0.225\textwidth}
        \centering
        \includegraphics[width=\textwidth]{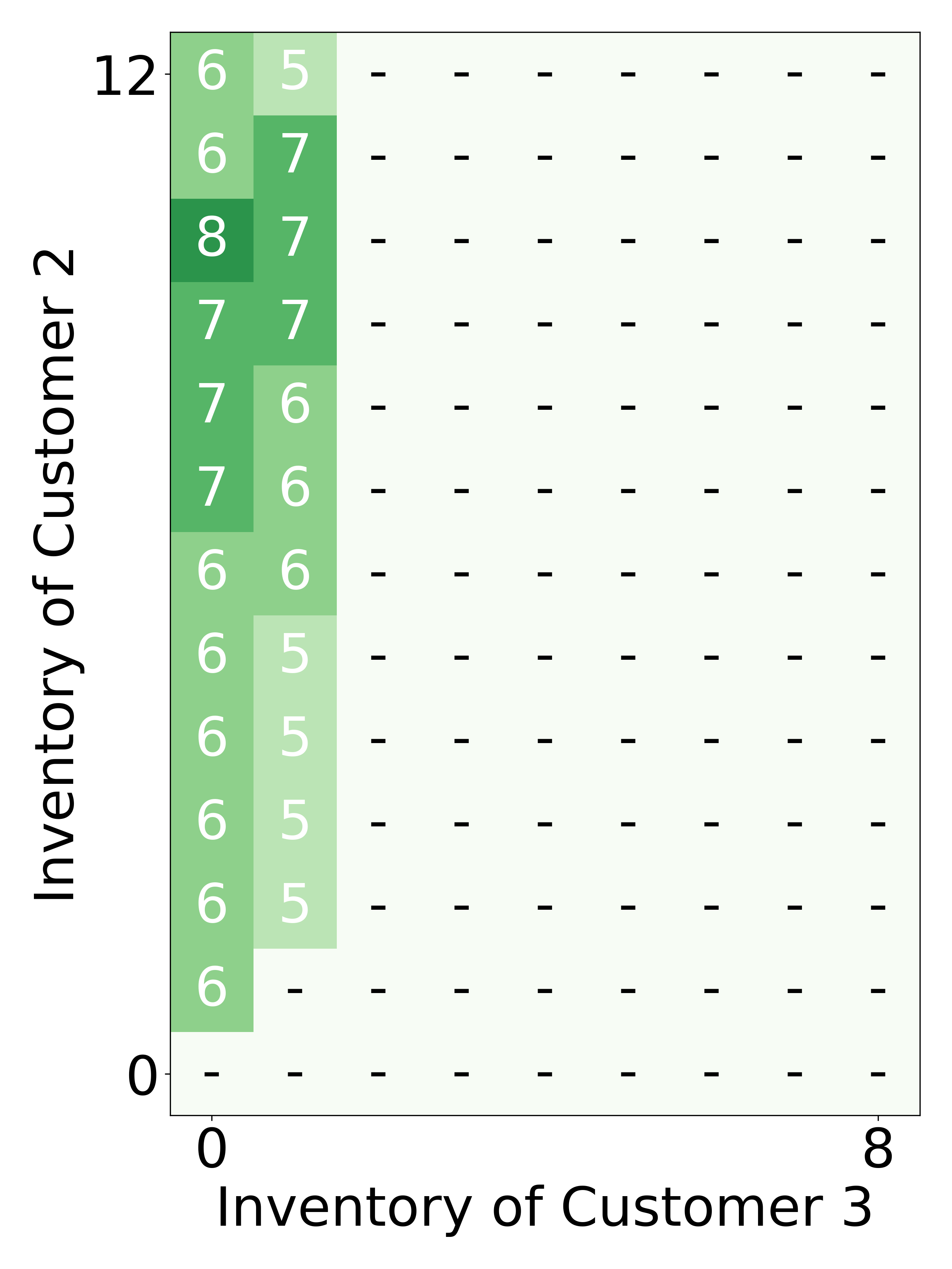}
        \caption{Customer 3 with HS-LD}
    \end{subfigure}
    \begin{subfigure}{0.225\textwidth}
        \centering
        \includegraphics[width=\textwidth]{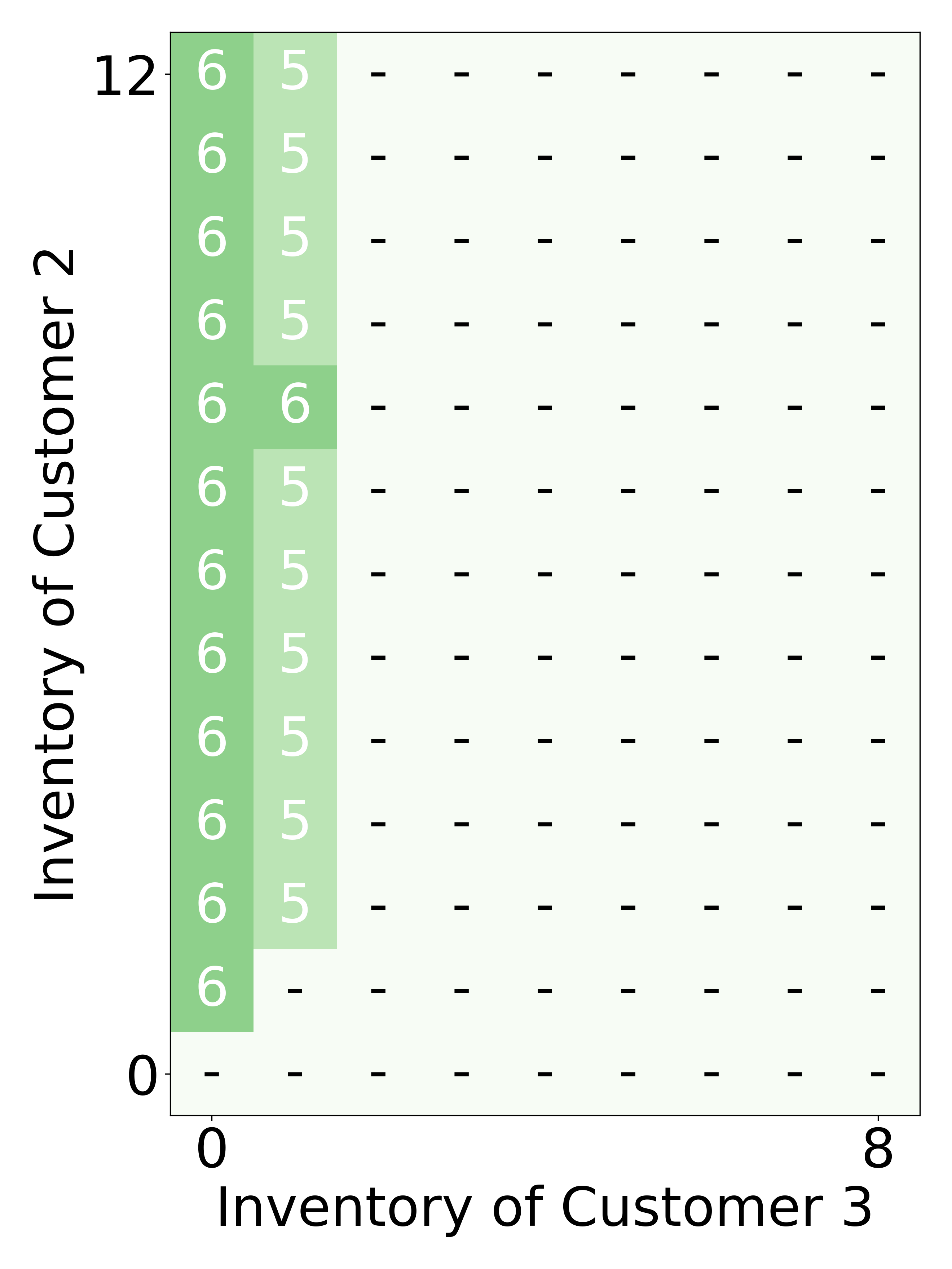}
        \caption{Customer 3 with HS-HD}
    \end{subfigure}
    
    \caption{Optimal replenishment policies at $(x_0, x_1) = (14,0)$ under varying combinations of demand and supply uncertainty. Rows are for customers, columns represent different uncertainty levels: LS, LD, HS, HD stands for Low/High uncertainty in Supply/Demand. }
    \label{figuncertainty}
\end{figure}

Additionally, we analyze the impact of demand and supply uncertainties on the optimal policy structure by solving four cases, each with a combination of low/high uncertainty in demand/supply, while holding all other parameters the same. Specifically, we adjust the parameters controlling the demand and supply uncertainty by factors of 0.5 and 1.5, respectively for low and high cases (see Appendix \ref{sec:instGen} for details). We adjust demand uncertainty across all customers, with no cases where demand uncertainty is increased for one customer while decreased for another. In Figure \ref{figuncertainty}, we provide the optimal policies for the states at $(x_0, x_1) = (14,0)$, comparing the values to Figure \ref{figoptpol2}. 
Sub-figure columns represent different uncertainty combinations. We observe that higher demand uncertainty increases replenishment to customer 1 while decreasing it for the other customers. This response is likely due to customer 1’s increased risk of lost sales relatively to the others, which raises the costs associated with under-replenishment. A similar, though less straightforward, pattern emerges with increased supply uncertainty. For example, for customer 3, while replenishments occur more frequently with higher supply uncertainty, their volumes decrease and this decreased amount is utilized again for customer 1.  Moreover, at $(x_2,x_3) = (2,8)$ in Figures \ref{5e} - \ref{5h}, we observe significant variation in customer 2's replenishment values, ranging within $[0,5,6,7]$, despite the system being in the same state; the only difference lies in the changes to stochasticity. These observations show that even in simplified, small-scale instances, stochasticity can sharply affect the policy structure that cannot be easily describable.

While these analyses provide insights for a simplified setup, predicting similar replenishment boundaries becomes infeasible in larger, practical scenarios. This is due to some challenges: (1) value iteration is not solvable with a large state space, and (2) complex interactions among multiple customers and vehicles lead to behaviors that are even more difficult to interpret when considering more customers simultaneously. In the next section, we analyze larger instances with up to 15 customers and six vehicles, where we compare the solution quality and computational time of our learning-based methods against benchmark algorithms.

\subsection{Comparison of benchmark algorithms on large instances} \label{sec:SSPO2}

This section presents an analysis of our constrained reinforcement learning approaches, CRL and LCRL, against the $(s,S)$-based heuristic and Power-of-Two (PO2) methods. We generate 10 instances for each combination of \(N \in \{9,12,15\}\) and \(q \in \{4,5,6\}\). The  results are presented in Table \ref{tab:crlsspo2}.

\renewcommand{\arraystretch}{1.25}
\setlength{\tabcolsep}{7pt}
\begin{table}[htp]
\centering
\caption{Comparison of benchmark algorithms on solution quality and time for instances on average (times in sec.)} \label{tab:crlsspo2}
\resizebox{0.85\textwidth}{!}{
\begin{tabular}{llrrrrrrrr} \toprule
    &     & \multicolumn{2}{c}{LCRL}  & \multicolumn{2}{c}{CRL}    & \multicolumn{2}{c}{$(s,S)$}        & \multicolumn{2}{c}{PO2}       \\ \cmidrule(r){3-4} \cmidrule(r){5-6} \cmidrule(r){7-8} \cmidrule(r){9-10}
$N$ & $q$ & $\text{time}_\text{train}$ & $\text{time}_\text{sim}$ & $\Delta(\%)$ & $\text{time}_\text{sim}$ & $\Delta(\%)$ & $\text{time}_\text{train}$ & $\Delta(\%)$ & $\text{time}_\text{train}$ \\ \midrule  
9  & 4    & 6,797  & 10.24 & 13 & 0.18 & 28 & 2,059 & 104 & 0.01 \\
   & 5    & 8,745  & 10.23 & 18 & 0.24 & 34 & 2,005 & 192 & 0.01 \\
   & 6    & 9,130  & 26.09 & 15 & 0.25 & 30 & 2,005 & 258 & 0.02 \\
12 & 4    & 9,713  & 18.12 & 7  & 0.24 & 16 & 2,629 & 50  & 0.01 \\
   & 5    & 17,333 & 19.57 & 13 & 0.44 & 22 & 2,734 & 88  & 0.02 \\
   & 6    & 19,239 & 17.95 & 16 & 0.51 & 26 & 2,739 & 146 & 0.02 \\
15 & 4    & 11,955 & 31.10 & 13 & 0.32 & 12 & 3,497 & 31  & 0.01 \\
   & 5    & 19,073 & 26.72 & 6  & 0.47 & 14 & 3,499 & 44  & 0.01 \\
   & 6    & 23,876 & 30.87 & 11 & 0.55 & 18 & 3,494 & 69  & 0.02 \\ \midrule
\multicolumn{2}{c}{avg.} & 13,985 & 21.21 & 13 & 0.36 & 22 & 2,740 & 109 & 0.01         \\ \bottomrule   
\end{tabular}
}
\end{table}

All time measurements are reported in seconds, averaged over 10 problem instances. The gaps are reported on average compared to LCRL. For LCRL, $\text{time}_\text{train}$ denotes the total duration spent executing Algorithm \ref{RLalgo}, which is also the training time for CRL (see Section \ref{sect:LCRL} for its discussion). For $(s,S)$, training time includes all pre-calculation and model solution times, excluding simulation to derive average costs per period. For PO2, it is the total time required to solve the model given in Eq. \eqref{po2start} to \eqref{po2end}, and then to derive a cyclic schedule for the simulation. We separate the time for deriving policies, $\text{time}_\text{train}$, from the simulation time used to estimate the policy's average costs per period. The simulation time per period, $\text{time}_\text{sim}$, is specified for CRL and LCRL only, as it is negligibly small for $(s,S)$ and PO2 ($\leq 10^{-5}$ seconds).

The analysis indicates that our methods outperform both benchmark algorithms across the tested scenarios. 
Furthermore, LCRL improves the solution of CRL by about $10\%$, though it requires significantly more solution time—around 60 times that of CRL—while utilizing 16 CPU cores for its computations. In fact, LCRL often reaches the 120-second time limit (see Appendix \ref{sec:instGen}), so we set $T=1$, as increasing this would exponentially raise the problem's complexity.
As we move from LCRL to CRL to $(s,S)$ and then to PO2, there is a reduction in both training and simulation times, but inversely, solution quality deteriorates, highlighting a trade-off between solution quality and computational time. However, CRL's performance decreases when the $q/N$ ratio is quite low, notably in cases with $N=15$ and $q=4$. This decrease is primarily due to the $\epsilon$-greedy training approach potentially missing the proper learning for some customers when the number of vehicles is limited. This occurred two out of $10$ instances on this set, resulting in slightly worse performance than $(s,S)$. 
However, even when these instances are included, CRL shows approximately  $10\%$ improvement over $(s,S)$ on average of all instances. This issue is not seen in LCRL, as its expanded short-term planning model has inherently calculations of holding and lost sales costs, preventing missing out on these customers. Its performance is consistently better than all other solution methods. Training time for $(s,S)$ seems independent of the number of vehicles $q$, as Algorithm \ref{SSalgo} has two for loops, both of which depend only on the customers and their inventory capacities. Finally, PO2 is significantly worse in solution quality as highly uncertain environments are hard to maintain with true cyclic approaches. However, its solution takes no time and policies are direct to find.

We further analyze the average cost components across solution methods in Table  \ref{tab:crlsspo2COST}. $ \textsc{t}$, $ \textsc{h}$, $ \textsc{l}$, and $\textsc{s}$  represent the costs of transportation, holding, lost sales, and sales, respectively. Values are presented as marginal percentage differences relative to LCRL. Results indicate that CRL performs worse than LCRL; despite comparable inventory holding, CRL uses more vehicles yet fails to reduce lost sales. This is likely due to its short-sighted planning, missing out on potential lost sales for the upcoming periods. In contrast, other benchmark methods show decreased transportation costs relative to LCRL. This is due to their individual customer-focused policies, which overlook vehicle capacity constraints in aggregate. This approach leads to missed deliveries and, thus, higher lost sales. This issue is seen more often when $q/N$ decreases, imposing stricter vehicle limits and missing more often. Fewer shipments, however, allow more inventory for sale at the supplier, increasing revenue, though this does not compensate for the overall cost increase. Moreover, with larger $N$, both benchmarks see reduced lost sales costs relative to learning methods due to demand aggregation, which reduces overall demand uncertainty. Nevertheless, learning-based methods achieve more balanced, cost-effective outcomes.

\renewcommand{\arraystretch}{1.25}
\setlength{\tabcolsep}{7pt}
\begin{table}[htp]
\centering
\caption{Comparison of benchmark algorithms on costs for instances on average (values in $\Delta\%$)} \label{tab:crlsspo2COST}
\resizebox{0.85\textwidth}{!}{
\begin{tabular}{llrrrrrrrrrrrr} \toprule
    &     & \multicolumn{4}{c}{CRL}    & \multicolumn{4}{c}{$(s,S)$}        & \multicolumn{4}{c}{PO2}     \\ \cmidrule(r){3-6} \cmidrule(r){7-10} \cmidrule(r){11-14} 
$N$ & $q$ & $\textsc{t}$ & $ \textsc{h}$ & $ \textsc{l}$ & $ \textsc{s}$ & $\textsc{t}$ & $ \textsc{h}$ & $ \textsc{l}$ & $ \textsc{s}$& $\textsc{t}$ & $ \textsc{h}$ & $ \textsc{l}$ & $ \textsc{s}$ \\ \midrule  
9  & \multicolumn{1}{r|}{4} & 32                   & -1 & 4   & \multicolumn{1}{r|}{2}  & -6  & -34 & 141 & \multicolumn{1}{r|}{131} & -27 & -25 & 422  & 390  \\
  & \multicolumn{1}{r|}{5} & 42                   & 5  & -4  & \multicolumn{1}{r|}{-6} & 8   & -36 & 154 & \multicolumn{1}{r|}{143} & -30 & -41 & 797  & 733  \\
  & \multicolumn{1}{r|}{6} & 40                   & 14 & -10 & \multicolumn{1}{r|}{-6} & 8   & -27 & 141 & \multicolumn{1}{r|}{138} & -38 & -49 & 1110 & 1057 \\
12 & \multicolumn{1}{r|}{4} & 14                   & -8 & 18  & \multicolumn{1}{r|}{11} & -29 & -27 & 126 & \multicolumn{1}{r|}{102} & -27 & -17 & 229  & 189  \\
& \multicolumn{1}{r|}{5} & 28                   & -5 & 14  & \multicolumn{1}{r|}{-4} & -15 & -36 & 156 & \multicolumn{1}{r|}{114} & -30 & -27 & 420  & 324  \\
 & \multicolumn{1}{r|}{6} & 40                   & 0  & 1   & \multicolumn{1}{r|}{-8} & -1  & -40 & 155 & \multicolumn{1}{r|}{125} & -31 & -37 & 655  & 547  \\
15 & \multicolumn{1}{r|}{4} & 6                    & -6 & 39  & \multicolumn{1}{r|}{29} & -35 & -13 & 82  & \multicolumn{1}{r|}{71}  & -16 & -9  & 115  & 100  \\
 & \multicolumn{1}{r|}{5} & 12                   & -6 & 13  & \multicolumn{1}{r|}{15} & -30 & -26 & 120 & \multicolumn{1}{r|}{112} & -27 & -15 & 208  & 192  \\
 & \multicolumn{1}{r|}{6} & 23                   & -4 & 12  & \multicolumn{1}{r|}{7}  & -21 & -33 & 138 & \multicolumn{1}{r|}{121}                    & -30 & -22 & 323  & 285  \\  \midrule \multicolumn{2}{c|}{avg.}                    & 26                   & -1 & 10  & \multicolumn{1}{r|}{5}                    & -13 & -30 & 135 & \multicolumn{1}{r|}{117}                      & -28 & -27 & 475  & 424  \\ \bottomrule        
\end{tabular}
}
\end{table}

\subsection{Sensitivity analysis of CRL learning parameters} \label{sec:sensitivity}

In this subsection, we analyze the impact of learning parameters on the performance of the CRL algorithm and compare these results with those from value iteration (VI). We use the instance of Section \ref{sec:opt}.

We first solve the instance using VI, obtaining post-decision state value predictions, $V_0$, from Algorithm \ref{VI}. These predictions are then used in a linear regression model as the dependent variable, and independent variables are selected as the same features used in CRL; $s$, $s^2$, $s^3$, and $\sqrt{s}$. The resulting regression coefficients are selected as VI’s feature weights. Simultaneously, CRL is run to derive its own feature weights. We use the weights obtained from both VI and CRL to predict the long-run profit for each inventory level in $[0,U]$, expressed as $- w^\intercal \psi(s)$. This approach allows us to compare the weights obtained from VI's regression with those learned by CRL. 

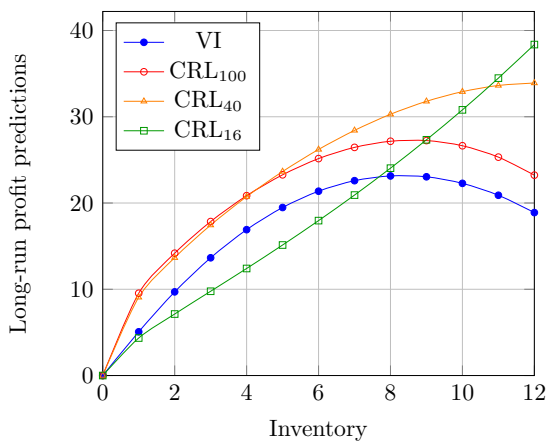
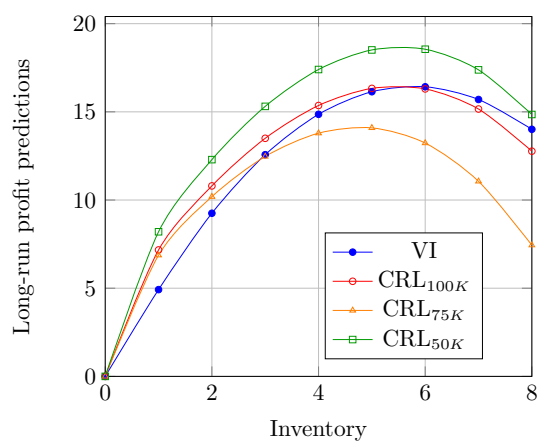
\begin{figure}[h]
    \centering
    \begin{subfigure}{0.46\textwidth}
        \centering
        \resizebox{\textwidth}{!}{
        \begin{tikzpicture}
    \begin{axis}[
        width=8cm, 
        height=7cm, 
        xlabel={Inventory}, 
        ylabel={Long-run profit predictions}, 
        legend pos=north west, 
        smooth, 
        grid=both, 
        xmin=0, xmax=12, 
        ymin=0, 
        legend style={font=\small}, 
        label style={font=\small}, 
        tick label style={font=\small} 
    ]
    \addplot[color=blue, mark=*, mark options={scale=0.7}] 
        coordinates {
        (0,0)(1,5.07)(2,9.70)(3,13.65)(4,16.91)(5,19.48)(6,21.37)(7,22.59)(8,23.14)(9,23.04)(10,22.28)(11,20.90)(12,18.89)
    };
    \addlegendentry{VI}

    \addplot[color=red, mark=o, mark options={scale=0.7}] 
        coordinates {
        (0,0)(1,9.54)(2,14.18)(3,17.85)(4,20.85)(5,23.27)(6,25.15)(7,26.45)(8,27.16)(9,27.24)(10,26.64)(11,25.32)(12,23.23)
    };
    \addlegendentry{$\text{CRL}_{100}$}

    \addplot[color=orange, mark=triangle, mark options={scale=0.7}] 
        coordinates {
        (0,0)(1,9.07)(2,13.66)(3,17.44)(4,20.74)(5,23.65)(6,26.21)(7,28.42)(8,30.29)(9,31.79)(10,32.91)(11,33.63)(12,33.91)
    };
    \addlegendentry{$\text{CRL}_{40}$}

    \addplot[color=green!60!black, mark=square, mark options={scale=0.7}] 
        coordinates {
        (0,0)(1,4.35)(2,7.13)(3,9.77)(4,12.41)(5,15.13)(6,17.96)(7,20.92)(8,24.04)(9,27.33)(10,30.80)(11,34.48)(12,38.37)
    };
    \addlegendentry{$\text{CRL}_{16}$}

    \end{axis}
\end{tikzpicture}
        }
        \caption{Effect of learning rate numerator; for $16$, $40$, and $100$.}
        \label{fig:alpha}
    \end{subfigure}
    \hspace*{0.05\textwidth} 
    \begin{subfigure}{0.45\textwidth}
        \centering
        \resizebox{\textwidth}{!}{
        \begin{tikzpicture}
    \begin{axis}[
        width=8cm, 
        height=7cm, 
        xlabel={Inventory}, 
        ylabel={Long-run profit predictions}, 
        smooth, 
        grid=both, 
        xmin=0, xmax=8, 
        ymin=0, 
        legend style={at={(0.7,0.4)}, anchor=north, font=\small}, 
        label style={font=\small}, 
        tick label style={font=\small} 
    ]
    \addplot[color=blue, mark=*, mark options={scale=0.7}] 
        coordinates {
        (0,0)(1,4.92)(2,9.25)(3,12.57)(4,14.86)(5,16.14)(6,16.42)(7,15.70)(8,14.01)
    };
    \addlegendentry{VI}

    \addplot[color=red, mark=o, mark options={scale=0.7}] 
        coordinates {
        (0,0)(1,7.17)(2,10.80)(3,13.50)(4,15.36)(5,16.33)(6,16.30)(7,15.16)(8,12.77)
    };
    \addlegendentry{$\text{CRL}_{100K}$}

    \addplot[color=orange, mark=triangle, mark options={scale=0.7}] 
        coordinates {
        (0,0)(1,6.87)(2,10.19)(3,12.47)(4,13.79)(5,14.09)(6,13.23)(7,11.06)(8,7.44)
    };
    \addlegendentry{$\text{CRL}_{75K}$}

    \addplot[color=green!60!black, mark=square, mark options={scale=0.7}] 
        coordinates {
        (0,0)(1,8.20)(2,12.29)(3,15.31)(4,17.40)(5,18.51)(6,18.55)(7,17.38)(8,14.85)
    };
    \addlegendentry{$\text{CRL}_{50K}$}

    \end{axis}
\end{tikzpicture}
        }
        \caption{Effect of training length; for 50K, 75K, and 100K periods.}
        \label{fig:period}
    \end{subfigure}
    \caption{Impact of learning parameters on long-run profit predictions of CRL against VI.}
    \label{fig:parameter}
\end{figure}

We first examine how different learning rates affect the function approximation. Initially, we employ a generalized harmonic learning rate formula \citep{powell2011approximate}, expressed as $\alpha_t = \frac{40}{5000 + t - 1}$ at the $t^{th}$ period of the algorithm. In addition to the numerator of 40, we test values of 16 and 100, referred to as $\text{CRL}_{16}$, $\text{CRL}_{40}$, and $\text{CRL}_{100}$. In Figure \ref{fig:alpha}, we show the long-run profit for varying inventory levels, where increasing the learning rate numerator from 16 to 100 noticeably improves approximation quality. Although and $\text{CRL}_{100}$ do not fully converge to the VI, the marginal changes across inventory levels are nearly identical, except for the zero-inventory case. Such behavior is characteristic of undiscounted value iteration, where the set of $V_0$ diverges in the long run, but only the marginal differences across states remain meaningful. 

In addition to analyzing the effects of the learning rate, we also examine how the number of training periods impacts cost prediction. CRL is normally run with 100K periods. We compare this with runs of 50K and 75K periods while maintaining the initial learning rates as $\alpha_t = \frac{40}{5000 + t - 1}$. Figure \ref{fig:period} illustrates the long-run profit predictions for this analysis, which is given for another customer of the same instance. We observe that while the results from 50K and 75K periods show some volatility, the model with 100K periods converges near to values that align with those obtained from VI. This behavior can be attributed to the inherent exploration-exploitation dynamics in RL, where the model initially wombles around local optima—exploring a range of suboptimal solutions—before converging. The near-convergence with 100K periods suggests that, with sufficient simulation time, CRL is able to approximate VI’s predictions more accurately. 

In Figure \ref{fig:alpha}, increasing the default learning rate numerator from $40$ to $100$ noticeably improves the approximation quality. However, this also leads to more complex optimization problems for Gurobi, increasing the solution time. Specifically, the training times for $\text{CRL}_{100}$, $\text{CRL}_{40}$, and $\text{CRL}_{16}$ were 8.7, 6.7, and 4.5 minutes, respectively. This pattern is typical of our algorithm, where higher learning rates improve weight prediction accuracy but also leads to longer solution times. Besides, the same numerator appears to be sufficient for another customer in Figure \ref{fig:period}. Based on preliminary experiments, we select the numerator value of $40$ as it achieves a balance solution quality and computation time. 

While both figures show similar convex function approximations, we often observe a linear increase in long-run profit predictions occurs with inventory for suppliers. Additionally, for customers with high inventory capacities relative to vehicle capacities, the weight predictions tend to show piece-wise convex approximations, which appear more often in mid-sized instances as in Section \ref{sec:SSPO2}. These suggest CRL’s flexibility in capturing varying inventory dynamics in capacities and constraints.

\subsection{Sensitivity analysis of algorithm components}
\label{sec:sensitivity_components}

This section assesses the sensitivity of the proposed algorithms to three key components: the basis functions of the value function approximation, the exploration decay rate, and the lookahead horizon.

First, we analyze the impact of the basis functions $\psi(s)$ used in the post-decision value function approximation. In the base case, $\psi(s)$ comprises linear combinations of $s_i, s_i^2, \sqrt{s_i}$, and $s_i^3$. To test robustness, we re-ran CRL on the test set of instances with $N=9$ and $q=5$ using alternative specifications where these terms are selectively enabled if $1$, and $0$ otherwise. Table~\ref{tab:vfa_robustness} reports the average percentage deviation in objective value relative to the best specification, along with the average training time in seconds.

\begin{table}[h]
\centering
\caption{Effect of basis functions in the value function approximation for $N=9, q=5$ (average of 10 instances).}
\label{tab:vfa_robustness}
\begin{tabular}{ccccrr}
\toprule
$s$ & $s^2$ & $\sqrt{s}$ & $s^3$ & Costs ($\Delta\%$) & Time (sec.) \\
\midrule
1 & 1 & 1 & 1 &   0.0 &  9,141 \\
1 & 0 & 1 & 1 &   1.6 &  6,427 \\
0 & 1 & 1 & 1 &   2.5 & 10,723 \\
1 & 1 & 1 & 0 & {15.3$^{*}$} &  2,848 \\
1 & 1 & 0 & 1 & {20.5$^{*}$} &  2,310 \\
1 & 0 & 0 & 0 & {38.8$^{*}$} &   216 \\
1 & 1 & 0 & 0 & {41.8$^{*}$} &   396 \\
0 & 1 & 0 & 0 & {89.3$^{*}$} &   650 \\
\bottomrule
\end{tabular}
\end{table}

The results indicate that the full specification $(s_i, s_i^2, \sqrt{s_i}, s_i^3)$ yields the lowest average cost. For the asterisk-marked specifications, the cost increase is observed for all 10 instances. On the other hand, the differences for omitting only $s_i^2$ or only $s_i$, which are on average $1.6\%$ and $2.5\%$ respectively, are not statistically significant under $p=0.05$. A paired two-sided Wilcoxon signed-rank test over the 10 instances indicates $T=24$, $p=0.770$ and $T=11$, $p=0.105$ for these settings, respectively. Specifications omitting $\sqrt{s_i}$ or $s_i^3$, however, yield substantially larger deviations, confirming that a richer basis is required to capture the complex state-value relationships in the DIRP. While simpler specifications reduce training time by simplifying the underlying mixed-integer quadratic programs, they do so at the expense of solution quality.

Second, we evaluate the effect of the exploration decay rate $\epsilon_t$ on the performance of CRL. We compare the base setting ($\epsilon_t = 0.999983^t$) against a faster decay ($\epsilon_t = 0.999966^t$) and a slower decay ($\epsilon_t = 0.999991^t$) on the same set of instances. We observed that the base setting achieves the best solution quality. The faster decay rate increases the average cost by $5.3\%$, which is likely due to premature exploitation. Furthermore, because the algorithm selects the best action (requiring a MIQP solution) more frequently with lower $\epsilon$ values, the training time increases by approximately $33\%$ (from 9,142 to 12,126 seconds). Conversely, the slower decay rate reduces the training time by $47\%$ (to 4,827 seconds) by favoring random action selection, but results in a slightly higher average cost of $0.3\%$. We therefore maintain the base setting to prioritize solution quality.

Finally, regarding the lookahead horizon $T$ in LCRL, we test larger values compared to the base $T=1$ while using the feature weights trained via CRL. We observe that extending the simulation horizon beyond $T=1$ does not yield performance benefits. In experiments on instances with $N=9$ and $q=5$, increasing $T$ to 2 and 3 increased the average cost by $5.3\%$ and $8.6\%$, respectively, while increasing the average simulation time per period from $10.23$~sec. ($T=1$) to $17.77$~sec. ($T=2$) and $20.71$~sec. ($T=3$). This degradation is likely attributable to the disparity between the value function (trained on 1-step transitions via CRL) and the multi-period lookahead. The solution quality dropped further when testing more complex instances ($N=12$, $q=4$), likely due to the amplified misalignment in larger state and action spaces. Ideally, the model would be trained directly using the LCRL formulation. However, even at $T=1$, this is computationally intractable; training an instance at $N=9$ and $q=5$ is estimated to require almost a full week with 16 CPU cores, compared to a few hours for the CRL training with a single core. Consequently, we restrict our analysis to the CRL-trained weights and simulate with LCRL only at $T=1$.

\section{Conclusions} \label{conclusion}

In this study, we address the dynamic inventory routing problem (DIRP) with uncertainty of both supply and demand, motivated by the green hydrogen logistic challenges identified in the HEAVENN initiative. Unlike existing approaches that often assume infinite supply, our model explicitly incorporates both limited and stochastic supply, creating a hard-coupled DIRP where customer replenishment decisions are interdependent. To solve this problem, we develop a Constrained Reinforcement Learning (CRL) method that integrates a Mixed-Integer Quadratic Programming (MIQP) formulation for action selection directly into the learning process, ensuring state-dependent feasibility. We further enhance this approach with Lookahead-based CRL (LCRL), which embeds our CRL method within an n-step lookahead policy. This hybrid structure improves solution quality by using the lookahead model to refine the current action. While CRL itself is strong, LCRL offers additional benefits in handling complex planning scenarios, emphasizing the adaptability and strength of our proposed frameworks.

Our computational experiments demonstrate the efficacy of the proposed CRL and LCRL methods in solving the DIRP. We show that uncertainty in both supply and demand significantly influences optimal replenishment policies across all customers, highlighting the importance of considering these factors. For small-scale instances, where optimal solutions are obtained via value iteration, our methods achieve near-optimal solutions to optimality. In larger, realistic-sized instances, where value iteration becomes computationally intractable, CRL and LCRL outperform deep learning techniques, such as PPO, and heuristic approaches commonly implemented in inventory management, including $(s,S)$-policy based and Power-of-Two based replenishment policies, demonstrating their scalability and effectiveness. While LCRL further improves the solution of CRL, CRL provides solutions in significantly shorter time, showing a trade-off between proposed methods. 

While the proposed model provides a focused framework for early-stage distribution, its scope is subject to constraints that limit the immediate general applicability of the current contribution. Future research can extend our DIRP model in several directions to address these limitations. First, incorporating a vehicle routing system in which routes consider more than one customer would provide a more comprehensive approach to some logistics, addressing both routing and replenishment challenges simultaneously. Second, multi-supplier settings could be introduced, enabling an exploration of potential global cost reductions through supplier collaboration, and customer decomposition to the available suppliers. Third, heterogeneous vehicle fleets, with variations in capacity and possibly in costs, could be considered to reflect other practical settings and trade-offs in transportation planning. Fourth, to further reflect the specific physical and infrastructural context of green hydrogen, future models could incorporate constraints such as hydrogen inventory decay, vehicle emission limits, or dependencies on pipeline development phases. Finally, while this study assumes static supply and demand random distributions over the infinite planning horizon, future work could model dynamic random distributions over a finite horizon setting, capturing the evolving nature of supply chains as economies grow or contract over time.

Future research could also explore methodological extensions of CRL and LCRL to other MDP problems with state-dependent constraints. The proposed algorithms are generically designed and can be directly applied across other MDPs with state-dependent constraints, although the definition of the feasible action set $A(x)$, which constrains actions through problem-specific MIP regions, must indeed be explicitly tailored for each problem. This could be either one-step lookahead as in CRL, or a finite-horizon, finite-scenario model could be integrated as in LCRL to improve longer term planning and capture more realistic decision-making processes. Furthermore, LCRL relies on CRL-trained parameters in our computational experiments, to mitigate computational complexity. Future work could investigate training DIRP directly with LCRL to incorporate long-term planning during the learning phase. Although computationally demanding, such an approach holds significant potential to enhance solution quality by aligning the value function training with the lookahead policy, particularly for horizons at $T > 1$.

\ACKNOWLEDGMENT{This project has received funding from the Fuel Cells and Hydrogen 2 Joint Undertaking (now Clean Hydrogen Partnership) under Grant Agreement No 875090. This Joint Undertaking receives support from the European Union's Horizon 2020 research and innovation programme, Hydrogen Europe and Hydrogen Europe Research. Albert H. Schrotenboer has received support from the Dutch Science Foundation (NWO) through grant VI.Veni.211E.043. We thank the Center for Information Technology of the University of Groningen for their support and for providing access to the Hábrók high performance computing cluster.}











\newpage
\bibliographystyle{informs2014trsc}
\bibliography{main}

\newpage
\begin{APPENDICES} 

\section{Deep Reinforcement Learning} \label{sec:DRL}
To evaluate the performance of our methodology, we implement a Deep Reinforcement Learning (DRL) benchmark based on the Proximal Policy Optimization (PPO), which has become the default policy gradient techniques \citep{schulman2017proximal}. 

This DRL implementation operates in a continuous action space, contrasting with the discrete action space of CRL. This adaptation is necessary as we have exponentially many feasible actions in even modest-sized problems, making  traditional discrete action evaluations on deep learning intractable, such as policy gradients. In our DRL, continuous action space is defined as $a \coloneqq \{a_0, a_1,\ldots,a_N\}$, which is normalized relative to the inventory capacities of their respective locations. A similar normalization is done on states according to the corresponding location's inventory capacity. 

The actor neural network generates a mean and a standard deviation for each action $a_i$, utilizing shared neurons and hidden layers across the network. We duplicate only in the final layer where mean and standard deviation are separately generated. The outputs for both the mean and the standard deviation are transformed via the sigmoid function to ensure they remain close within the unit interval $[0,1]$. We control the standard deviation by decreasing it during training through a scaling factor, maintaining a balanced exploration-exploitation. At each training period, actions are stochastically selected based on Gaussian distributions defined by the neural network's outputs for mean and standard deviation. The resulting action, then, implies $a_0 U_0$ units to be sold, with each customer receiving $a_i U_i$ units of product.

Addressing constraints in Eq. \eqref{totAction} - \eqref{MDP-end}, the DIRP with DRL involves a different approach than CRL. Unlike CRL, which utilizes MIP region of $A(x)$ for action-selection process, the neural networks does not implicitly have such hard constraints. In order to manage constraint violations, DRL applies a Lagrangian relaxation method which is a common approach for DRL literature on constrained MDPs \citep{miryoosefi2019reinforcement,li2021augmented}. For instance, excess demands that exceed inventory are penalized per unit, see Eq. \eqref{totAction}. Any dispatched quantity exceeding customer inventory capacities is considered lost without cost, see Eq. \eqref{capCust}. The number of vehicles needed for the corresponding action is always selected as the minimum number of vehicle possible, so Eq. \eqref{capTruck} is satisfied. However, the total number of required vehicles for the action may be more than available, which is penalized with a fixed cost of emergent shipment per vehicle, see Eq. \eqref{limTruck}.

Training updates are performed in batches with traditional PPO clipping. We employ Generalized Advantage Estimation (GAE) to mirror the eligibility traces used in our CRL approach. The actor-loss and critic-loss are computed where an entropy bonus is added on actor-loss to encourage exploration and to prevent premature policy convergence. Each batch iteration adjusts the standard deviation multiplier, reducing exploration as the algorithm converges towards a policy. Moreover, the learning rate is gradually decreased towards the end of the training phase to refine the policy updates. Upon completion of the training phase, corresponding neural networks are tested on a simulation to estimate a long run total cost of the policy. The parameters of our DRL approach are summarized in Table \ref{tab:drl_parameters}.

\renewcommand{\arraystretch}{1}
\begin{table}[h]
\centering
\caption{Deep Reinforcement Learning Parameters}
\begin{tabular}{lll}
\hline
\textbf{Parameter} & \textbf{Description} & \textbf{Value} \\ \hline
Optimizer & Optimizer used for neural networks & Adam \\
Learning rate & Learning rates for actor and critic networks & 1e-4, 5e-5 \\
Learning rate decay & Decay rate of learning rates per 1000 periods & 0.997 \\
Hidden layers & Number of neurons of both networks & $\{256, 128, 128, 128\}$ \\
Activation & Activation function used in networks & ReLU \\
$N$ & Number of actions & 20 \\
$M$ & Number of periods per rollout & 200 \\
Mini-batch size & Number of samples per mini-batch & 200 \\
$T_{\text{train}}$  & Number of training periods & 1M \\
Std. deviation multiplier & Initial multiplier for exploration & 0.25 \\
Std. deviation decay  & Decay rate of std. multiplier at each batch & 0.997 \\
Entropy Bonus & Weight of entropy in the actor loss calculation & 0.01 \\
Discount factor & Discount factor for future rewards & 1 (No discount) \\
$\lambda$ & $\lambda$ for Generalized Advantage Estimation & 0.8 \\
Clipping rate & Clipping rate for PPO & 0.15 \\
Epoch & Number of passes over the batch per update & 5 \\
$T_{\text{sim}}$ & Number of simulation periods & 100K \\
$L$ & Length of warm-up on simulation & 1K \\
\hline
\end{tabular}
\label{tab:drl_parameters}
\end{table}

As it has recently started to gain attention in operations research literature \citep{greif2024combinatorial, harsha2025deep}, DRL methods face significant challenges when handling both complex state and action spaces. We observe the same, that the PPO setting underperforms in our problem context, and this is likely due to several factors. In traditional PPO, probability weights are assigned to each feasible action, with separate calculations made for each one. In our case, however, actions are not explicitly listed but implicitly defined by constraints. 
Even if we could list them, the sheer size of the action space makes distinct calculations or neuron assignments for each action computationally intractable. While \citet{kaynov2023deep} successfully assign neurons to each replenishment quantity, their problem is smaller in scale, both in terms of the number of customers and inventory capacities. To address our vast action space, we opted for a continuous action space normalized to inventory capacities in PPO, though continuous approaches are often less stable and harder to train \citep{dehaybe2024deep}.

Moreover, neural networks cannot handle constraints like limited vehicle availability, often producing infeasible solutions. We implemented a Lagrangian penalty for this purpose, but we further tested projection to the nearest feasible solution \citep{stranieri2024performance}, which did not deliver satisfactory results either. Projections often led to divergence in trained weights, and Lagrangian methods complicate the learning process by introducing additional terms in the reward structure. Ultimately, the implicit action definitions, along with the vast and constrained action space, likely caused PPO’s poor performance. We further explored alternative DRL methods, including REINFORCE, TD($\lambda$), and A3C, but none yielded satisfactory results. 

A natural question is whether a DRL policy can simply be restricted to output only feasible actions, thereby avoiding infeasibility. While this is possible for simple, independent constraints, it is computationally difficult for problems with coupled constraints, such as our DIRP. For example, whether a suggested action is feasible depends on jointly allocating a limited vehicle fleet and finite supplier inventory. This system-wide feasibility cannot be enforced by independent clipping or normalization of per-customer quantities. Indeed, this fundamental challenge is why common approaches for such MDPs often rely on indirect methods, such as Lagrangian relaxation \citep{miryoosefi2019reinforcement,li2021augmented}, applying post-hoc allocation rules to infeasible actions \citep{kaynov2023deep, stranieri2024performance}, or solving a mathematical program during each learning step \citep{rivera2017anticipatory, al2020approximate}. More generally, with greater system complexity, it is practical to cast the per-state decision in the underlying dynamic program as a mathematical program, as in the DIRP. 

Table \ref{tab:drl_comparison} summarizes the differences between standard DRL methods and CRL/LCRL. In particular, standard DRL methods typically address constraints at the policy level, whereas CRL/LCRL enforce feasibility at each decision step through MIP-based action selection.
\renewcommand{\arraystretch}{1}
\begin{table}[h]
\centering
 
\caption{Comparison of representative DRL methods and CRL/LCRL across key design dimensions.}
\label{tab:drl_comparison}
\footnotesize
\setlength{\tabcolsep}{4pt}
\begin{tabular}{l
                >{\raggedright\arraybackslash}p{2.1cm}
                >{\raggedright\arraybackslash}p{3.3cm}
                >{\raggedright\arraybackslash}p{1.8cm}
                >{\raggedright\arraybackslash}p{5.6cm}}
\toprule
\textbf{Method} &
\textbf{Action space} &
\textbf{Constraint handling} &
\textbf{Feasibility guarantee} &
\textbf{Applicable setting} \\
\midrule
PPO / TRPO &
Continuous or discrete &
Lagrangian / penalty-based &
Policy-level only$^{\dagger}$ &
On-policy control with differentiable objectives where safety can be managed as an expected long-run cost. \\
SAC &
Continuous &
Lagrangian / penalty-based &
Policy-level only$^{\dagger}$ &
Off-policy continuous control problems where high sample efficiency is required due to limited interaction data.\\
CQL &
Continuous or discrete &
Lagrangian / penalty-based &
Policy-level only$^{\dagger}$ &
Offline problems where online exploration is impractical (e.g., cost); learning from static logged data. \\
CRL / LCRL &
Mixed-integer &
Exact enforcement via MIP at each action; LCRL adds lookahead &
Yes &
Online/offline decision problems with dynamic action sets and hard, state-dependent constraints. \\
\bottomrule
\multicolumn{5}{>{\raggedright\arraybackslash}p{0.98\linewidth}}{\scriptsize
$^{\dagger}$ Constraints are handled at the policy level, typically as expected cumulative cost constraints in the CMDP sense, and do not guarantee feasibility of each action at each decision step.} \\
\end{tabular}
\end{table}

Nevertheless, we observe that PPO performs better in cases where supply and vehicle fleet are abundant. For example, in the $N=8$, $q=24$ case with supply set to three times the expected demand, PPO yields about $15\%$ better solutions than CRL. In these cases, the underlying MDP becomes weakly coupled per customer, and the probability of the neural network generating infeasible actions sharply decreases, resulting in improved performance. This supports our argument that the poor performance in other instances is driven by the infeasible solutions returned by neural networks. However, since these instances are not the focus of this study, and as all DRL methods consistently performed worse than other benchmarks for our instances, we excluded them from our computational experiments.

\section{Benchmark instances and computational parameters} \label{sec:instGen}

For each instance, we generate a graph where locations are uniformly randomized within a $10 \times 10$ square area. Demand at each point follows a normal distribution, with discrete means uniformly distributed in $[6,12]$, and standard deviations between $[25\%,75\%]$ of the mean. The given number of vehicles, $q$, influences vehicle capacity, which is selected as $C \coloneqq 2  \mathbb{E}\left[ \sum_{i \in \mathscr{N}} \phi_i\right] / q$, where $\mathbb{E}\left[ \sum_{i \in \mathscr{N}} \phi_i\right]$ is the expected total demand per period. Supply follows a normal distribution, where $\mathbb{E}\left[ \phi_0 \right] \coloneqq \mathbb{E}\left[ \sum_{i \in \mathscr{N}} \phi_i\right]$. Its standard deviation is set to $60\%$ of the mean. The inventory capacities, $U_i$ are selected as $10$ times of the expected demand per customer, and $2.5$ times of the expected supply for supplier. For discrete support of supply and demand realizations, $\Phi$, we round probabilities to the nearest integer within the range $\left[ \max (0, \mathbb{E}\left[ \phi_i \right] - 3 \sigma_i), \min (U_i, \mathbb{E}\left[ \phi_i \right] + 3 \sigma_i) \right]$, where $\sigma_i$ denotes the standard deviation associated with location $i \in \mathscr{N}^0$. Finally, we select a cost vector of $(W,w,h_s,h_c,\ell,\rho) = (15,1.5,0.1,0.2,30,2.5)$.

For CRL, we train Algorithm \ref{RLalgo} over $100K$ periods. We set forgetting factor $\lambda = 0.9$, and employ a generalized harmonic learning rate formula \citep{powell2011approximate}, expressed as $\alpha_t = \frac{40}{5000 + t - 1}$ at $t^{th}$ period of the algorithm. The $\epsilon$-greedy approach is controlled by $\epsilon_t = 0.999983^t$. At the end of training, the assigned weights are tested on $60K$ simulation periods, with the first $1K$ periods excluded as a warm-up period.

For LCRL, we implement $T=1$, as slight increase in the value causes explosion in solution time. We implement a two-stage decision tree with $|\Omega| = 20$. For its training, the same weights of CRL is implemented. For simulation, since the solution time increases substantially, we test on $3K$ simulation periods, with the first $20$ periods excluded as a warm-up.

For the $(s,S)$-based heuristic, simulations on a pair of $(s,S)$ values are done over 1M training periods in the for loop, divided into $50$ episodes. Each episode excludes the first 1K periods as a warm-up period. Additionally, the simulation in the while loop is done over 2M training periods in the for loop, divided into $50$ episodes, with the initial 1K periods of each excluded for warm-up. We set the parameters $\xi_0 \coloneqq 0.01 q$ and $m = 1.1$.

For Power-of-two (PO2) based solution, we set $\tau = 4$. The simulation of the cyclic policy is done the same of $(s,S)$, with 2M periods, $50$ episodes, and a warm-up of 1K periods per episode. Both $(s,S)$ and PO2 methods utilize significantly more periods in simulation since their simulation process is comparatively simpler than of (L)CRL. This is because learning algorithms solve a mixed-integer quadratic programming (MIQP) model within a MIP feasible region whenever rand $>\epsilon$. This is more time consuming due to the complexity of the underlying optimization problem.

\end{APPENDICES} 

\end{document}